\documentclass[12pt,leqno]{article}

\usepackage{amsmath,amssymb,latexsym,theorem,epsfig}
\usepackage[all]{xy}
\CompileMatrices

\newcommand{\bbz}{{\mathbb Z}}
\newcommand{\bbf}{{\mathbb F}}

\newcommand{\cp}[1]{{\mathbb P}^{#1}}
\newcommand{\op}[1]{\operatorname{#1}}
\newcommand{\MW}{{\mathbb M}{\mathbb W}}
\newcommand{\mw}[1]{[#1]}
\newcommand{\Pic}{\op{Pic}}
\newcommand{\FM}{\boldsymbol{F}{\boldsymbol{M}}}
\newcommand{\fm}{\boldsymbol{f}\boldsymbol{m}}
\newcommand{\T}{\boldsymbol{T}}
\newcommand{\ct}{\boldsymbol{t}}
\newcommand{\D}{\boldsymbol{D}}

\newcommand{\rhom}{R^{\bullet}{\mathcal H}om}
\newcommand{\cL}{{\mathcal L}}
\newcommand{\cV}{{\mathcal V}}
\newcommand{\cW}{{\mathcal W}}

\newtheorem{theo}{Theorem}[section]
\newtheorem{lem}[theo]{Lemma}
\newtheorem{cor}[theo]{Corollary}
\newtheorem{prop}[theo]{Proposition}

{\theorembodyfont{\rmfamily} \newtheorem{rem}[theo]{Remark}}
{\theorembodyfont{\rmfamily} }

\numberwithin{equation}{section}

\begin{document}

\title{   \vspace*{-5em} 
\  \hfill {\normalsize CERN-TH/2000-203, UPR-894T, RU-00-5B} \\[0.5em]
\ \\
{\sf\LARGE Spectral involutions on rational
elliptic surfaces}}
\author{Ron Donagi$^1$, Burt A.~Ovrut$^2$, Tony Pantev$^1$ 
      and Daniel Waldram$^{3}$ \\ [0.5em]
      \ \\
   {\normalsize $^1$Department of Mathematics, 
      University of Pennsylvania} \\[-0.2em]
      {\normalsize Philadelphia, PA 19104--6395, USA}\\
   {\normalsize $^2$Department of Physics, 
      University of Pennsylvania} \\[-0.2em]
      {\normalsize Philadelphia, PA 19104--6396, USA}\\
    {\normalsize $^3$Theory Division, CERN CH-1211, Geneva 23, 
Switzerland, and}  \\[-0.2em]
{\normalsize Department of Physics, The Rockfeller University} \\[-0.2em]
{\normalsize New York, NY 10021}
}
\date{}
\maketitle

\begin{abstract} In this paper we describe a four dimensional
family of special rational elliptic surfaces admitting an involution with
isolated fixed points.  For each surface in this family we calculate
explicitly the action of a spectral version of the involution (namely
of its Fourier-Mukai conjugate) on global line bundles and on spectral
data. The 
calculation is carried out both on the level of cohomology and in the
derived category. We find that the spectral
involution behaves like a
fairly simple affine transformation away from the union of those
fiber components which do not intersect the zero section. These
results are the key ingredient in the construction of Standard-Model bundles
in \cite{dopw-ii}.

\smallskip

\noindent
{\bf MSC 2000:} 14D20, 14D21, 14J60
\end{abstract}

\section{Introduction} \label{s-intro}

Let $Z \to S$ be an elliptic fibration on a smooth variety $Z$, i.e. a
flat morphism whose generic fiber is a curve of genus one, and which
has a section $S \to Z$. The choice of such a section defines a
Poincare sheaf ${\mathcal P}$ on $Z\times_{S} Z$. 
The corresponding Fourier-Mukai transform $\FM : D^{b}(Z) \to
D^{b}(Z)$  is then an
autoequivalence of the derived category $D^{b}(Z)$ of complexes of
coherent sheaves on $Z$. It sets up an equivalence between
$SL(r,{\mathbb C})$-bundles on $Z$ and spectral data consisting of
line bundles 
(and their degenerations) on spectral covers $C \subset Z$ which are
of degree $r$ over $S$. This equivalence has been used extensively to
construct vector bundles on elliptic fibrations and to study their
moduli \cite{fmw,donagi,bjps}. 

For many applications it is important to remove the requirement of the
existence of a section, i.e. to allow genus one fibrations. This could
be done in two ways.

The `spectrum' of a degree zero semistable rank $r$ bundle on a genus
one curve $E$ consists of $r$ points in the Jacobian
$\op{Pic}^{0}(E)$, rather than in $E = \op{Pic}^{1}(E)$ itself. So one
approach is to consider spectral covers $C$ contained in the relative
Jacobian $\op{Pic}^{0}(Z/S)$. But the spectral data in this case no
longer involves a line bundle on $C$; instead, it lives in a certain
non-trivial gerbe, or twisted form of $\op{Pic}(C)$. So the essential
problem becomes the analysis of this gerbe. 

The second approach is to find an elliptic fibration $\pi : X \to B$
together with a group $G$ acting compatibly on $X$ and $B$ (but not
preserving the section of $\pi$) such that the action on $X$ is fixed
point free and the quotient is the original $Z \to S$. One can then
use  the Fourier-Mukai transform to construct vector bundles on
$X$. The problem becomes the determination of conditions for such a
bundle on $X$ to be $G$-equivariant, hence to descend to
$Z$. Equivalently we need to know the action of each $g \in G$ on spectral
data. This is the restriction of the action on $D^{b}(X)$ of the
Fourier-Mukai conjugate $\FM^{-1}\circ g^{*} \circ \FM$ of
$g^{*}$. This will be referred to as the {\em spectral action} of
$g$. Unfortunately, the spectral action can be quite complicated: both global
vector bundles on $X$ and sheaves supported on $C$ can go to complexes
on $X$ of amplitude greater than one. 

In this paper, 
we work out such a spectral action in one
class of examples consisting of special rational elliptic surfaces. 
In the second part \cite{dopw-ii} of this paper we
use this analysis to construct special bundles  on certain non-simply
connected smooth Calabi-Yau  threefolds. 
These special bundles in turn are the
main ingredient for the construction of Heterotic M-theory vacua
having the Standard Model symmetry group $SU(3)\times SU(2)\times
U(1)$ and three generations of quarks and leptons.  
The physical significance  of such vacua
is explained in \cite{usnew} and was the original motivation of this work.

Here is an outline of the paper.
We begin in section~\ref{ss-dp9-general} with a review of the basic
properties of rational elliptic surfaces.  Within the eight
dimensional moduli space of all rational elliptic surfaces we focus
attention on a five dimensional family of rational elliptic surfaces
admitting a particular involution $\tau$, and then we restrict further
to a four dimensional family of surfaces with reducible fibers. This
seems to be the simplest family of surfaces for which one needs the
full force of Theorem~\ref{prop-T}: for general surfaces in the five
dimensional family, the spectral involution $\T := \FM^{-1}\circ
\tau^{*} \circ \FM$ of $\tau$ takes line bundles to line
bundles, while in the four dimensional subfamily it is possible for
$\T$ to take a line bundle to a complex which can not be represented by
any single sheaf. We study the five dimensional family in
section~\ref{ss-dp9-special} and the four dimensional subfamily in
section~\ref{s-4dim-family}. This section concludes, in
subsection~\ref{sss-synthetic}, with a synthetic construction of the
surfaces in the four dimensional subfamily. This construction maybe
less motivated than the original a priori analysis we use, but it is more
concise and we hope it will make the exposition more accessible. 

In the remainder of the paper we work out the actions of $\tau$,
$\FM$, $\T$, first at the level of cohomology in
sections~\ref{ss-dp9-cohoaction} and \ref{s-cohoFM}, and then on the
derived category in  section~\ref{ss-dp9-bundleaction}. The main
result is Theorem~\ref{prop-T}, which says that $\T$ behaves like a
fairly simple affine transformation {\em away from the union of those
fiber components which do not intersect the zero section}. A corollary
is that for spectral curves which do not intersect the extra vertical
components, all the complications disappear. This fact together with
the cohomological formulas from sections \ref{ss-dp9-cohoaction} and
\ref{s-cohoFM} will be used in \cite{dopw-ii} to build
invariant vector bundles on a family of Calabi-Yau threefolds
constructed from the rational elliptic surfaces in our four
dimensional subfamily.

\

\bigskip

{\bf Acknowledgements:} We would like to thank Ed Witten, Dima Orlov,
and Richard Thomas for valuable conversations on the subject of
this work. 

R.~Donagi is supported in part by an NSF grant DMS-9802456 as well as a
UPenn Research Foundation Grant. 
B.~A.~Ovrut is supported in part by a Senior Alexander von Humboldt
Award, by the DOE under contract No. DE-AC02-76-ER-03071 and by a
University of Pennsylvania Research Foundation Grant. 
T.~Pantev is supported in part by an NSF grant DMS-9800790 and by an
Alfred P. Sloan Research Fellowship. 
D.~Waldram would like to thank Enrico Fermi Institute at The
University of Chicago and the Physics Department of The Rockefeller
University for hospitality during the completion of this work.

\tableofcontents

\section{Rational elliptic surfaces} \label{ss-dp9-general}

A rational elliptic surface is a rational surface $B$ which admits an
elliptic fibration  $\beta : B \to \cp{1}$. It can be described as
the blow-up of the plane $\cp{2}$ at nine points $A_{1}, \ldots, A_{9}$
which are the base
points of a pencil $\{ f_{t} \}_{t \in \cp{1}}$ of cubics. The map
$\beta$ is recovered as the anticanonical map of $B$ and the proper
transform of $f_{t}$ is $\beta^{-1}(t)$. 

In particular the topological Euler characteristic of $B$ is $\chi(B)
= \chi(\cp{2}) + 9 = 12$. For a generic $B$ the map $\beta$ has twelve
distinct singular fibers each of which has a single node. For future
use we denote by $B^{\#} \subset B$ the open set of regular points of
$\beta$ and we set  $\beta^{\#} := \beta_{|B^{\#}}$.

Under mild general position requirements \cite{lnm777}
each subset of eight of these points determines the pencil of cubics
and hence the ninth point. In particular we see that the rational
elliptic surfaces depend on $2\cdot 8 - \dim {\mathbb P}GL(3,{\mathbb
C}) = 8$ parameters.

Let $e_{1}, \ldots, e_{9}$ be the exceptional divisors in $B$
corresponding to the $A_{i}$'s. Let $\ell$ be the preimage of the class of
a line in $\cp{2}$ and let $f := \beta^{*}{\mathcal O}_{\cp{1}}(1)$.
Note that
\[
f = - K_{B} = 3\ell - \sum_{i=1}^{9}e_{i}
\]
and that $\ell, e_{1}, \ldots, e_{9}$ form a basis of
$H^{2}(B,\bbz)$.

The curves $e_{1}, e_{2}, \ldots, e_{9}$ are sections  of the map
$\beta : B \to \cp{1}$. Choosing a section $e : \cp{1} \to B$
determines a group law on the fibers of $\beta^{\#}$. The inversion
for this group law is an involution on $B^{\#}$ which for a general
$B$ extends to a well defined involution $(-1)_{B,e} :  B \to B$.
When $B$ or $e$
are understood from the context we will just write $(-1)_{B}$ or
$(-1)$. The involution $(-1)_{B,e}$ fixes the section $e$ as
well as  a tri-section of $\beta$ which parameterizes the non-trivial
points of order two. The quotient $W_{\beta}/(-1)_{B,e}$ is a smooth
rational surface which is ruled over the base $\cp{1}$. For a general
$B$ this quotient is the Hirzebruch surface $\bbf_{2}$ and the image of
$e$  is the exceptional section of $\bbf_{2}$. This
gives yet another realization of $B$ as a branched double cover of $\bbf_{2}$.

A convenient way to describe the involution $(-1)_{B,e}$ is
through the Weierstrass model $w : W_{\beta} \to \cp{1}$ of  
$
\xymatrix@1{B
\ar[r]^-{\beta} & {\mathbb P}^{1} \ar@/^0.5pc/[l]^-{e}}.
$

\medskip

The model $W_{\beta}$ is described explicitly as follows. By relative
duality $R^{1}\beta_{*}{\mathcal O}_{B} \cong {\mathcal
O}_{\cp{1}}(-1)$. This implies that $\beta_{*}{\mathcal
O}_{B}(3e) = ({\mathcal O}_{\cp{1}} \oplus {\mathcal O}_{\cp{1}}(2)
\oplus {\mathcal O}_{\cp{1}}(3))^{\vee}$. Let
\[
p : P := {\mathbb P}({\mathcal
O}_{\cp{1}} \oplus {\mathcal O}_{\cp{1}}(2) \oplus {\mathcal
O}_{\cp{1}}(3))  \to {\mathbb P}^{1}.
\]
be the natural projection. The linear system ${\mathcal O}_{B}(3e)$
defines a map $\nu : B \to P$ compatible with the projections.
The Weierstrass model $W_{\beta}$ is defined to be the image of
this map. It is given explicitly by an equation
\[
y^{2}z = x^{3} + (p^{*}g_{2})xz^{2} +
(p^{*}g_{3})z^{3}
\] 
where $g_{2} \in H^{0}({\mathcal O}_{{\mathbb P}^{1}}(4))$ and $g_{3}
\in H^{0}({\mathcal O}_{{\mathbb P}^{1}}(6))$ and $x$, $y$ and $z$ are the
natural sections of ${\mathcal O}_{P}(1)\otimes p^{*}{\mathcal
O}_{{\mathbb P}^{1}}(2)$, ${\mathcal O}_{P}(1)\otimes p^{*}{\mathcal
O}_{{\mathbb P}^{1}}(3)$ and ${\mathcal O}_{P}(1)$ respectively.

In terms of $W_{\beta}$ the section $e$ is given by $x = z = 0$ and
the involution $(-1)_{B,e}$ sends $y$ to $-y$. The tri-section of fixed
points of $(-1)_{B,e}$ is given by $y = 0$.

\bigskip

The Mordell-Weil group $\MW = \MW(B,e)$ is the group of sections of
$\beta$. As a set $\MW$ is the collection of all sections of $\beta : B \to
\cp{1}$ or equivalently all sections of $\beta^{\#} : B^{\#}
\to \cp{1}$. The group law on $\MW$ is induced from the addition law on
the group scheme $\beta^{\#} : B^{\#}
\to \cp{1}$ and so $e$ corresponds to the neutral element in
$\MW(B,e)$. For a section $\xi \subset B$  we will put
$\mw{\xi}$ for the corresponding element of $\MW$. Note that the natural map 
\[
c_{1} : \MW(B,e) \to \Pic (B), \qquad \mw{\xi} \mapsto {\mathcal O}_{B}(\xi).
\] 
is not a group homomorphism. When written out in
coordinates, it involves both a  linear part and a quadratic
term (see e.g. \cite{manin}). However, when $B$ is smooth the map
$c_{1}$ induces a linear map to a quotient of $\op{Pic}(B)$ which
describes $\MW(B,e)$ completely. Indeed, let $B$ be smooth and let 
${\mathcal T} \subset \op{Pic}(B)$ be the sublattice generated by $e$
and all the components of the fibers of $\beta$. Then $c_{1}$ induces a map
\[
\bar{c}_{1} : \MW(B,e) \to \Pic (B)/{\mathcal T}, 
\qquad \mw{\xi} \mapsto ({\mathcal O}_{B}(\xi)\!\!\!\mod{\mathcal T}) 
\]
which is a linear isomorphism \cite[Theorem~1.3]{shioda}
\

\bigskip

There is a natural group homomorphism $t : \MW \to \op{BirAut}(B)$
assigning to each  section $\xi \in \MW$ the birational automorphism
$t_{\xi} : B \dashrightarrow B$, which on the open set $B^{\#}$ 
is just translation by $\xi$
with respect to the group law determined by $e$. When $\beta : B \to
\cp{1}$ is relatively minimal the map $t_{\xi}$
extends canonically to a biregular automorphism of $B$
\cite[Theorem~2.9]{kodaira-casIII}.

\section{Special rational elliptic surfaces} \label{ss-dp9-special}

In the second part of this paper \cite{dopw-ii} we will work with
Calabi-Yau threefolds $X$ which are elliptically fibered over a
rational elliptic surface $B$. 
Any involution $\tau_X$ on an elliptic CY $\pi:X \to B$ commuting with
$\pi$ induces 
(either the identity or) an involution $\tau_B$ on the base $B$. 
In order for
$\tau_X$ to act freely on $X$ we need the fixed points of $\tau_{B}$ to
be disjoint from the discriminant of $\pi$. If $B$ is a rational
elliptic surface, then the discriminant of $\pi$ is a section in
$K_{B}^{-12} = {\mathcal O}_{B}(12f)$ and so $(-1)_{B}$ will not do.  
We want to describe some special rational elliptic surfaces which 
admit additional involutions. Within the $8$ dimensional family of
rational elliptic surfaces we describe first a $5$ dimensional family of
surfaces which  admit an involution $\alpha_{B}$. The fixed locus of
$\alpha_{B}$ has the right properties but it turns out that $\alpha_{B}$
does not lift to a free involution on $X$. However, one can easily
show that each $\alpha_{B}$ can be corrected by a translation 
$t_{\zeta}$  (for a special type   
of section $\zeta$) to obtain  an additional involution $\tau_{B}$
which does the job.  Unfortunately the general member of the $5$
dimensional family leads to a Calabi-Yau manifold which does not admit
any bundles satisfying all the constraints required by the Standard
Model of particle physics (see \cite{dopw-ii}). 
We therefore specialize further to a $4$ dimensional
family of surfaces for which the extra involution $\tau_{B}$ can be
constructed in an explicit geometric way. This provides some extra
freedom which enables us to carry out the construction.  The involution
$\alpha_{B}$ fixes one fiber of $\beta$ and four points in another
fiber. The involution $\tau_{B}$ fixes only four points in one fiber.
A special feature of the $4$ dimensional family is that it consists of
$B$'s for which $\beta$ has at least two $I_{2}$ fibers. This translates
into a special position requirement on the nine points in $\cp{2}$. Another
special feature of the $4$ dimensional family is seen
in the double cover realization of $B$ where the quotient $B/(-1)$
becomes $\bbf_{0} = \cp{1}\times \cp{1}$ instead of $\bbf_{2}$. 

We thank Chad Schoen for pointing out that essentially the same
surfaces and threefolds were constructed in section 9 of
\cite{schoenCY}. The explicit example he gives there for what he calls
the ``$m=2$ case'' happens to exactly coincide with our
four-dimensional family of rational elliptic surfaces. With small
modifications, his construction could have given our full
five-dimensional family as well. Schoen's construction technique is
rather different than ours. He constructs the equivalent of our
rational elliptic surface $B$ and involution $\alpha_B$ directly (the
surfaces we call $B, B/ {\alpha_B}$ are called $Y_0, T_0$ in
\cite{schoenCY}); then he invokes a general result of Ogg and
Shafarevich for the existence of a logarithmic transform $Y$ with
quotient $T$; and finally, results from classification theory are used
to deduce existence of an abstract isomorphism of $Y$ with $Y_0$ such
that his $T$ becomes our $B/ {\tau_B}$.

In the next several sections we will describe the structure
of the rational elliptic surfaces that admit additional
involutions. This rather extensive geometric analysis is ultimately
distilled into a fairly simple synthetic
construction of our surfaces which is explained in
section~\ref{sss-synthetic}.  The impatient reader who is interested
only in the end result of the construction and wants to avoid the tedious
geometric details is advised to skip directly to section \ref{sss-synthetic}.

\subsection{Types of involutions on a rational elliptic surfaces}
\label{sss-dp9-involutions}
Consider a smooth rational elliptic surface $\xymatrix@1{B
\ar[r]^-{\beta} & {\mathbb P}^{1} \ar@/^0.5pc/[l]^-{e}}$
with a fixed section.  For any automorphism
$\tau_{B}$ of $B$ we have $\tau_{B}^{*}K_{B} \cong
K_{B}$. Since $K_{B}^{-1} = \beta^{*}{\mathcal O}_{{\mathbb P}^{1}}(1)$
this implies that $\tau_{B}$  induces an automorphism 
$\tau_{{\mathbb P}^{1}} : {\mathbb P}^{1}
\to {\mathbb P}^{1}$. If $\tau_{B}$ is an involution we have
two possibilities: either $\tau_{{\mathbb P}^{1}} = \op{id}_{{\mathbb P}^{1}}$
or $\tau_{{\mathbb P}^{1}}$ is an involution of ${\mathbb P}^{1}$. 

Both of these cases occur and lead to Calabi-Yau manifolds with freely
acting involutions. For concreteness here we only treat the case when
$\tau_{\cp{1}}$ is an involution. The case $\tau_{{\mathbb P}^{1}} =
\op{id}_{{\mathbb P}^{1}}$ can be analyzed easily in a similar fashion.

If $\tau_{{\mathbb P}^{1}}$ is an involution, then  $\tau_{{\mathbb
P}^{1}}$ will have two fixed points on
${\mathbb P}^{1}$ which we will denote by $0, \infty \in {\mathbb
P}^{1}$. Note that every involution on $\cp{1}$ is uniquely determined by its
fixed points and so specifying $\tau_{\cp{1}}$ is equivalent to
specifying the points $0, \infty \in \cp{1}$. Next we classify the
types of involutions on $B$ that lift a given involution $\tau_{\cp{1}}$.

\begin{lem} \label{lem-inv-types} Let $\beta : B \to \cp{1}$ be a
rational elliptic surface and let $\tau_{\cp{1}} : \cp{1} \to \cp{1}$
be a fixed involution. There is a canonical bijection
\[
\left\{  
\begin{minipage}[c]{2.5in}
Involutions $\tau_{B} : B \to B$, satisfying $\tau_{\cp{1}}\circ \beta
= \beta\circ \tau_{B}$.
\end{minipage}
\right\}
\leftrightarrow
\left\{ 
\begin{minipage}[c]{2.5in}
Pairs $(\alpha_{B},\zeta)$ consisting of: 
\begin{itemize}
\item An involution $\alpha_{B} : B \to B$, satisfying
$\tau_{\cp{1}}\circ \beta = \beta\circ \alpha_{B}$ which leaves the zero
section invariant, i.e.  $\alpha_{B}(e) = e$.
\item A section $\zeta$ of $\beta$ satisfying $\alpha_{B}(\zeta) =
(-1)_{B}(\zeta)$. 
\end{itemize}
\end{minipage}
\right\}
\]
\end{lem}
{\bf Proof.} Let $\tau_{B} : B \to B$ be such that $\tau_{\cp{1}}\circ \beta
= \beta\circ \tau_{B}$.
Put $\zeta = \tau_{B}(e)$ for the image of the zero section
under $\tau_{B}$ and let $\alpha_{B} = t_{-\zeta}\circ \tau_{B}$. 

Then $\alpha_{B}$ is an automorphism of $B$ which induces
$\tau_{{\mathbb P}^{1}}$ on ${\mathbb P}^{1}$ and preserves the zero
section $e \subset B$. So $\alpha_{B}^{2} : B \to B$ will be an
automorphism of $B$ which acts trivially on ${\mathbb P}^{1}$. But
\[
t_{-\zeta}\circ \tau_{B} =
\tau_{B}\circ t_{- \tau_{B/\cp{1}}^{*-1}(\zeta)}
\] 
where $\tau_{B/\cp{1}}^{*} :
\op{Pic}^{0}(B/{\mathbb P}^{1}) \to \op{Pic}^{0}(B/{\mathbb P}^{1})$
is the involution  on the relative Picard scheme induced from $\tau_{B}$. 
In particular we have that $\alpha_{B}^{2}$ must be a translation by a
section. Indeed we  have
\begin{equation} \label{eq-alpha2}
\alpha_{B}^{2} = t_{-\zeta}\circ \tau_{B} \circ \tau_{B}
\circ t_{- \tau_{B/\cp{1}}^{*-1}(\zeta)} = t_{- \zeta -
\tau_{B/\cp{1}}^{*-1}(\zeta)}. 
\end{equation}
Combined with the fact that $\alpha_{B}^{2}$ preserves $e$ 
\eqref{eq-alpha2} implies that $\alpha_{B}^{2}
= \op{id}_{B}$. 
On the other hand, if we use the zero section $e$ to identify 
$\op{Pic}^{0}(B/{\mathbb P}^{1}) \to \cp{1}$ with $\beta^{\#} :
B^{\#} \to \cp{1}$, 
then  $\tau_{B/\cp{1}}^{*} = \alpha_{B}$. Indeed, 
let $\xi \in \op{Pic}^{0}(B/{\mathbb P}^{1})$
and  let $x \in {\mathbb P}^{1}$ be the projection of the point $\xi$.  
Let $f_{x} \subset B$  be the fiber of $\beta$ over $x$. Denote by 
$m_{\xi} \in f_{x}$ the unique smooth point in $f_{x}$ for which
${\mathcal O}_{f_{x}}(m_{\xi}) = \xi\otimes {\mathcal
O}_{f_{x}}(e(x))$. Then by definition $\tau_{B}^{*}(\xi)$ is a line
bundle of degree zero on $f_{x}$ such that 
\[
{\mathcal O}_{f_{x}}(\tau_{B}(m_{\xi})) = \tau_{B}\xi\otimes {\mathcal
O}_{f_{x}}(\tau_{B}(e(x)))  = \tau_{B}\xi\otimes {\mathcal
O}_{f_{x}}(\zeta(x)).
\] 
In other words 
under the identification of $\op{Pic}^{0}(f_{x})$ with the smooth 
locus of $f_{x}$ via $e(x)$ the line bundle
$\tau_{B}^{*}\xi \to f_{x}$  corresponds to the unique point $p_{\xi}$
of $f_{x}$ such that
\[
{\mathcal O}_{f_{x}}(p_{\xi}) = 
{\mathcal O}_{f_{x}}(\tau_{B}(m_{\xi}))\otimes {\mathcal
O}_{f_{x}}(e(x) - \zeta(x)).
\] 
But the right hand side of this identity equals 
${\mathcal O}_{f_{x}}(\alpha_{B}(m_{\xi}))$
by definition and so $p_{\xi} = \alpha_{B}(m_{\xi})$.

Combined with the identity \eqref{eq-alpha2} and the fact that $t :
\MW(B) \to \op{Aut}(B)$ is injective this yields
\[
\alpha_{B}(\zeta) = (-1)_{B}(\zeta).
\]
Conversely, given a pair $(\alpha_{B},\zeta)$ we set $\tau_{B} =
t_{\zeta}\circ \alpha_{B}$. Clearly $\tau_{B}$ is an automorphism of
$B$ which induces $\tau_{\cp{1}}$ on $\cp{1}$. Furthermore we
calculate $\tau_{B}^{2} = t_{\zeta}\circ \alpha_{B}\circ t_{\zeta}
\circ \alpha_{B} = t_{\zeta}\circ \alpha_{B}\circ \alpha_{B}\circ
t_{-\zeta} = \op{id}_{B}$. The lemma is proven. \hfill $\Box$

\bigskip

\noindent
The above lemma implies that  in order to understand all
involutions $\tau_{B}$ 
it suffices to understand all pairs $(\alpha_{B},\zeta)$. Since the
involutions $\alpha_{B}$ stabilize $e$ it follows that $\alpha_{B}$ 
will have to necessarily act on
the Weierstrass model of $B$. In the next section we analyze this
action in more detail.

\subsection{The Weierstrass model of $B$}
\label{subsec-weierstrass-case2}  
Let as before  
$\tau_{{\mathbb P}^{1}} :
{\mathbb P}^{1} \to {\mathbb P}^{1}$ be an involution and let
$(t_{0}:t_{1})$ be homogeneous coordinates on ${\mathbb P}^{1}$ such
that $\tau_{{\mathbb P}^{1}}((t_{0}:t_{1})) = (t_{0}:-t_{1})$ and $0 =
(1:0)$ and $\infty = (0:1)$. Since 
$t_{0}$ and $t_{1}$ are a basis of $H^{0}({\mathbb P}^{1}, {\mathcal
O}_{{\mathbb P}^{1}}(1))$ and since ${\mathcal
O}_{{\mathbb P}^{1}}(1)$ is generated by global sections we can lift
the action of $\tau_{{\mathbb P}^{1}}$ to ${\mathcal O}_{{\mathbb
P}^{1}}(1)$. For concreteness choose the lift $t_{0} \mapsto t_{0}$,
$t_{1} \mapsto - t_{1}$. Since $H^{0}({\mathbb P}^{1},{\mathcal
O}_{{\mathbb P}^{1}}(k)) = S^{k} H^{0}({\mathbb P}^{1}, {\mathcal
O}_{{\mathbb P}^{1}}(1))$ we get
a lift of the action of $\tau_{{\mathbb P}^{1}}$ to the line bundles 
${\mathcal O}_{{\mathbb P}^{1}}(k)$ for all $k$. We will call this
action {\em the standard action} of $\tau_{{\mathbb P}^{1}}$ on
${\mathcal O}_{{\mathbb P}^{1}}(k)$. Via the standard action the
involution $\tau_{{\mathbb P}^{1}}$ acts also on the vector bundle
${\mathcal O}_{{\mathbb P}^{1}}\oplus {\mathcal O}_{{\mathbb
P}^{1}}(2) \oplus {\mathcal O}_{{\mathbb P}^{1}}(3)$ and hence we get
a standard lift  $\tau_{P} : P \to P$ of $\tau_{{\mathbb P}^{1}}$
satisfying $\tau_{P}^{*}{\mathcal O}_{P}(1) \cong {\mathcal O}_{P}(1)$.

Assume that we are given an involution $\alpha_{B} : B \to B$ which
induces $\tau_{{\mathbb P}^{1}}$ on ${\mathbb P}^{1}$ and preserves
the section $e$. We have the following 

\begin{lem} \label{lem-alphaW} 
\begin{itemize}
\item[{\em (i)}]  There exists a unique involution $\alpha_{W_{\beta}}
: W_{\beta} \to W_{\beta}$ such that the natural map $\nu : B \to
W_{\beta}$ satisfies $\alpha_{W_{\beta}}\circ \nu =
\nu\circ\alpha_{B}$.  
\item[{\em (ii)}] Let $W \subset P$ be a Weierstrass rational elliptic
surface. Then the involution $\tau_{\cp{1}}$ lifts to an involution on
$W$ which preserves the zero section if and only if $\tau_{P}(W) = W$.
\item[{\em (iii)}] If $w : W_{\beta} \to \cp{1}$ is not isotrivial,
then $\alpha_{W_{\beta}}$ is either $\tau_{P|W_{\beta}}$ or
$\tau_{P|W_{\beta}}\circ (-1)_{W_{\beta}}$.
\end{itemize}
\end{lem}
{\bf Proof.} Since $\alpha_{B}^{*}({\mathcal
O}_{B}(e)) \cong {\mathcal O}_{B}(e)$, there
exists an
involution on the total space of the bundle ${\mathcal O}_{B}(e)$ which
acts linearly on the fibers and induces the involution $\alpha_{B}$ on
$B$. Indeed - the square $\gamma\circ
\alpha_{B}^{*}\gamma$ of the
isomorphism $\gamma : \alpha_{B}^{*}({\mathcal
O}_{B}(e)) \widetilde{\to} {\mathcal O}_{B}(e)$ is a bundle automorphism of
${\mathcal O}_{B}(e)$ (acting trivially on the base) and so is given
by multiplication by some non-zero complex number $\lambda \in
{\mathbb C}$. Rescaling the isomorphism $\gamma$ by
$\sqrt{\lambda^{-1}}$ then gives the desired lift.

In this way the involution $\alpha_{B}$ induces an involution on
${\mathcal O}_{e}(-e) = {\mathcal O}_{{\mathbb P}^{1}}(1)$ which lifts
the action of $\tau_{{\mathbb P}^{1}}$. Let us normalize  the lift of
$\alpha_{B}$ to ${\mathcal O}_{B}(e)$ so that the induced action on 
${\mathcal O}_{e}(-e) = {\mathcal O}_{{\mathbb P}^{1}}(1)$ coincides
with the standard action of $\tau_{{\mathbb P}^{1}}$. Thus the
Weierstrass model $W_{\beta} \subset P$ must be stable under the
corresponding 
$\tau_{P}$ and the restriction of $\tau_{P}$ to $W_{\beta}$ is an
involution that preserves the zero section of $w$ and induces
$\tau_{{\mathbb P}^{1}}$ on the base. By construction
$\tau_{P|W_{\beta}}$ coincides with the involution induced from
$\alpha_{B}$ up to a composition with $(-1)_{W_{\beta}}$. 
This finishes the proof of the lemma. \hfill $\Box$

\bigskip

\noindent
We are now ready to construct the Weierstrass models of all surfaces
$B$ that admit an involution $\alpha_{B}$.  Similarly to the proof of
Lemma~\ref{lem-alphaW}, the fact that 
$\tau_{P}^{*}{\mathcal O}_{P}(1) \cong {\mathcal O}_{P}(1)$ implies
that the action of $\tau_{P}$ can be lifted to an action on ${\mathcal
O}_{P}(1)$. Since there are two possible such lifts and they differ by
multiplication by $\pm 1 \in {\mathbb C}^{\times}$ 
we can use the identification ${\mathcal
O}_{P}(1)_{|B} = {\mathcal O}_{B}(3e)$ to choose the unique lift that
will induce the standard action of $\tau_{{\mathbb P}^{1}}$ on
${\mathcal O}_{{\mathbb P}^{1}}(3) = {\mathcal O}_{e}(-3e)$. With
these choices we define an action 
\[
\tau_{P}^{*} : H^{0}(P,{\mathcal O}_{P}(r)\otimes p^{*}{\mathcal
O}_{P}(s)) \to H^{0}(P,{\mathcal O}_{P}(r)\otimes p^{*}{\mathcal
O}_{P}(s))
\]
of $\tau_{P}$ on the global sections of any line bundle on $P$. Note
that by construction we have $\tau_{P}^{*}x = x$, $\tau_{P}^{*}y = y$
and $\tau_{P}^{*}z = z$.

Consider the general equation of the Weierstrass model $W_{\beta}$ of $B$:
\begin{equation} \label{eq-weierstrass-case2}
y^{2}z  = x^{3} + (p^{*}g_{2})xz^{2} + (p^{*}g_{3})z^{3}.
\end{equation}
Here $g_{2} \in H^{0}({\mathcal O}_{{\mathbb P}^{1}}(4))$ and $g_{3}
\in H^{0}({\mathcal O}_{{\mathbb P}^{1}}(6))$. The fact $W_{\beta}  
\subset P$ is stable under $\tau_{P}$ implies that the image of
the Weierstrass equation (\ref{eq-weierstrass-case2}) 
under $\tau_{P}^{*}$ must be a
proportional Weierstrass equation. In particular we ought to have 
$\tau_{{\mathbb P}^{1}}^{*}g_{2} = g_{2}$ and $\tau_{{\mathbb
P}^{1}}^{*}g_{3} = g_{3}$.

Conversely, for any $g_{2} \in H^{0}({\mathcal O}_{{\mathbb
P}^{1}}(4))$ and $g_{3} \in H^{0}({\mathcal O}_{{\mathbb P}^{1}}(6))$
which are invariant  for the standard action of $\tau_{{\mathbb
P}^{1}}$ it follows that $\tau_{P}$ will preserve the Weierstrass
surface $W$ given by the equation (\ref{eq-weierstrass-case2}). Note that
for a generic choice of $g_{2}$ and $g_{3}$ the surface $W$ will
be smooth and so $B = W$, $\alpha_{B} = \tau_{P|W}$. When $W$ is
singular, the surface $B$ is the minimal resolution of singularities
of $W$ and hence $\alpha_{W} = \tau_{P|W}$ determines uniquely 
$\alpha_{B}$ by the universal property of the minimal resolution.

\bigskip

Next we describe the fixed locus of $\alpha_{B}$. Note that since
$\alpha_{B}$ induces $\tau_{{\mathbb P}^{1}}$ on ${\mathbb P}^{1}$ the
fixed points of $\alpha_{B}$ will necessarily sit over the two fixed
points of $\tau_{{\mathbb P}^{1}}$. So in order to understand the fixed
locus of $\alpha_{B}$ it suffices to understand the action of
$\alpha_{B}$ on the two $\alpha_{B}$-stable fibers of $\beta$ -
namely $f_{0} = \beta^{-1}(0)$ and $f_{\infty} = \beta^{-1}(\infty)$.

\begin{lem} \label{lem-alphaB-fixed} Let $\alpha_{B}$ be the involution
on $B$ induced from $\tau_{P|W_{\beta}}$ (with the above
normalizations). Then  $\alpha_{B}$ fixes $f_{0}$ pointwise and has
four isolated fixed points on $f_{\infty}$, namely the points of order two.
\end{lem}
{\bf Proof.} The curve $f_{0}$ is a smooth cubic in the projective plane 
\[
P_{0} = {\mathbb P}({\mathcal O}_{0}\oplus {\mathcal O}(2)_{0}
\oplus {\mathcal O}(3)_{0}),
\] 
Where ${\mathcal O}(k)_{0}$ denotes the fiber of the line bundle
${\mathcal O}_{{\mathbb P}^{1}}(k)$ at the point $0 \in {\mathbb P}^{1}$.
Note that $1$, $t_{0}(0)^{2}$ and $t_{0}(0)^{3}$ span the lines
${\mathcal O}_{0}$, ${\mathcal O}(2)_{0}$ and ${\mathcal O}(3)_{0}$
respectively and so $\tau_{{\mathbb P}^{1}}$ acts trivially on those
lines via its standard action. So if we identify those lines with
${\mathbb C}$ via the basis $1$, $t_{0}(0)^{2}$ and $t_{0}(0)^{3}$, then
$X_{0} := x_{|P_{0}}$, $Y_{0} := y_{|P_{0}}$ and $Z_{0} := 
z_{|P_{0}}$ become identified with
sections of the line bundle ${\mathcal O}_{P_{0}}(1)$ and can be used
as homogeneous coordinates on $P_{0}$ in which $\tau_{P|P_{0}} : P_{0}
\to P_{0}$ is given by $(X_{0}:Y_{0}:Z_{0}) \mapsto
(X_{0}:Y_{0}:Z_{0})$. In other words $\tau_{P|P_{0}}$ acts as the
identity on $P_{0}$ and hence $\alpha_{B}$ preserves pointwise
the cubic 
\[
f_{0} \; : \quad Y_{0}^{2}Z_{0} = X_{0}^{3} + g_{2}(1:0)X_{0}Z_{0}^{2}
+ g_{3}(1:0)Z_{0}^{3} \quad \subset B.
\]
In a similar fashion $f_{\infty}$ is a cubic in  the projective plane
\[
P_{\infty} = {\mathbb P}({\mathcal O}_{\infty}\oplus {\mathcal
O}(2)_{\infty} \oplus {\mathcal O}(3)_{\infty}).
\]
In this case the lines ${\mathcal O}_{\infty}$, ${\mathcal
O}(2)_{\infty}$ and ${\mathcal O}(3)_{\infty}$ have frames $1$,
$t_{1}^{2}$ and $t_{1}^{3}$ respectively and so $\tau_{{\mathbb
P}^{1}}$ acts trivially on ${\mathcal O}_{\infty}$ and ${\mathcal
O}(2)_{\infty}$ and by multiplication by $-1$ on ${\mathcal
O}(3)_{\infty}$. This means that if we use these frames to identify
${\mathcal O}_{\infty}$, ${\mathcal
O}(2)_{\infty}$ and ${\mathcal O}(3)_{\infty}$ with ${\mathbb C}$ we
get projective coordinates $X_{\infty} := x_{|P_{\infty}}$,
$Y_{\infty} := y_{|P_{\infty}}$ and $Z_{\infty} := z_{|P_{\infty}}$ in
which $\tau_{P|P_{\infty}}$ acts as
$(X_{\infty}:Y_{\infty}:Z_{\infty}) \mapsto
(X_{\infty}:-Y_{\infty}:Z_{\infty})$ and $f_{\infty}$ has equation
\[
Y_{\infty}^{2}Z_{\infty} = X_{\infty}^{3} +
g_{2}(0:1)X_{\infty}Z_{\infty}^{2} + g_{3}(0:1)Z_{\infty}^{3}.
\]
In other words $\alpha_{B|f_{\infty}} = (-1)_{B|f_{\infty}}$ and so
$\alpha_{B}$ has four isolated fixed points on $f_{\infty}$ coinciding
with the points of order two on $f_{\infty}$. \hfill $\Box$

\bigskip

\noindent
 Note that if we consider the
involution $\alpha_{B}\circ (-1)_{B}$ instead of $\alpha_{B}$ we will
get the same distribution of fixed points with $f_{0}$ and
$f_{\infty}$ switched, i.e. we will
get four isolated fixed points on $f_{0}$ and a trivial action on
$f_{\infty}$. 

\subsection{The quotient $B/\alpha_{B}$.}
\label{subsec-quotient-case2} 
Let $\beta : B \to
{\mathbb P}^{1}$ be a rational elliptic surface whose 
Weierstrass  model is given by (\ref{eq-weierstrass-case2}), 
with $g_{2} \in H^{0}({\mathcal O}_{{\mathbb P}^{1}}(4))$ and $g_{3}
\in H^{0}({\mathcal O}_{{\mathbb P}^{1}}(6))$ being invariant for
the standard action of $\tau_{{\mathbb P}^{1}}$. For the
time being we  will
assume that $g_{2}$ and $g_{3}$ are chosen generically so that $B = W$ is
smooth and $\beta$ has twelve $I_{1}$ fibers necessarily permuted by
$\tau_{{\mathbb P}^{1}}$.  

We have a commutative diagram
\[
\xymatrix{
B \ar[r] \ar[d]_-{\beta} & B/\alpha_{B} \ar[d]\\
{\mathbb P}^{1} \ar[r]_{\op{sq}} & {\mathbb P}^{1}
}
\]
where $\op{sq} : {\mathbb P}^{1} \to {\mathbb P}^{1}$ is the squaring
map $(t_{0}:t_{1}) \mapsto (t_{0}^{2}:t_{1}^{2})$.

Now by the analysis of the fixed points of $\alpha_{B}$ above we have
that $B/\alpha_{B} \to {\mathbb P}^{1}$ is a genus one fibration which
has six $I_{1}$ fibers. Furthermore we saw that the only singularities
of $B/\alpha_{B}$ are four singular points of type $A_{1}$  
sitting on the fiber over $\infty = (0:1) \in {\mathbb P}^{1}$. 

\begin{lem} \label{lem-B/alphaB-hat} Assume that $B$ is Weierstrass.
\begin{itemize}
\item[{\em (i)}] The minimal resolution $\widehat{B/\alpha_{B}}$ 
of $B/\alpha_{B}$ is a 
rational elliptic surface with a $6 I_{1} + I_{0}^{*}$ configuration
of singular fibers and  $B/\alpha_{B} \to {\mathbb P}^{1}$ is its
Weierstrass model.
\item[{\em (ii)}] The surface $B$ is the unique double cover of $B/\alpha_{B}$
whose branch locus consists of the fiber of $B/\alpha_{B} \to {\mathbb
P}^{1}$ over $0 = (1:0) \in {\mathbb P}^{1}$ and the four singular
points of $B/\alpha_{B}$.
\end{itemize}
\end{lem}
{\bf Proof.} By construction $\widehat{B/\alpha_{B}} \to {\mathbb
P}^{1}$ is a genus one fibered surface with seven singular fibers -
six fibers of type $I_{1}$ (i.e. the images of the twelve $I_{1}$
fibers of $\beta$ under the quotient map $B \to B/\alpha_{B}$) and one
$I_{0}^{*}$ fiber (i.e. the fiber of $\widehat{B/\alpha_{B}} \to {\mathbb
P}^{1}$ over $\infty \in {\mathbb P}^{1}$). Moreover since the section
$e : {\mathbb P}^{1} \to B$ is stable under $\alpha_{B}$ we see that
$e({\mathbb P}^{1})/\alpha_{B} \subset B/\alpha_{B}$ will again be
a section of the genus one fibration that passes through one of the
singular points. So the proper transform of $e({\mathbb
P}^{1})/\alpha_{B}$ in $\widehat{B/\alpha_{B}}$ will be a section of
$\widehat{B/\alpha_{B}} \to {\mathbb P}^{1}$ which intersects the
$I_{0}^{*}$ fiber at a point on one of the four non-multiple
components. \hfill $\Box$

\bigskip

\noindent
In fact the quotient $B \to B/\alpha_{B}$ can be constructed directly
as a double cover of the quadric $Q \cong
{\mathbb F}_{0} = \cp{1}\times \cp{1}$. In particular this gives a
geometric construction of $B$ as an iterated double cover of $Q$.

\begin{lem} \label{lem-cover-of-Q} Every rational elliptic surface with $6
I_{1} + I_{0}^{*}$ configuration of singular fibers can be obtained
as a minimal resolution of a double cover of the quadric $Q$ 
branched along a 
curve $M \in {\mathcal O}_{Q}(2,4)$ which splits as a union of two
curves of bidegrees $(1,4)$ and $(1,0)$ respectively.
\end{lem}
{\bf Proof.} Indeed consider a curve $T \subset Q$ of
bidegree $(1,4)$ and a ruling $r \subset Q$ of type $(1,0)$. Assume
for simplicity that $T$ is smooth and that $T$ and $r$ intersect
transversally. The
double cover $W_{M}$ of $Q$ branched along $M := T\cup r$ is singular 
at the ramification
points sitting over the four points in $T\cap r$. The curve $T$ is of
genus zero and so for a general $T$ the four sheeted covering map
$p_{1|T} : T \to {\mathbb P}^{1}$ will have six simple ramification
points. Thus 
\[
W_{M} \to Q \stackrel{p_{1}}{\to} {\mathbb P}^{1}
\] 
has six singular fibers
of type $I_{1}$ and one fiber passing trough the four singularities
of $W_{M}$. 

Let $s \subset Q$ be any ruling of type $(0,1)$ that
passes  trough one of the points in $T\cap r$. Then $s$ intersects $M$
at one double point and so the preimage of $s$ in $W_{M}$ splits into
two sections of the elliptic fibration $W_{M} \to {\mathbb
P}^{1}$ that intersect at one of the  singular points of $W_{M}$. This
implies (as promised) 
that the minimal resolution $\widehat{W}_{M}$ of $W_{M}$ is a
rational elliptic surface of type $6I_{1} + I_{0}^{*}$ and that 
$W_{M}$ is its Weierstrass form.

Alternatively we can construct  $\widehat{W}_{M}$ as follows. Label
the four points in $T\cap r$ as $\{ P_{1}, P_{2}, P_{3}, P_{4} \}$. Consider
the  blow-up $\phi : \widehat{Q} \to Q$ of $Q$ at the points 
$\{ P_{1}, P_{2}, P_{3}, P_{4} \}$ and let $\widehat{T}$ and $\hat{r}$
be the proper transforms of $T$ and $r$ under $\phi$. We have
\[
{\mathcal O}_{\widehat{Q}}(\widehat{T} + \hat{r}) = \phi^{*}{\mathcal
O}_{Q}(T + r)\otimes {\mathcal O}_{\widehat{Q}}\left(- 2\sum_{i=1}^{4}
E_{i}\right) 
\]
where $E_{i} \subset \widehat{Q}$ is the exceptional divisor
corresponding to the point $P_{i}$. This shows that the line bundle
${\mathcal O}_{\widehat{Q}}(\widehat{T} + \hat{r})$ is uniquely
divisible by two in $\op{Pic}(\widehat{Q})$ and so we may consider the
double cover of $\widehat{Q}$ branched along $\widehat{T} + \hat{r}$.
Since each of the rational curves $E_{i}$ intersects the branch
divisor $\widehat{T}\cup \hat{r}$ at exactly two points it follows
that the preimage $D_{i}$ of $E_{i}$ in the double cover of $\widehat{Q}$ is a
smooth rational curve of self-intersection $-2$. But if we contract
the curves $D_{i}$ we will obtain a surface with four $A_{1}$
singularities which doubly covers $Q$ with branching along $M = T\cup
r$, i.e. we will get the surface $W_{M}$. In other words the double cover of
$\widehat{Q}$ branched along $\widehat{T} + \hat{r}$ must be the
surface $\widehat{W}_{M}$. Let $\psi : W_{M} \to Q$ and $\hat{\psi} : 
\widehat{W}_{M} \to \widehat{Q}$ denote the covering maps and let
$\hat{\phi} : \widehat{W}_{M} \to W_{M}$ be the blow-up that resolves
the singularities of $W_{M}$. Hence the elliptic fibrations on $W_{M}$
and $\widehat{W}_{M}$ are given by the composition maps $\omega :=
p_{1}\circ \psi : W_{M} \to {\mathbb P}^{1}$ and $\hat{\omega} :=
p_{1}\circ \psi \circ \hat{\phi} : \widehat{W}_{M} \to {\mathbb
P}^{1}$ respectively.

Finally to write $W_{M}$ as a quotient $W_{M} = B/\alpha_{B}$
(respectively $\widehat{W}_{M}$ as a quotient $\widehat{W}_{M} =
\widehat{B/\alpha_{B}}$ we proceed as follows. If there exists a
Weierstrass rational elliptic surface $\beta : B \to {\mathbb P}^{1}$
so that $W_{M} = B/\alpha_{B}$, then $\kappa : B\to W_{M}$ will be the
unique double
cover of $W_{M}$ branched along the fiber $(W_{M})_{0} :=
\omega^{-1}(0)$ and at the four singular points of $W_{M}$. In view of
the universal property of the blow-up we may instead consider the
unique double cover $\hat{\kappa} : \widehat{B} \to \widehat{W}_{M}$
which is branched along the divisor $(\widehat{W}_{M})_{0} + \sum_{i =
1}^{4} D_{i}$. To see that such a cover exists observe that 
$\hat{\omega}^{-1}(\infty)$ is a Kodaira fiber of type $I_{0}^{*}$ and
we have $\hat{\omega}^{-1}(\infty) = 2V + \sum_{i = 1}^{4} D_{i}$,
where $2V = \hat{\psi}^{*}(\hat{r})$ is the double component of
$\hat{\omega}^{-1}(\infty)$. This yields
\[
{\mathcal O}_{\widehat{W}_{M}}\left((\widehat{W}_{M})_{0} + \sum_{i =
1}^{4} D_{i}\right) = \hat{\omega}^{*}{\mathcal O}_{{\mathbb
P}^{1}}(2)\otimes {\mathcal O}_{\widehat{W}_{M}}(-2V)
\]
and so ${\mathcal O}_{\widehat{W}_{M}}((\widehat{W}_{M})_{0} + 
\sum_{i = 1}^{4} D_{i})$ is divisible by two in $\op{Pic}(\widehat{W}_{M})$.
But from the construction of $\widehat{W}_{M}$ it follows immediately
that $\pi_{1}(\widehat{W}_{M}) = 0$ and so $\op{Pic}(\widehat{W}_{M})$
is torsion-free. Due to this there is a unique square root of the
line bundle ${\mathcal O}_{\widehat{W}_{M}}((\widehat{W}_{M})_{0} +
\sum_{i = 1}^{4} D_{i})$ and we get a unique root cover
$\hat{\kappa} : \widehat{B} \to \widehat{Q}$ as desired.

Let $\widehat{D}_{i} \subset \widehat{B}$ denote the component of the
ramification divisor of $\hat{\kappa}$ which maps to $D_{i}$. Note
that each $\widehat{D}_{i}$ is a smooth rational curve and that since 
$\hat{\kappa}^{*}D_{i} = 2 \widehat{D}_{i}$ we have 
\[
\widehat{D}_{i}\cdot \widehat{D}_{i} = \frac{1}{4}
\hat{\kappa}^{*}(D_{i}^{2}) = \frac{1}{4}\cdot 2 \cdot D_{i}^{2} = 
\frac{1}{4}\cdot 2 \cdot (-2) = -1.
\]
Therefore we can contract the disjoint $(-1)$ curves $\{
\widehat{D}_{i} \}_{i=1}^{4}$ to obtain a smooth surface $B$ which
covers $W_{M}$ two to one with branching exactly along $(W_{M})_{0}$
and the the four singular points of $W_{M}$. If we now denote the
covering involution of $\kappa : B \to W_{M}$ by $\alpha_{B}$ we
have $W_{M} = B/\alpha_{B}$ and $\widehat{W}_{M} =
\widehat{B/\alpha_{B}}$. This construction is clearly invertible, so
the lemma a is proven. \hfill $\Box$

\begin{cor} \label{cor-5parameters} All rational
elliptic surfaces $\beta : B \to \cp{1}$ which admit an involution 
$\alpha_{B}$, which preserves the zero section $e$ of $\beta$ and induces
an involution on $\cp{1}$, form a five dimensional irreducible family.
\end{cor}
{\bf Proof.} According to lemma~\ref{lem-cover-of-Q} every such
surface $B$ determines and is determined by the curve $M = T\cup r
\subset Q$ and by the choice of a smooth fiber $(W_{M})_{0}$ of
$W_{M}$. The curve $M$ depends on $\dim |{\mathcal O}_{Q}(1,4)| + \dim
|{\mathcal O}_{Q}(1,0)| - \dim \op{Aut}(Q) = 9 + 1 -  6 = 4$
parameters. Adding one more parameter for the choice of $(W_{M})_{0}$
we obtain the statement of the corollary. \hfill $\Box$

\bigskip

\noindent
It is
convenient to assemble all the
surfaces and maps  described above in the following commutative
diagram:
\[
\xymatrix@!0{
{\widehat{B}} \ar[dd]_-{\varepsilon} \ar@/_2pc/[dddd]|(.3){{\hat{\beta}}}
\ar[rr]^-{{\hat{\kappa}}} & & {\widehat{W}}_{M} 
\ar[dd]_-{\hat{\phi}} \ar@/_2pc/[dddd]|(.3){\hat{\omega}}
\ar[rr]^-{\hat{\psi}} & & {\widehat{Q}} \ar[dd]^-{\phi}
\ar@/_2pc/[dddd]|(.3){\hat{p}_{1}} \\
& & & & \\
B \ar[dd]_-{\beta} \ar'[r][rr]^(.3){\kappa} & & W_{M} 
\ar[dd]_-{\omega} \ar'[r][rr]^(.3){\psi} & & Q \ar[dd]^-{p_{1}} \\
& & & & \\
{\mathbb P}^{1} \ar[rr]_-{\op{sq}} & & {\mathbb P}^{1}
\ar[rr]_-{\op{id}} & & {\mathbb P}^{1}
}
\]
where the maps $\phi$, $\hat{\phi}$ and $\varepsilon$ are
blow-ups. The maps $\psi$, $\hat{\psi}$, $\kappa$ and $\hat{\kappa}$
are double covers and $\omega$, $\hat{\omega}$, $\beta$ and
$\hat{\beta}$ are elliptic fibrations.

Now we are ready to look for the involutions $\tau_{B}$.

Let $B$ and $\alpha_{B}$ be as in the previous section. As explained
in  Section~\ref{sss-dp9-involutions}, in order to describe all
possible involutions $\tau_{B}$ we need to describe all sections
$\zeta : {\mathbb P}^{1} \to B$ such that $\alpha_{B}^{*}\zeta = (-1)_{B}^{*}
\zeta$. 

\begin{rem} \label{rem-find-zeta} The existence of such a section
$\zeta$ can be shown by solving an equation in the group
$\MW$. For this, observe that since $\alpha_{B}$ preserves the fibers
of $\beta$ it must send a section to a section. Thus $\alpha_{B}$
induces a bijection $\alpha_{\MW} : \MW \to \MW$, which is uniquely
characterized by the property 
\[
c_{1}(\alpha_{\MW}(\mw{\xi})) = {\mathcal O}_{B}(\alpha_{B}(\xi)).
\]
Also, by the definition of $(-1)_{B}$ we know that $c_{1}(-\mw{\xi}) =
(-1)_{B}(\xi)$ and hence we need to show the existence of a section
$\zeta$, such that $\alpha_{\MW}(\mw{\zeta}) = - \mw{\zeta}$. 

The first step is to observe that since the isomorphism
$\tau_{\cp{1}}^{*}B \widetilde{\to} B$ preserves the group structure
on the fibers, the induced bijection $\alpha_{\MW}$ on sections is
actually a group automorphism.

Next note that for the general $B$ in the five dimensional
family from Corollary~\ref{cor-5parameters}, the lattice ${\mathcal T}$
has rank two since the general such $B$ has only singular  fibers of
type $I_{1}$ and so ${\mathcal T} = {\mathbb Z}e\oplus {\mathbb
Z}f$. Moreover $\alpha_{B|{\mathcal T}} = \op{id}_{{\mathcal T}}$, and
so the space of anti-invariants of 
$\alpha_{B}^{*}$ acting on $\op{Pic}(B)\otimes {\mathbb Q}$ injects
into the space of anti-invariants of $\alpha_{\MW}$. But in
Section~\ref{subsec-quotient-case2} we showed that $B/\alpha_{B}$ is
again a rational elliptic surface which has four $A_{1}$
singularities. In particular $\op{rk}(\op{Pic}(B/\alpha_{B})) = 6$ and
so  there is a $4$-dimensional space of anti-invariants for the
$\alpha_{B}^{*}$ action on $\op{Pic}(B)\otimes {\mathbb Q}$. 

This implies that $\alpha_{\MW}$ has a $4$ dimensional space of
anti-invariants on $\MW\otimes {\mathbb Q}$ and hence we can find a
section $\zeta \neq e$ with $\alpha_{\MW}([\zeta]) = -[\zeta]$. 
The involution $\tau_{B}$ corresponding to $(\alpha_{B},\zeta)$ will
have only four isolated fixed points. 
\end{rem}

\section{The four dimensional subfamily of special rational elliptic
surfaces}  \label{s-4dim-family}

From now on we will restrict our attention to a 4-dimensional
subfamily of the 5-dimensional family of surfaces of
Corollary~\ref{cor-5parameters}. We do this for two
reasons:

\begin{itemize}
\item  Mathematically, this seems to be the simplest family where the
full range of possible behavior of the spectral involution
$\T = \FM^{-1}\circ \tau_{B}^{*} \circ \FM$ is present, see
Proposition~\ref{prop-T}. Indeed, for a generic surface in the five
dimensional family, $\T$ takes line bundles to line bundles, so
everything can be rephrased without the use of the derived category. 
\item In terms of our motivation from the physics, this specialization
is needed for the construction of the Standard Model bundles. 
By taking fiber products of 
surfaces from the five dimensional family 
one indeed gets a smooth Calabi-Yau with a
freely acting involution.  
However, it turns out that for a
generic such $B$, the cohomology of the resulting
Calabi-Yau is not rich enough to lead to invariant vector bundles
satisfying the Chern class constraints from \cite{dopw-ii}. 
\end{itemize}

\subsection{The quotient $B/\tau_{B}$}
\label{subsec-quotient-tauB-case2}

The starting point of the construction of the four dimensional family is
the following simple observation: since $\zeta$ must satisfy  
$\alpha_{B}^{*}(\zeta) = (-1)_{B}^{*}(\zeta)$ it will help to work with
rational elliptic surfaces $B$ for which we know the geometric
relationship between the two involutions $\alpha_{B}$ and $(-1)_{B}$. 
In the previous section we
interpreted  the involution
$\alpha_{B}$ as the covering involution of the map $\kappa$. On the
other hand the involution $(-1)_{B}$ was the group inversion along the
fibers of $\beta$ corresponding to a zero section $e : {\mathbb P}^{1}
\to B$ which was chosen to be one of the two components of the
preimage in $B$ of a ruling of type $(0,1)$ in $Q$ which passes trough
one of the four points in $T\cap r$. Since in this setup the
involutions $\alpha_{B}$ and $(-1)_{B}$ are generically unrelated it
is natural to look for a special configuration of the curves $T$ and
$r$ for which $(-1)_{B}$ can be related to the maps $\kappa$ and
$\psi$. 

\begin{lem} \label{lem-4parameters}Consider the family of 
rational elliptic surfaces $B$ obtained as an
iterated double cover $B \to W_{M} \to Q$ for which the component $T$ of the
branch curve $M$ is split further into a union $T = s\cup {\mathfrak
T}$ where $s$ is a ruling of $Q$ of type $(0,1)$ and ${\mathfrak T}$
is a curve of type $(1,3)$. Let as before $e$ be the section of $B$
mapping to $s \subset Q$. Then we have:
\begin{itemize}
\item[{\em (i)}] The involution $(-1)_{B,e}$ is a lift of the covering
involution 
of the double cover $\psi : W_{M} \to {\mathbb P}^{1}$. 
\item[{\em (ii)}] For a general pair $(B,\alpha_{B})$ corresponding to a
branch curve $M = s \cup {\mathfrak T} \cup r$  there exist three
pairs of sections of $\beta$ labeled by the non-trivial points of order two on
$f_{0}$ and such that the two members of each pair are interchanged
both by $\alpha_{B}$ and $(-1)_{B}$. 
\end{itemize}
\end{lem}
{\bf Proof.}  
If the curve ${\mathfrak T}$ is chosen to be general and
smooth, then the branch curve $M$ has five nodes $\{ P, P_{1},
P_{2}, P_{3}, P_{4} \}$. Here as before $\{ P_{1}, P_{2}, P_{3}, P_{4}
\} = T\cap r$  and the extra point $P$ is the intersection point
of the curves ${\mathfrak T}$ and $s$. 

Let $\{ p, p_{1}, p_{2}, p_{3}, p_{4} \} \subset W_{M}$ denote the
corresponding singularities of $W_{M}$. Observe that for a general
choice of the curve ${\mathfrak T}$ and the point $0 \in {\mathbb
P}^{1}$ the singularity $p \in W_{M}$ is not contained in the branch
locus $(W_{M})_{0}\cup \{ p_{1}, p_{2}, p_{3}, p_{4} \}$ of the map
$\kappa$. In particular the double cover of $W_{M}$ branched along 
$(W_{M})_{0}\cup \{ p_{1}, p_{2}, p_{3}, p_{4} \}$ will have two
$A_{1}$ singularities at the two preimages $\bar{p}_{1}$ and
$\bar{p}_{2}$  of the point
$p$. In order to get a smooth rational elliptic surface we have to
to blow up this two points. Abusing slightly the notation we will
denote by $B$ the resulting smooth surface and by $\kappa : B \to
W_{M}$ the composition of the blow-up map with the double cover
of $W_{M}$ branched along $(W_{M})_{0}\cup \{ p_{1}, p_{2}, p_{3},
p_{4} \}$. Let $n_{1}, n_{2} \subset B$ denote the exceptional curves
corresponding to $\bar{p}_{1}$ and $\bar{p}_{2}$ and let $o_{1},
o_{2}$  denote proper
transforms in $B$ of the two preimages of the fiber 
$\omega^{-1}(\omega(p))$ in the double cover of $W_{M}$ branched along 
$(W_{M})_{0}\cup \{ p_{1}, p_{2}, p_{3}, p_{4} \}$. Here we have
labeled $o_{1}$ and $o_{2}$ so that $\bar{p}_{1} \in o_{1}$ and
$\bar{p}_{2} \in o_{2}$. From
this picture it is clear
that $\beta : B \to {\mathbb P}^{1}$ is a smooth rational elliptic
surface with a $8I_{1} + 2I_{2}$ configuration of singular fibers
which is symmetric with respect to the involution $\tau_{{\mathbb
P}^{1}}$. Furthermore the two $I_{2}$ fibers of $\beta$ are just the
curves $o_{1}\cup n_{1}$ and $o_{2}\cup n_{2}$ and the two fixed points $\{0,
\infty\}$ of $\tau_{{\mathbb P}^{1}}$ correspond to two smooth fibers
$f_{0}$ and $f_{\infty}$ of $\beta$. Note also that the proper
transform of the section $s \subset Q$ via the generically finite map
$\psi\circ \kappa : B \to Q$ is an irreducible rational curve $e
\subset B$ which is a section of $\beta : B \to {\mathbb
P}^{1}$. Moreover the inversion $(-1)_{B}$ with respect to $e$
commutes with the covering involution $\alpha_{B}$ for the map
$\kappa$ and descends to an inversion $(-1)_{W_{M}}$ along the fibers
of the elliptic fibration $\omega : W_{M} \to {\mathbb P}^{1}$ which
fixes the image of $e$ pointwise. But by construction the image of $e$
in $W_{M}$ is just the component of the ramification divisor of the
cover $\psi : W_{M} \to Q$ sitting over $s \subset Q$. In particular
$(-1)_{W_{M}}$ is just the covering involution for the map $\psi$.

We are now ready to construct a section $\zeta : {\mathbb P}^{1} \to
B$ of $\beta$ satisfying $\alpha_{B}^{*}(\zeta) =
(-1)_{B}(\zeta)$. Indeed, assume that such a section exists.

\begin{figure}[!ht]
\begin{center}
\psfig{file=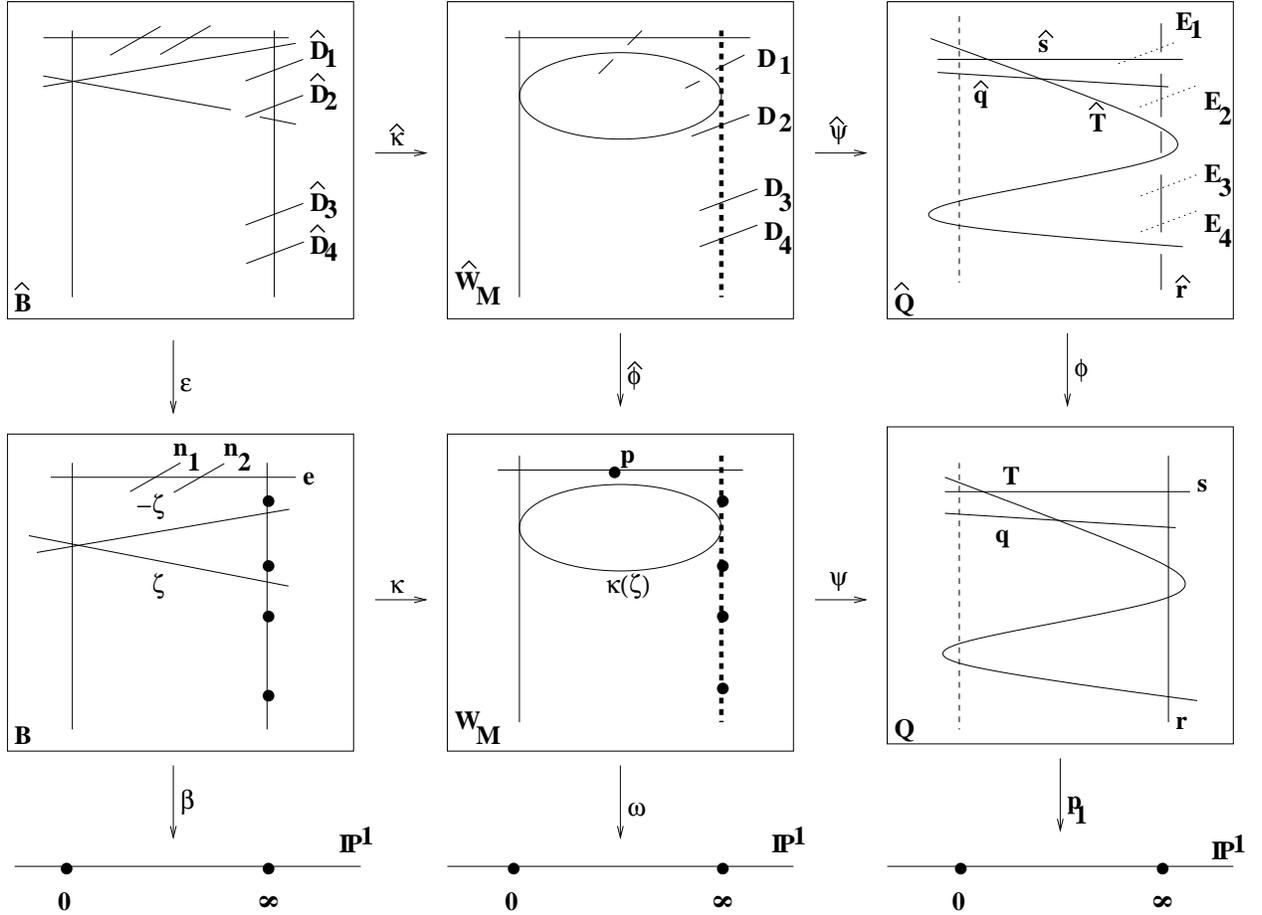,width=6.5in} 
\end{center}
\caption{The section $\zeta$}
\label{fig7} 
\end{figure}

Due to the fact that
$\alpha_{B|f_{0}} = \op{id}_{f_{0}}$  we have $\zeta(0) = - \zeta(0)$
i.e. $\zeta(0)$ is a point of order two on $f_{0}$. Now from the
Weierstrass equation (\ref{eq-weierstrass-case2}) of $B$ it is clear
that the general $B$ cannot have monodromy $\Gamma_{0}(2)$ and so
without a loss of generality we may assume that $\zeta \neq -
\zeta = \alpha_{B}^{*}\zeta$. 
Consider now the image $\kappa(\zeta) \subset W_{M} = B/\alpha_{B}$ 
of $\zeta$ in $W_{M}$. We have $\kappa^{-1}(\kappa(\zeta)) = \zeta
\cup \alpha_{B}^{*}\zeta$. On the other hand the preimage of the general
elliptic fiber of $\omega : W_{M} \to {\mathbb P}^{1}$ via $\kappa$ splits as a
disjoint union of two fibers of $\beta$ and so
$\alpha_{B|f_{0}} = \op{id}_{f_{0}}$  we have $\zeta(0) = - \zeta(0)$
i.e. $\zeta(0)$ is a point of order two on $f_{0}$. 
Consider now the image $\kappa(\zeta) \subset W_{M} = B/\alpha_{B}$ 
of $\zeta$ in $W_{M}$. We have $\kappa^{-1}(\kappa(\zeta)) = \zeta
\cup \alpha_{B}^{*}\zeta$. On the other hand the preimage of the general
elliptic fiber of $\omega : W_{M} \to {\mathbb P}^{1}$ via $\kappa$ splits as a
disjoint union of two fibers of $\beta$ and so
\[
\kappa(\zeta)\cdot \omega^{-1}(\op{pt}) =
\frac{1}{2}\kappa^{*}(\kappa(\zeta)\cdot \omega^{-1}(\op{pt})) =
\frac{1}{2} (\zeta + \alpha^{*}\zeta)\cdot (2 \beta^{-1}(\op{pt})) = 2
\]
i.e. the smooth rational curve $\kappa(\zeta)$ is a double section of
$\omega$. Moreover the condition $\alpha_{B}^{*}\zeta = - \zeta$
combined with the property $\alpha_{B|B_{\infty}} = (-1)_{B_{\infty}}$
implies that $(\alpha_{B}^{*})\zeta(\infty) = \zeta(\infty)$ and so the double
cover $\omega_{|\kappa(\zeta)} : \kappa(\zeta) \to {\mathbb P}^{1}$ is
branched exactly over the points $0, \infty$. Furthermore since
$\zeta(0)$ is a point of order two on $f_{0}$ it must lie on the
preimage of $T$ in $B$ and so the two ramification points of the cover 
$\omega_{|\kappa(\zeta)} : \kappa(\zeta) \to {\mathbb P}^{1}$ must
both lie on the ramification divisor of the double cover $\psi : W_{M}
\to Q$ as depicted on Figure~\ref{fig7}.

Also note that if we pullback to $B$ the involution of $W_{M}$  
acting along the fibers of $\psi$ we will get precisely
$(-1)_{B}$. Combined with the fact that $\alpha_{B}^{*}\zeta =
(-1)_{B}^{*}\zeta$ this shows that $\kappa(\zeta)$ is stable under the
involution of $W_{M}$  
acting along the fibers of $\psi$ and so $\psi^{-1}(\psi(\kappa(\zeta)))
= \kappa(\zeta)$. Put $q := \psi(\kappa(\zeta))$. Then $q$ is a smooth
rational curve which intersects each of the curves 
$T$ and $r$ at a single point so that the double cover 
$\psi_{|\kappa(\zeta)} : \kappa(\zeta) \to q$ is branched exactly at
$q \cap (T\cup r)$. So $q$ is the unique ruling of type $(0,1)$ on $Q$
which passes trough the point $\psi(\kappa(\zeta(0))) \in T\cap Q_{0}$.

Conversely if we start with any ruling $q$ of type $(0,1)$ 
that passes trough one of the four points in $T\cap f_{0}$ we see that
$\psi^{-1}(q)$ is a smooth rational curve which is a double cover of
$q$ with branch divisor $q \cap (T\cup r)$. Since the rulings of type
$(1,0)$ pull back to a single fiber of $\omega$ via $\psi$ we see that
\[
\psi^{-1}(q)\cdot \omega^{-1}(\op{pt}) = 
\psi^{*}(q\cdot p_{1}^{-1}(\op{pt})) = 2
q\cdot p_{1}^{-1}(\op{pt}) = 2,
\]
and so $q$ is a double
section of the elliptic fibration $\omega : W_{M} \to {\mathbb P}^{1}$
which is tangent to the fibers $(W_{M})_{0}$ and
$(W_{M})_{\infty}$. Also it is clear that for $T$ and $r$ in general
position the point $q\cap r$ is not one of the four points in $T\cap
r$ and so the point of contact of $\psi^{-1}(q)$ and
$(W_{M})_{\infty}$ is not one of the four isolated branch points of
the covering $\kappa : B \to W_{M}$. So $\psi^{-1}(q)$ intersects the
branch locus of $\kappa$ at a single point with multiplicity two -
namely the point of contact of $(W_{M})_{0}$ and $\psi^{-1}(q)$. This
implies that the preimage of $\psi^{-1}(q)$ in $B$ splits into two
sections of $\beta$ that intersect at a point on the fiber $f_{0}$ and
are exchanged both by $\alpha_{B}$ and $(-1)_{B}$. The lemma is
proven. \hfill $\Box$

\bigskip

\noindent
Finally, let $\tau_{B}$ be the involution of $B$ corresponding to the
pair $(\alpha_{B},\zeta)$ constructed in the previous lemma. Then
the quotient $B/\tau_{B}$ is again a genus one fibered
rational surface which similarly to $B/\alpha_{B}$ has four $A_{1}$
singularities all sitting on fiber over $\infty \in \cp{1}$. However
$B/\tau_{B}$ has also a smooth double
fiber and so is only genus one fibered. The minimal resolution of
$B/\tau_{B}$ in this case has a $4I_{1} + I_{2} + I_{0}^{*} + 
_{2}\!I_{0}$ configuration of singular fibers.

\subsection{The basis in $H^{2}(B,{\mathbb Z})$} 
\label{subsec-case2-basis-in-h2} 

In order to describe an integral basis of the cohomology of $B$
we need to find a description of our $B$ as a blow-up of ${\mathbb
P}^{2}$ in the base points of a pencil of cubics. 

To achieve this we will use a different fibration on $B$, namely the
fibration 
\[
\xymatrix@1{
B \ar[r]^-{\psi\circ\kappa} \ar@/_1pc/[rr]_-{\delta} &  Q
\ar[r]^-{p_{2}} & {\mathbb P}^{1}.
}
\]
induced from the projection of the quadric $Q$ onto its {\em second}
factor. The fibers of $\delta$ can be studied directly in terms of the
degree four map $\psi\circ \kappa : B \to Q$ but it is much more
instructive to use instead an alternative description of $B$ as a {\em
double cover} of a quadric. 

In section~\ref{subsec-quotient-tauB-case2} we saw that the
description of $B$ as an iterated double cover 
\[
B \stackrel{\kappa}{\to} W_{M} \stackrel{\psi}{\to} Q
\]
of the quadric $Q$ yields two commuting involutions $\alpha_{B}$ and
$(-1)_{B}$ on $B$. By construction the quotient $B/\alpha_{B}$ can be
identified with the blow-up of the rational elliptic surface $W_{M}$
at the $A_{1}$ singularity $p \in W_{M}$ sitting over the unique
intersection point $\{ P \} = s \cap {\mathfrak T}$. In particular if 
we consider the Stein factorization of the generically finite map 
$\kappa : B \to W_{M}$ we get
\[
B \to W_{\beta} \to W_{M}.
\]
Here $W_{\beta}$ is the Weierstrass model of $\beta : B \to {\mathbb
P}^{1}$ and $B \to W_{\beta}$ is the blow-up the two 
$A_{1}$ singularities of $W_{\beta}$ and the map 
$W_{\beta} \to W_{M}$ is the double cover branched at 
$(W_{M})_{0}\cup \{ p_{1}, p_{2}, p_{3}, p_{4} \}$.

Similarly we can describe the quotients $B/(-1)_{B}$
and $B/((-1)_{B}\circ \alpha_{B})$ as blow-ups of appropriate double  
covers of $Q$. Indeed the curves on $Q$ that
play a special role in the description of $B$ as an iterated double
cover are: the $(1,3)$ curve ${\mathfrak T}$, the $(0,1)$
ruling $s$ and the $(1,0)$ rulings $r = r_{\infty} =
p_{1}^{-1}(\infty)$ and $r_{0} = p_{1}^{-1}(0)$.

Consider the double cover $\omega' : W_{M'} \to Q$ branched along the
curve $M' =  T \cup r_{0} = s \cup {\mathfrak T} \cup r_{0}$ and the
double cover $\op{Sq} : \widetilde{Q} \to Q$ branched along the union
of rulings $r_{0} \cup r_{\infty}$. Clearly $\widetilde{Q}$ is 
again a quadric which is just a the fiber product of $p_{1} : Q \to
{\mathbb P}^{1}$ with the squaring map $\op{sq} : {\mathbb P}^{1} \to
{\mathbb P}^{1}$, i.e. we have a fiber-square
\[
\xymatrix{
{\widetilde{Q}} \ar[r]^-{\op{Sq}} \ar[d]_-{\tilde{p}_{1}} & Q
\ar[d]^-{p_{1}} \\
{\mathbb P}^{1} \ar[r]_-{\op{sq}} & {\mathbb P}^{1}
}
\]
The preimage $\widetilde{{\mathfrak T}} := \op{Sq}^{-1}({\mathfrak T})  
\subset \widetilde{Q}$ of ${\mathfrak T}$ in $\widetilde{Q}$ is a
genus two curve doubly covering ${\mathfrak T}$ with branching at the
six points ${\mathfrak T}\cap (r_{0}\cup r_{\infty})$. Also, the
preimage $\tilde{s} = \op{Sq}^{-1}(s)$ is a rational curve doubly
covering the ruling $s$ branched at the two points $s \cap (r_{0} \cup
r_{\infty})$. In particular, $\tilde{s}$ is a ruling of type $(0,1)$
on $\widetilde{Q}$. Similarly, if we denote by $\tilde{r}_{0}$ and
$\tilde{r}_{\infty}$ the two components of the ramification divisor of
$\op{Sq} : \widetilde{Q} \to Q$, then $\tilde{r}_{0}$ and
$\tilde{r}_{\infty}$ are rulings of type $(1,0)$ on $\widetilde{Q}$.

Now it is clear that the Weierstrass model $W_{\beta}$ of $B$
can be described as either of the following

\begin{itemize}
\item $W_{\beta} \to W_{M}$ is the double cover branched at the fiber
$(W_{M})_{0}$ and the four points $\{p_{1}, p_{2}, p_{3}, p_{4} \}$ of
order two of the fiber $(W_{M})_{\infty}$. 
\item $W_{\beta} \to W_{M'}$ is the double cover branched at the fiber
$(W_{M})_{\infty}$ and the four points of
order two of the fiber $(W_{M})_{0}$.
\item $W_{\beta} \to \widetilde{Q}$ is the double cover branched at
the curve $\tilde{s}\cup \widetilde{{\mathfrak T}}$.
\end{itemize} 
Furthermore 
\begin{itemize}
\item The quotient $B/\alpha_{B} \to W_{M}$ is the blow-up of $W_{M}$ at the
$A_{1}$ singularity $p$ sitting over the point $P \in Q$ of
intersection of $s$ and ${\mathfrak T}$. The map $B \to B/\alpha_{B}$
is the double cover of $B/\alpha_{B}$ branched at the fiber
$(B/\alpha_{B})_{0}$ and the four points of order two of
$(B/\alpha_{B})_{\infty}$. 
\item The quotient $B/(\alpha_{B}\circ (-1)_{B}) \to W_{M'}$ is the  
blow-up of $W_{M'}$ at the
$A_{1}$ singularity sitting over the point of
intersection of $s$ and ${\mathfrak T}$. The map $B \to
B/(\alpha_{B}\circ (-1)_{B})$
is the double cover of $B/(\alpha_{B}\circ (-1)_{B})$ branched at the fiber
$(B/\alpha_{B})_{\infty}$ and the four points of order two of
$(B/\alpha_{B})_{0}$.
\item The quotient $B/(-1)_{B}$ is the blow-up of $\widetilde{Q}$ at
the two intersection points of $\tilde{s}$ and $\widetilde{{\mathfrak
T}}$. The map $B \to B/(-1)_{B}$ is the double cover branched at the
strict transform of $\tilde{s}\cup \widetilde{{\mathfrak T}}$.
\end{itemize}

\

\bigskip

The action of the Klein group $\langle \alpha_{B}, (-1)_{B}
\rangle$ on $B$ and all of the above maps are most conveniently
recorded in the commutative diagram
\begin{equation} \label{eq-allmaps}
\xymatrix@C=4pt@R=6pt{
& &  B \ar@/_1pc/[ddll] \ar[dd] \ar@/^1pc/[ddrr] \ar@{.>}[dr] & & \\
& &  & W_{\beta}\ar@{.>}@/_1pc/[dddr] \ar@{.>}@/^0.5pc/[dddl]
\ar@/_2pc/@{.>}[dddlll]  & \\
B/\alpha_{B} \ar[dd] & &  B/(-1)_{B} \ar[dd] & &
B/(\alpha_{B}\circ (-1)_{B})
\ar[dd] \\
&& && \\ 
W_{M} \ar[ddrr]_-{\psi} & & {\widetilde{Q}} \ar[dd]_-{\op{Sq}} 
& & W_{M'} \ar[ddll]^-{\psi'} \\
&&&& \\
& & Q & &
}
\end{equation}
where the solid arrows in the first and third rows are all double covers, the
solid arrows in the  middle row are blow-ups and the dotted arrows are
Stein factorization maps.

\begin{figure}[!ht]
{\Large
\[
\xymatrix@C=2.5in@R=1.5in{ & W_{\beta}
\ar[ld]|-{\boxed{\epsfig{file=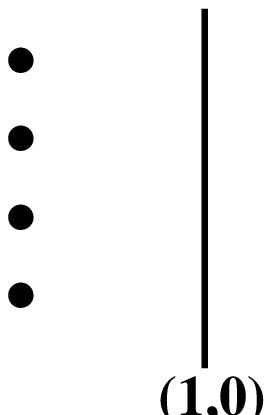,height=.75in}}} 
\ar[d]|-{\boxed{\epsfig{file=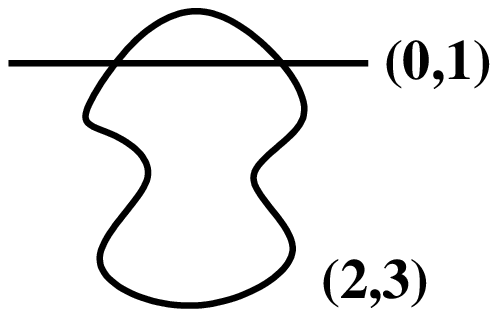,height=.75in}}}
\ar[rd]|-{\boxed{\epsfig{file=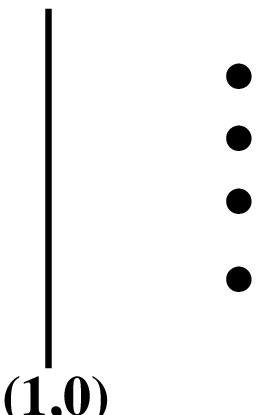,height=.75in}}}  &\\
W_{M} \ar[rd]|-{\boxed{\epsfig{file=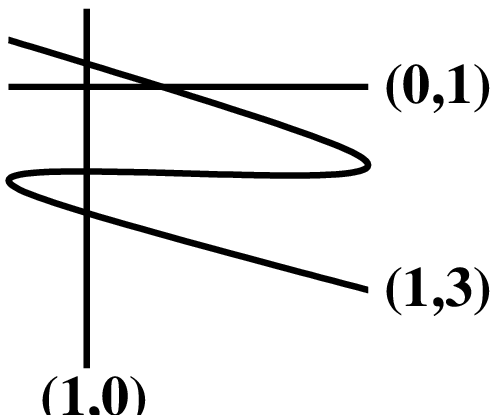,height=.75in}}}  &  
{\widetilde{Q}} \ar[d]|-{\boxed{\epsfig{file=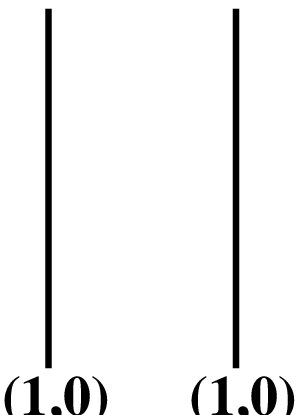,height=.75in}}} & 
W_{M'} \ar[ld]|-{\boxed{\epsfig{file=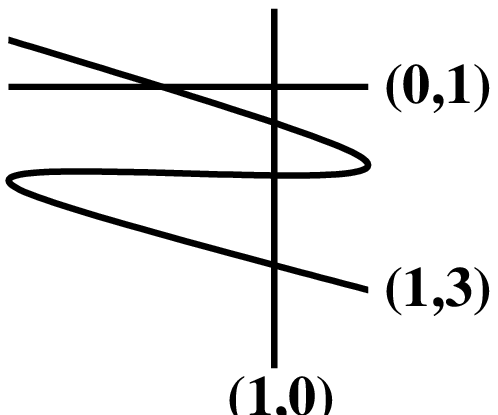,height=.75in}}}\\ 
& Q &
}
\]
}
\caption{$W_{\beta}$ as a double cover of a quadric}
\label{fig-maps} 
\end{figure}

In order to visualize the system of maps (\ref{eq-allmaps}) better it
is instructive to label all the double cover maps  appearing in
(\ref{eq-allmaps}) by a picture of their branch loci. This is
recorded in the diagram in Figure~\ref{fig-maps}.

There is a definite advantage in interpreting geometric questions on $B$
or $W_{\beta}$ on all three surfaces $W_{M}$, $W_{M'}$ and
$\widetilde{Q}$. For example, by viewing  $W_{\beta}$ as a
double  cover of the quadric $\widetilde{Q}$ we can easily describe
the fibers of the rational curve fibration $\delta : B \to {\mathbb
P}^{1}$ defined in the beginning of the section. Indeed, due to the
commutativity of (\ref{eq-allmaps}) the map $\delta = p_{2}\circ
\psi\circ \kappa$ decomposes also as
\[
\xymatrix@1{
B \ar[r] \ar@/_1pc/[rrr]_-{\delta} & B/(-1)_{B} \ar[r] &
{\widetilde{Q}} \ar[r]^-{\tilde{p}_{2}} & {\mathbb P}^{1},
}
\] 
where $\tilde{p}_{2} : \widetilde{Q} \to {\mathbb P}^{1}$ is the
projection onto the ruling of type $(1,0)$. 
In particular we can view 
each  fiber $\delta^{-1}(x)$ of the map $\delta : B \to {\mathbb
P}^{1}$ as the double cover 
of the fiber $\tilde{p}_{2}^{-1}(x)$ of $\tilde{p}_{2} : \widetilde{Q} \to
{\mathbb P}^{1}$ branched along the degree two divisor
$\widetilde{{\mathfrak T}} \cap \tilde{p}_{2}^{-1}(x) \subset
\tilde{p}_{2}^{-1}(x)$. This shows that the singular fibers of
$\delta$ are precisely the preimages under the map $B \to
\widetilde{Q}$ of $\tilde{s}$ and of those $(0,1)$ rulings of
$\widetilde{Q}$  which happen to be tangent to the curve
$\widetilde{{\mathfrak T}}$.

\begin{figure}[!ht]
\begin{center}
\epsfig{file=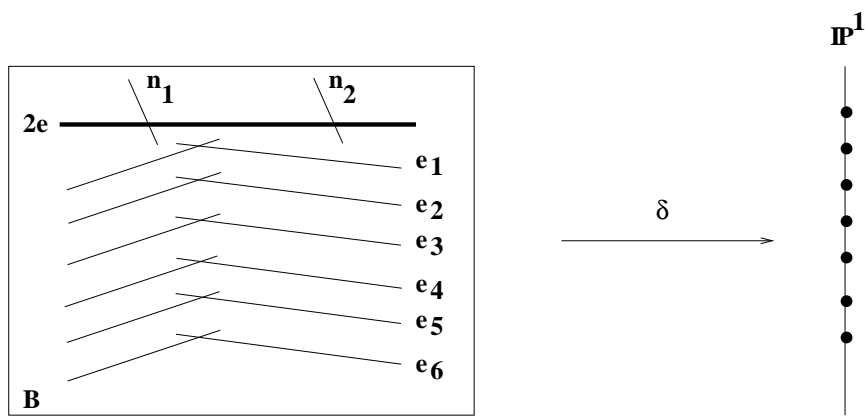,width=4.5in} 
\end{center}
\caption{The singular fibers of $\delta$}\label{fig8} 
\end{figure}

Since the curve $\widetilde{{\mathfrak T}}$ is of type $(2,3)$ on
$\widetilde{Q}$ we see by adjunction that $\widetilde{{\mathfrak T}}$
must have genus two and so by the Hurwitz formula the double cover map
$\tilde{p}_{2} : \widetilde{{\mathfrak T}} \to {\mathbb P}^{1}$ will
have six ramification points. This means that there are six rulings of
$\widetilde{Q}$ of type $(0,1)$ which are tangent to
$\widetilde{{\mathfrak T}}$, i.e. generically 
$\delta$  will have seven singular fibers (see
Figure~\ref{fig8}). Six of those
will be unions of two rational curves meeting at a point and the
seventh one will have one rational component occurring with
multiplicity two (the preimage in $B$ of the strict transform of
$\tilde{s}$ in $B/(-1)_{B}$) and two reduced
rational components $n_{1}$ and $n_{2}$ (the exceptional divisors of
the blow-up $B \to W_{\beta}$). Notice moreover that
(\ref{eq-allmaps}) implies that the preimage in
$B$ of the strict transform of $\tilde{s}$ in $B/(-1)_{B}$ is
precisely the zero section $e$ of the elliptic fibration $\beta : B
\to {\mathbb P}^{1}$ and so the non-reduced fiber of $\delta$ is just
the divisor $2e + n_{1} + n_{2}$ on $B$.
In fact, one can describe explicitly the $(0,1)$ rulings of
$\widetilde{Q}$ that are tangent to the curve $\widetilde{{\mathfrak
T}}$. Indeed let $\op{pt} \in r_{0}\cap {\mathfrak T}$ be one of the
three intersection points of $r_{0}$ and ${\mathfrak T}$. Choose
(analytic) local coordinates $(x,y)$ on a neighborhood $\op{pt}\in U
\subset Q$  so that $\op{pt} = (0,0)$, $r_{0}$ has equation $x = 0$ in
$U$ and the $(0,1)$ ruling through $\op{pt} \in Q$ has equation $y = 0$
in $U$. Let
$\widetilde{U} \subset \widetilde{Q}$ be the preimage of $U$ in
$\widetilde{Q}$. Then there are unique coordinates $(u,v)$ on
$\widetilde{U}$ such that the double cover  $\widetilde{U}
\to U$ is given by $(u,v) \mapsto (u^{2},v) = (x,y)$. Due to our
genericity assumption\footnote{We are assuming that ${\mathfrak T}$
meets $r_{0}$ and $r_{\infty}$ transversally.} the local equation of
${\mathfrak T}$ in $U$ will be $x = ay + (\text{higher order terms})$
for some number $a$. Thus the pullback of $r_{0}$ to $\widetilde{U}$ will
be given by $u = 0$ and  $\widetilde{\mathfrak T}$ will have equation
$u^{2} = av + (\text{higher order 
terms})$. Since by construction $v = 0$ is the local equation of a
$(0,1)$ ruling of $\widetilde{Q}$ it follows that
$\widetilde{\mathfrak T}$ is tangent to the three $(0,1)$ rulings of
$\widetilde{Q}$ passing through the three intersection points in
$\widetilde{{\mathfrak T}}\cap \tilde{r}_{0}$. In the same way one
sees that $\widetilde{\mathfrak T}$ is tangent to the three $(0,1)$ rulings of
$\widetilde{Q}$ passing through the three intersection points in
$\widetilde{{\mathfrak T}}\cap \tilde{r}_{\infty}$. This accounts for
all six $(0,1)$ rulings of $\widetilde{Q}$ that are tangent to
$\widetilde{{\mathfrak T}}$. 

We are now ready to describe $B$ as the blow-up of ${\mathbb P}^{2}$
at the base locus of a pencil of cubics. Each component of a
reduced singular fiber of $\delta$ is a curve of self-intersection
$(-1)$ on $B$. For every such fiber choose one of the components and
label it by $e_{i}$, $i = 1, 2, \ldots, 6$ (see
Figure~\ref{fig8}). Now $e, e_{1}, e_{2}, \ldots , e_{6}$ is a
collection of seven disjoint $(-1)$ curves on the rational elliptic
surface $B$. The curves $n_{1}$ and $n_{2}$ are rational $(-2)$ curves
on $B$ and so if we contract $e$ each of them will become a $(-1)$
curve. So if we contract $e, e_{1}, e_{2}, \ldots , e_{6}$ and after
that we contract $n_{1}$ we will end up with a Hirzebruch
surface. Moreover numerically $e, e+n_{1}, e_{1}, e_{2}, \ldots ,
e_{6}$ behave like eight disjoint $(-1)$ curves on $B$ and so the
result of the contraction of $e, n_{1}, e_{1}, e_{2}, \ldots ,
e_{6}$ should be ${\mathbb F}_{1}$. Contracting the infinity section 
of ${\mathbb F}_{1}$ we will finally obtain ${\mathbb P}^{2}$ as the
blow down of nine $(-1)$ divisors on $B$. Let
$e_{7}$ denote the infinity section of ${\mathbb F}_{1}$. To make
things explicit let us identify $e_{7}$ as a curve  coming from
$\widetilde{Q}$. Denote by ${\mathfrak e} \subset \widetilde{Q}$ the
image of $e_{7}$ in $\widetilde{Q}$. Then ${\mathfrak e}$ is an
irreducible curve which intersects the generic  $(0,1)$ ruling at one
point. This implies that ${\mathfrak e}$ is of type $(1,k)$ on
$\widetilde{Q}$ and so ${\mathfrak e}$ must be a rational curve. In
particular the map $e_{7} \to {\mathfrak e}$ ought to be an
isomorphism and $e_{7}\cup (-1)_{B}^{*}(e_{7})$ is the preimage in $B$
of the strict transform of ${\mathfrak e}$ in
$B/(-1)_{B}$. Equivalently $e_{7}\cup (-1)_{B}^{*}(e_{7})$ is the
strict transform in $B$ of the preimage of ${\mathfrak e}$ in
$W_{\beta}$. This implies that the preimage of ${\mathfrak e}$ in
$W_{\beta}$ is reducible and so ${\mathfrak e}$ must have order of
contact two with the branch divisor $\tilde{s} \cup
\widetilde{{\mathfrak T}}$ of the covering $W_{\beta} \to
\widetilde{Q}$  at each point where ${\mathfrak e}$ and $\tilde{s}\cup
\widetilde{{\mathfrak T}}$ meet. Since ${\mathfrak e}\cdot\tilde{s} =
(1,k)\cdot (0,1) = 1$ this implies that ${\mathfrak e}$ must pass
through one of the two intersection points of $\tilde{s}\cap
\widetilde{{\mathfrak T}}$ and be tangent to $\widetilde{{\mathfrak
T}}$ at $({\mathfrak e}\cdot \widetilde{{\mathfrak T}} -1)/2$
points. But 
\[
\frac{{\mathfrak e}\cdot \widetilde{{\mathfrak T}} -1}{2} = \frac{(1,k)\cdot
(2,3) - 1}{2} = k + 1
\]
and so $e_{7}\cdot (-1)_{B}^{*}e_{7} = k +1$. From here we can
calculate $k$. Indeed, on one hand we know that $e_{7}^{2} = -1$ and
so 
\[
(e_{7} + (-1)_{B}^{*}e_{7})^{2} = -2 + 2 + 2k = 2k.
\]
On the other hand $e_{7} + (-1)_{B}^{*}e_{7}$ is the preimage in $B$
of the strict transform of ${\mathfrak e}$ in $B/(-1)_{B}$. But
$B/(-1)_{B}$ is simply the blow-up of $\widetilde{Q}$ at the two
intersection points of $\tilde{s}$ and $\widetilde{{\mathfrak T}}$ and
${\mathfrak e}$ passes trough only one of those points and so the 
strict transform of ${\mathfrak e}$ in $B/(-1)_{B}$ has
self-intersection ${\mathfrak e}^{2} - 1$. In other words
\[
(e_{7} + (-1)_{B}^{*}e_{7})^{2} = 2({\mathfrak e}^{2} - 1) = 2(2k -1)
= 4k -2,
\]
and so $k = 1$. 

Therefore, in order to reconstruct $e_{7}$ starting from
$\widetilde{Q}$ we need to find a $(1,1)$ curve ${\mathfrak e}$ on
$\widetilde{Q}$ which passes through one of the two points in
$\tilde{s}\cap \widetilde{{\mathfrak T}}$ and tangent to
$\widetilde{{\mathfrak T}}$ at two extra points. But curves like that
always exist. Indeed, the linear system $|{\mathcal
O}_{\widetilde{Q}}(1,1)|$ embeds $\widetilde{Q}$ in ${\mathbb
P}^{3}$. Pick a point $J \in \tilde{s}\cap
\widetilde{{\mathfrak T}}$ and let $j : \widetilde{Q} \dashrightarrow
{\mathbb P}^{2}$ be the linear projection of $\widetilde{Q}$ from that point. 
Now the $(1,1)$-curves passing through $J$ are precisely the preimages
via $j$  of all lines in ${\mathbb P}^{2}$ and so the curve
${\mathfrak e}$ will be just the preimage under $j$ of a line in
${\mathbb P}^{2}$ which is bitangent to $j(\widetilde{{\mathfrak
T}})$. To understand  better the curve $j(\widetilde{{\mathfrak
T}}) \subset {\mathbb P}^{2}$ note that it has degree $(1,1)\cdot
(2,3) - 1 = 4$ and that the map $j : \widetilde{{\mathfrak T}} \to 
j(\widetilde{{\mathfrak T}})$ is a birational morphism. Furthermore
any $(1,1)$-curve passing trough $J$ and another point on the $(1,0)$
ruling  through $J$ will have to contain the whole $(1,0)$
ruling. Since the $(1,0)$ ruling trough $J$ intersects
$\widetilde{{\mathfrak T}}$ at $J$ and two extra points$J'$ and $J''$, 
it follows that $j(J') = j(J'')$. Therefore 
$j(\widetilde{{\mathfrak T}})$ is a nodal quartic in ${\mathbb P}^{2}$
and the curve ${\mathfrak e} \subset \widetilde{Q}$ corresponds to a
bitangent line of this nodal quartic. The normalization of this
nodal quartic is just the genus two curve $\widetilde{{\mathfrak T}}$
and the lines in ${\mathbb P}^{2}$ correspond just to sections in the
canonical class $\omega_{\widetilde{{\mathfrak T}}}$ that have poles
at the two preimages of the node. But a linear system of degree $4$ on
a genus two curve is always two dimensional and so the space of lines
in ${\mathbb P}^{2}$ is canonically isomorphic with
$|\omega_{\widetilde{{\mathfrak T}}}(J' + J'')|$. In other words, 
finding the bitangent lines to  $j(\widetilde{{\mathfrak T}})$ in
${\mathbb P}^{2}$ is equivalent to finding all divisors in
$|\omega_{\widetilde{{\mathfrak T}}}(J' + J'')|$ of the form
$2{\mathcal D}$ where ${\mathcal D}$ is an effective divisor of degree
two on $\widetilde{{\mathfrak T}}$. Since every degree two line bundle
on a genus two curve is effective we see that finding ${\mathfrak e}$
just amounts to choosing a non-trivial square 
root of the degree four line bundle $\omega_{\widetilde{{\mathfrak
T}}}(J' + J'')$.

\bigskip

Going back to the description of $B$ as the blow-up of ${\mathbb
P}^{2}$ at the base points of a pencil of cubics assume for
concreteness that  $J$ is the point in $\tilde{s}\cap
\widetilde{{\mathfrak T}}$ corresponding to the exceptional curve
$n_{1} \subset B$. Let ${\mathfrak e} \subset \widetilde{Q}$ be a
$(1,1)$ curve which passes trough $J$ and is bitangent to
$\widetilde{{\mathfrak T}}$ at two extra points. Let $e_{7} \subset B$
be one of the components of the preimage in $B$ of the strict
transform of $e_{7}$ in $B/(-1)_{B}$. Label by $e_{1}, \ldots, e_{6}$
the components of the reduced singular fibers of $\delta : B \to
{\mathbb P}^{1}$ which do not intersect $e_{7}$. Then $e_{1}, \ldots,
e_{6}$ and $e$ and $e_{7}$ are disjoint $(-1)$ curves on
$B$. After contracting these eight curves and the image of the curve $n_{1}$ 
we will get a ${\mathbb P}^{2}$.

Let $c : B \to {\mathbb P}^{2}$ denote this contraction map and let
$\ell = c^{*}{\mathcal O}_{{\mathbb P}^{2}}(1)$ be the pullback of the
class of a line via $c$. Thus $\op{Pic}(B)$ is generated over
${\mathbb Z}$ by the classes of the curves $\ell$, $e_{1}, \ldots, e_{6}$,
$e$, $e_{7}$ and $n_{1}$. In particular, if we put 
\[
\begin{split}
e_{9} & := e \\
e_{8} & := e + n_{1}
\end{split}
\]
we see that 
\[
H^{2}(B,{\mathbb Z}) = {\mathbb Z}\ell\oplus (\oplus_{i =
1}^{9}{\mathbb Z}e_{i}),
\]
with $\ell^{2} = 1$, $\ell\cdot e_{i} = 0$ and 
$e_{i}\cdot e_{j} = -\delta_{ij}$.

\bigskip

Note that in this basis we have 
\begin{equation} \label{eq-i2-fibers}
\begin{split}
n_{1} & = e_{8} - e_{9} \\
o_{1} & = f - e_{8} + e_{9} \\
n_{2} & = \ell - e_{7} - e_{8} - e_{9} \\
o_{2} & = 2\ell - e_{1} - e_{2} - e_{3} - e_{4} - e_{5} - e_{6}.
\end{split}
\end{equation}
\

\bigskip

\subsection{A synthetic construction} \label{sss-synthetic}
Before we proceed with the calculation of the action of $\tau_{B}$ on 
$H^{2}(B,{\mathbb Z})$ it will be helpful to analyze how the surface
$B$ and the map $c : B \to {\mathbb P}^{2}$ can be reconstructed
synthetically from geometric data on ${\mathbb P}^{2}$.

First we will need a general lemma describing a birational involution of
$\cp{2}$ fixing some smooth cubic pointwise.

\begin{lem} \label{lem-involution-p2} Let $\Gamma \subset \cp{2}$ be a smooth 
cubic and let $b \in \Gamma$. There exists a unique birational
involution  $\alpha : \cp{2} \dashrightarrow
\cp{2}$ which preserves the general line through $b$ and fixes the
general point of $\Gamma$. Let $b_{1}, b_{2}, b_{3}, b_{4} \in \Gamma$
be the four
ramification points for the linear projection of $\Gamma$ from $b$. Then
\begin{description}
\item[(i)] $\alpha$ sends a general line to a cubic which is nodal at $b$  and
passes through the $b_{i}$'s.
\item[(ii)] $\alpha$ sends the net of conics through $b_{1}, b_{2},
b_{3}$ to the net of cubics that are nodal at $b_{4}$ and pass through 
$b, b_{1}, b_{2}, b_{3}$.
\end{description}
\end{lem}
{\bf Proof.} Let $\alpha : \cp{2} \dashrightarrow \cp{2}$ be a birational
involution which fixes the general point of the cubic $\Gamma$ and
preserves the general line through $b \in \Gamma$. If $b
\in L \subset \cp{2}$ is a general line, then $L\cap \Gamma$ consists
of three distinct points $\{b, 0_{L}, \infty_{L} \}$. Since $\alpha$
preserves $L$ it follows that $\alpha_{|L}$ is a birational involution
of $L$ which fixes the points $0_{L}$ and $\infty_{L}$. But any
birational involution of $\cp{1}$ is biregular, has exactly two fixed
points and is uniquely determined by its fixed points. Thus the
restriction of $\alpha$ on the generic line through $b$ is uniquely
determined and so there can be at most one such $\alpha$. Conversely we can use
this uniqueness to show the existence of $\alpha$. Indeed, choose
coordinates $(x:y:z)$ in $\cp{2}$ so that $b = (0:0:1)$ and $\Gamma$ is
given by the equation $F(x,y,z) = 0$ with $F$ a homogeneous cubic
polynomial. Since $b \in \Gamma$ we can write $F = F_{1}z^{2} + F_{2}z
+ F_{3}$ with $F_{d}$ a homogeneous polynomial in $(x,y)$ of degree
$d$. Let $(x:y:z)$ be a point in $\cp{2}$ and let $L = \{(
x: y: z+ t)\}_{t \in \cp{1}}$ be the line
through $b$ and $(x:y:z)$.  The involution $\alpha_{|L}$ will have to
fix the two additional (besides $b$) intersection points of $L$ and
$\Gamma$. The values of $t$ corresponding to these points are just the
roots of the equation $F(x,y,z+t) = 0$, that is the solutions to
\begin{equation} \label{eq-0infty1}
F_{1}(x,y)t^{2} + F_{z}(x,y,z)t + F(x,y,z) = 0.
\end{equation} 
On the other hand since $t$ is the affine coordinate on $L$ the
involution $\alpha_{|L} : \cp{1} \to \cp{1}$ will be given by a
fractional linear transformation
\[
t \mapsto \frac{at +b}{ct + d}
\]
for some complex numbers $a$, $b$, $c$ and $d$. The condition that 
$\alpha_{|L} \neq \op{id}_{L}$ but $\alpha_{|L}^{2} = \op{id}_{L}$ 
is equivalent to $d = -a$.

In these terms the fixed points of $\alpha_{|L}$ correspond to the
values of $t$ for which 
\begin{equation} \label{eq-0infty2}
ct^{2} - 2a t - b = 0.
\end{equation} 
Comparing \eqref{eq-0infty1} with \eqref{eq-0infty2} we conclude that
$a = -(1/2)F_{z}(x,y,z)$, $b = -F(x,y,z)$ and $c = F_{1}(x,y)$ and so
\[
\alpha_{|L}((x:y:z+t)) = \left( x : y : z - \frac{F_{z}(x,y,z)t +
2F(x,y,z)}{2F_{1}(x,y)t + F_{z}(x,y,z)} \right).
\]
In particular for $t=0$ we must have
\begin{equation} \label{eq-alpha-formula}
\alpha((x:y:z)) = \alpha_{|L}((x:y:z)) = \left( x : y : z - 2
\frac{F(x,y,z)}{F_{z}(x,y,z)} \right).
\end{equation}
Now the formula \eqref{eq-alpha-formula} clearly defines a birational
automorphism $\alpha$ of $\cp{2}$ and it is straightforward to check that 
$\alpha^{2} = \op{id}_{\cp{2}}$. This shows the existence and
uniqueness of $\alpha$.

To prove the remaining statements note that the $\alpha$ that we have
just defined lifts to a biregular involution $\hat{\alpha}$ on the
blow-up $g : \widehat{\cp{2}} \to \cp{2}$ of
$\cp{2}$ at the points $b, b_{1}, b_{2}, b_{3}, b_{4}$. Let $\Sigma,
\Sigma_{1}, \Sigma_{2}, \Sigma_{3}, \Sigma_{4} \subset
\widehat{\cp{2}}$ denote the corresponding exceptional divisors and
let $\ell = g^{*}{\mathcal O}_{\cp{2}}(1)$ be the class of a line. By
definition $\alpha$ preserves the general line through $b$ and the
cubic $\Gamma$. Hence $\hat{\alpha}$ will preserve the proper
transforms of $\Gamma$ and the general line through $b$, i.e.
\[
\begin{split}
\hat{\alpha}(\ell - \Sigma) & = \ell - \Sigma \\
\hat{\alpha}\left(3\ell - \Sigma - \sum_{i = 1}^{4} \Sigma_{i}\right) & 
= 3\ell - \Sigma - \sum_{i = 1}^{4} \Sigma_{i}.
\end{split}
\]
Also it is clear (e.g. from \eqref{eq-alpha-formula}) that
$\hat{\alpha}$ identifies the proper transform of the
line through $b$ and $b_{i}$ with $\Sigma_{i}$ and so
\[
\hat{\alpha}(\Sigma_{i}) = \ell - \Sigma - \Sigma_{i}
\]
for $i = 1,2,3,4$. Therefore we get two equations for
$\hat{\alpha}(\ell)$ and $\hat{\alpha}(\Sigma)$:
\[
\begin{split}
\hat{\alpha}(\ell) - \hat{\alpha}(\Sigma) & = \ell - \Sigma \\
3\hat{\alpha}(\ell) - \hat{\alpha}(\Sigma) & = 7\ell - 5\Sigma - 2
\sum_{i = 1}^{4} \Sigma_{i},
\end{split}
\]
which yield $\hat{\alpha}(\ell) = 3\ell - 2\Sigma - \sum_{i = 1}^{4}
\Sigma_{i}$ and $\hat{\alpha}(\Sigma) = 2\ell - \Sigma - \sum_{i =
1}^{4} \Sigma_{i}$. 

If now $L$ is a line not passing through any of the points $b, b_{1},
b_{2}, b_{3}, b_{4}$ we see that the proper transform $\widehat{L}$ of
$L$ in $\widehat{\cp{2}}$ is an irreducible curve such that 
$\hat{\alpha}(\widehat{L})$ is in the linear system
$|3\ell - 2\Sigma - \sum_{i = 1}^{4} \Sigma_{i}|$.  In particular
$\hat{\alpha}(\widehat{L})$ intersects $\Sigma$ at two points and
intersects each $\Sigma_{i}$ at a point. So $\alpha(L) =
g(\hat{\alpha}(\widehat{L}))$ is a cubic which is nodal at $b$ and
passes through each of the $b_{i}$'s. This proves part (i) of the
lemma.

Similarly if $C$ is a conic through $b_{1}$, $b_{2}$ and $b_{3}$, then 
$\widehat{C}$ is an irreducible curve in the linear system $|2\ell -
\Sigma_{1} - \Sigma_{2} - \Sigma_{3}|$ on $\widehat{\cp{2}}$. Hence 
$\hat{\alpha}(\widehat{C})$ is an irreducible curve in the linear
system $|3\ell - \Sigma - \Sigma_{1} - \Sigma_{2} - \Sigma_{3} -
2\Sigma_{4}|$ and so $\alpha(C) = g(\hat{\alpha}(\widehat{C}))$ is a
cubic passing through $b, b_{1}, b_{2}, b_{3}$ which is nodal at
$b_{4}$. The lemma is proven. \hfill $\Box$
\

\bigskip

\begin{figure}[!ht]
\begin{center}
\epsfig{file=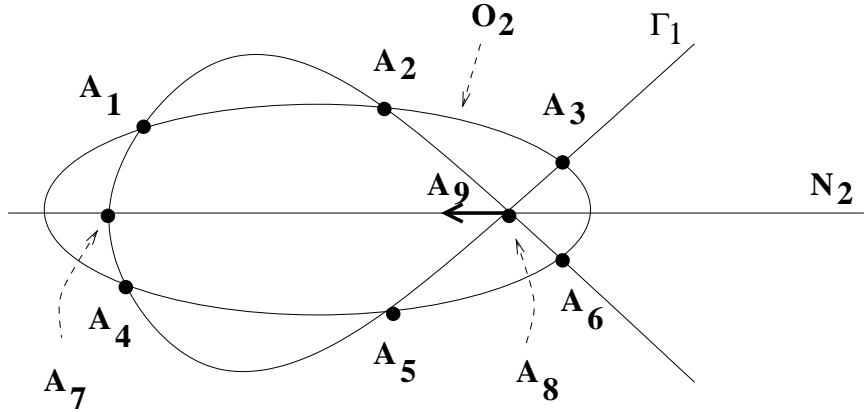,width=4.5in} 
\end{center}
\caption{The pencil of cubics determining $B$}\label{fig-synthetic} 
\end{figure}

For our synthetic construction of $B$ we will start with a nodal cubic 
$\Gamma_{1} \subset \cp{2}$ and will denote its node by $A_{8} \in
\Gamma_{1}$. Pick
four other points on $\Gamma_{1}$ and label them $A_{1}, A_{2}, A_{3},
A_{7}$. For generic such choices there is a unique smooth cubic
$\Gamma$ which passes through the points $A_{1}, A_{2}, A_{3}, A_{7},
A_{8}$ and is tangent to the line $\langle A_{7} A_{i} \rangle$ at the
point $A_{i}$ for $i = 1,2,3$ and $8$. Consider the pencil of cubics 
spanned by $\Gamma$  and $\Gamma_{1}$. All cubics in this pencil pass
through $A_{1}, A_{2}, A_{3}, A_{7}, A_{8}$ and are tangent to
$\Gamma$ at $A_{8}$. Let $A_{4}, A_{5}, A_{6}$ be
the remaining three base points. Each cubic in the pencil 
intersects the line $N_{2} := \langle
A_{7}A_{8} \rangle$ in the same divisor $A_{7} + 2 A_{8} \in
\op{Div}(N_{2})$.  Therefore there is a reducible cubic $\Gamma_{2} =
N_{2} \cup O_{2}$ in the pencil. Generically $O_{2}$ will be a smooth
conic as depicted on Figure~\ref{fig-synthetic}. 
By 
Lemma~\ref{lem-involution-p2} there is a birational
involution $\alpha$ of $\cp{2}$ corresponding to $\Gamma$ with $b =
A_{7}$. Note that by construction 
$b_{i} = A_{i}$ for $i = 1,2,3$ and $b_{4} = A_{8}$. By
Lemma~\ref{lem-involution-p2}(ii) we know that $\alpha(O_{2})$ is a
nodal cubic with a node at $A_{8}$ which passes through $A_{1}, A_{2},
A_{3}$ and $A_{7}$. Since the involution
$\alpha$ fixes $A_{4}, A_{5}, A_{6} \in \Gamma$ it also follows that
$\alpha(O_{2})$ contains $A_{4}, A_{5}, A_{6}$. The intersection
number $\alpha(O_{2})$ with $\Gamma_{1}$ is therefore at least $6 +
2\cdot 2 = 10$ and so $\alpha(O_{2}) = \Gamma_{1}$. Moreover $\alpha$
collapses $N_{2}$ to $A_{8}$. This shows that 
$\alpha$ preserves the pencil.

We define $B$ to be the blow-up of $\cp{2}$ at the points $A_{i}$, $i
=1, \ldots, 8$ and the point $A_{9}$ which is  infinitesimally near to
$A_{8}$ and corresponds to the tangent direction $N_{2}$. The pencil
of cubics becomes the anticanonical map $\beta : B \to {\mathbb P}^{1}$.
The reducible fibers are $f_{i} = n_{i}\cup o_{i}$, $i = 1,2$ where
$n_{2}, o_{2}$ are the proper transforms of $N_{2}, O_{2}$, $o_{1}$ is
the proper transform of $\Gamma_{1}$ and $n_{1}$ is the proper
transform of the exceptional
divisor corresponding to $A_{8}$. In order to conform with the
notation in Section~\ref{ss-dp9-general} we denote by $e_{i}$ 
for $i = 1, \ldots, 7$
and $9$ 
the exceptional divisors corresponding to $A_{i}$, $i = 1, \ldots,
7$ and $9$ and by  $e_{8}$ the reducible divisor $e_{9} + n_{1}$. 

The involution $\alpha : \cp{2} \dashrightarrow \cp{2}$ lifts to a
biregular involution $\alpha_{B} : B \to B$. The induced involution
$\tau_{\cp{1}}$ of $\cp{1}$ has two fixed points $0, \infty \in
\cp{1}$. One of them, say $0$, will be the image $\beta(\Gamma)$. We
will use $e_{9}$ as the zero section $e : \cp{1} \to B$. Note that 
$(-1)_{B}^{*}e_{i} = \alpha_{B}^{*}e_{i}$ for $i = 1, 2, 3$ and so we
can take $\zeta = e_{1}$.

\section{Action on cohomology} \label{ss-dp9-cohoaction}

First we describe  the action of the automorphisms
$(-1)_{B}$, $\alpha_{B}$, $t_{\zeta}$ and $\tau_{B}$ on
$H^{\bullet}(B,{\mathbb Z})$.

\subsection{Action of $(-1)_{B}$} 
From the discussion in section~\ref{subsec-case2-basis-in-h2} it is
clear that $(-1)_{B}$ preserves the fibers of $\delta : B \to \cp{1}$ and
exchanges the two components of the six singular 
fibers of $\delta$ which are unions of two rational
curves meeting at a point. Furthermore from the description of $B$ as
a blow-up of $\cp{2}$ at nine points (see
section~\ref{subsec-case2-basis-in-h2}) it follows that the class of
the fiber of $\delta$ is $\ell - e_{7}$. Hence $(-1)_{B}(\ell -e_{7})
= \ell -e_{7}$ and $(-1)_{B}(e_{i}) +
e_{i} = \ell - e_{7}$ for $i = 1, \ldots, 6$. Also, by the same
analysis we see that $(-1)_{B}$ preserves $n_{1}$ and $n_{2}$ and
since $(-1)_{B}$ preserves $f$ by definition, it follows that
$(-1)_{B}$ preserves $o_{1}$ and $o_{2}$ as well. Similarly
$(-1)_{B}$ preserves $e_{9}$ by definition and so $(-1)_{B}^{*}(e_{8}) =
(-1)_{B}^{*}(e_{9} + n_{1}) = e_{9} + n_{1} = e_{8}$. Finally we can solve 
the equations $(-1)_{B}^{*}(\ell -e_{7}) = \ell -e_{7}$ and
$(-1)_{B}^{*}(o_{2}) = o_{2}$ to get $(-1)_{B}(\ell) = f + \ell - 2e_{7} +
e_{8} + e_{9}$ and $(-1)_{B}^{*}(e_{7}) = f - e_{7} + e_{8} + e_{9}$.

\subsection{Action of $\alpha_{B}$} Again from the analysis in
section~\ref{subsec-case2-basis-in-h2} and the geometric description
of $B/\alpha_{B}$ and its Weierstrass model $W_{M}$ we see that
$\alpha_{B}$ preserves the classes of the fibers of the two fibrations
$\beta : B \to \cp{1}$ and $\delta : B \to \cp{1}$. In particular we
have $\alpha_{B}^{*}(f) = f$, $\alpha_{B}^{*}(\ell - e_{7}) =  \ell -
e_{7}$ and $\alpha_{B}^{*}(e_{9}) = e_{9}$. Also $\alpha_{B}$
interchanges $o_{1}$ and $o_{2}$ and hence interchanges $n_{1}$ and
$n_{2}$. From the relationship between the ramification divisors
defining $W_{M}$ and $\widetilde{Q}$ we see that $\alpha_{B}$ will
exchange the two components of the three singular fibers  of $\delta$
corresponding to the three intersection points in ${\mathfrak T}\cap
r_{0}$, i.e. $\alpha_{B}^{*}(e_{j}) + e_{j} = \ell -e_{7}$ for $j =
1,2,3$. Similarly $\alpha_{B}$ will preserve the two components of the
singular fibers of $\delta$ corresponding to the three intersection
points in ${\mathfrak T}\cap r_{\infty}$, that is
$\alpha_{B}^{*}(e_{i}) = e_{i}$ for $i = 4,5,6$. Finally, solving 
the equations $\alpha_{B}^{*}(\ell - e_{7}) =  \ell - e_{7}$ and
$\alpha_{B}^{*}(o_{1}) = o_{2}$ we get $\alpha_{B}^{*}(\ell) = 3\ell -
e_{1} - e_{2} - e_{3} - 2e_{7} - e_{8}$ and $\alpha_{B}^{*}(e_{7}) = 2\ell -
e_{1} - e_{2} - e_{3} - e_{7} - e_{8}$.

\subsection{Action of $t_{\zeta}^{*}$} By definition we have
$t_{\zeta}^{*}(f) = f$. In order to find the action of
$t_{\zeta}$ on the classes $e_{i}$ we will use the fact that
$t_{\zeta}$ is defined in terms of the addition law on $\beta^{\#}
: B^{\#} \to \cp{1}$. 

Since $t_{\zeta}$
preserves each fiber of $\beta : B \to \cp{1}$, the curve
$t_{\zeta}^{*}(n_{1})$ will have to be either $n_{1}$ or $o_{1}$. But
$\zeta = e_{1}$ and so $\zeta\cdot n_{1} = 0$ and $\zeta\cdot o_{1} =
1$, so since $n_{1}^{\#}$ is the identity component of the
disconnected group $n_{1}^{\#}\cup o_{1}^{\#} = (n_{1}\cup
o_{1})- (n_{1}\cap o_{1})$,  we must have 
$t_{\zeta}^{*}(n_{1}) = o_{1}$. In the same way one can argue that
$t_{\zeta}^{*}(n_{2}) = o_{2}$ and $t_{\zeta}^{*}(o_{i}) = n_{i}$ for
$i = 1, 2$.

Next note that since $t_{\zeta}$ is compatible with the group scheme
structure of $B^{\#}$ we must have $t_{\zeta}^{*}(\xi) =
c_{1}([\xi] - [\zeta])$ for any section $\xi$ of $\beta$. Using this
relation we calculate: 
\[
\begin{split}
t_{\zeta}^{*}(e_{1}) & = c_{1}([e_{1}] - [e_{1}]) = e_{9}, \\
t_{\zeta}^{*}(e_{9}) & = c_{1}([e_{9}] - [e_{1}]) = (-1)_{B}([e_{1}]) =
\ell - e_{1} - e_{7},
\end{split}
\]
which in turn implies $t_{\zeta}^{*}(e_{8}) = t_{\zeta}^{*}(e_{9} +
n_{1}) = \ell - e_{1} - e_{7} + o_{1} = f + \ell - e_{1} - e_{7} -
e_{8} + e_{9}$.

The previous formulas identify cohomology classes in $H^{2}(B,{\mathbb
Z})$ or equivalently line bundles on $B$. However
observe that the above formulas can also be viewed as equality of
divisors, due to the fact that the line bundles in question correspond
to sections of $\beta$, and so each of these is represented by a unique
(rigid) effective divisor.

Also since the addition law on an elliptic curve is defined in terms
of the Abel-Jacobi map we see that for a section $\xi$ of $\beta$, the
restriction of the line bundle $c_{1}([\xi] - [e_{1}])\otimes
{\mathcal O}_{B}(-e_{9})$  to the generic fiber of $\beta$ will be the
same as the restriction of ${\mathcal O}_{B}(\xi - e_{1})$. By the
see-saw principle the difference of these two line bundles will have
to be a combination of components of fibers of $\beta$, i.e.
\[
t_{\zeta}^{*}(\xi) = c_{1}([\xi] - [e_{1}]) = \xi - e_{1} + e_{9} +
a_{1}^{\xi}n_{1} + a_{2}^{\xi}n_{2} + a^{\xi}f.
\]
Intersecting both sides with $n_{1}$ and taking into account that
$(t_{\zeta}^{-1})^{*}(n_{1}) = o_{1}$ we get $o_{1}\cdot \xi =
\xi\cdot n_{1} + 1 - 2 a_{1}^{\xi}$. Similarly when we intersect with
$n_{2}$ we get $o_{2}\cdot \xi = \xi\cdot n_{2} + 1 - 2 a_{2}^{\xi}$. 
In particular since for $i = 2, \ldots, 6$ we have $e_{i}\cdot n_{1} =
e_{i}\cdot n_{2} = 0$ and  $e_{i}\cdot o_{1} =
e_{i}\cdot o_{2} = 1$ we get $a_{1}^{e_{i}} = a_{2}^{e_{i}} = 0$ and so
$t_{\zeta}^{*}(e_{i}) = e_{i} - e_{1} + e_{9} + a^{e_{i}}f$.
Using the fact that $(t_{\zeta}^{*}(e_{i}))^{2} = -1$ we find that
$a^{e_{i}} = 1$ and thus 
\[
t_{\zeta}^{*}(e_{i}) = e_{i} - e_{1} + e_{9} + f
\]
for $i = 2, \ldots, 6$. 

Finally, for $e_{7}$ we have $e_{7}\cdot n_{1} = e_{7}\cdot o_{2} = 0$
and $e_{7}\cdot n_{2} = e_{7}\cdot o_{1} = 1$ and so 
$t_{\zeta}^{*}(e_{7}) = e_{7} - e_{1} + e_{9} + n_{2} +
a^{e_{7}}f$. From $(t_{\zeta}^{*}(e_{7}))^{2} = -1$ we find $a^{e_{7}}
= 0$ and therefore $t_{\zeta}^{*}(e_{7}) = e_{7} - e_{1} + e_{9} +
n_{2}$. 

This completes the calculation of the action of $t_{\zeta}^{*}$ on
$H^{2}(B,{\mathbb Z})$. The action of $\tau_{B}^{*}$ is easily
obtained since by definition we have $\tau_{B}^{*} =
\alpha_{B}^{*}\circ t_{\zeta}^{*}$. 

All these actions are summarized in Table~\ref{table-auts} below.

\begin{table}[!ht]
\begin{center}
\begin{tabular}{|l||c|c|c|c|} \hline
 & $(-1)_{B}^{*}$ & $t_{\zeta}^{*}$ & $\alpha_{B}^{*}$ &
 $\tau_{B}^{*}$ \\ \hline \hline
$f$ & $f$ & $f$ & $f$ & $f$ \\ \hline
$e_{1}$ & $\ell - e_{1} - e_{7}$ & $e_{9}$ & $\ell - e_{1} - e_{7}$ &
 $e_{9}$ \\ \hline
$e_{j}$, & $\ell - e_{j} - e_{7}$ & $f + e_{j} - e_{1} +
 e_{9}$ & $\ell - e_{j} - e_{7}$ & $f - e_{j} + e_{1} + e_{9}$ \\
$j = 2, 3$ & &  & & \\ \hline
$e_{i}$,  & $\ell - e_{i} - e_{7}$ & $f + e_{i} - e_{1} +
 e_{9}$ & $e_{i}$ & $f - \ell + e_{i} + $ \\
$i = 4, 5, 6$ & & & & $+ e_{1} + e_{7} + e_{9}$ \\
 \hline
$e_{7}$ & $f - e_{7} + e_{8} + e_{9}$ & $\ell - e_{1} - e_{8}$ &
 $2\ell - (e_{1} + e_{2} + $ & $\ell - e_{2} -
 e_{3}$ \\ 
& & & $ + e_{3} + e_{7} + e_{8})$ & \\ \hline
$e_{8}$ & $e_{8}$ & $f + \ell + e_{9} - $ & $\ell
 - e_{7} - e_{8}$ & $f - \ell + e_{1} + $ \\
&  & $ - e_{1} - e_{7} - e_{8}$ & & $+ e_{7} + e_{8} + e_{9}$ \\
\hline
$e_{9}$ & $e_{9}$ & $\ell - e_{1} - e_{7}$ & $e_{9}$ & $e_{1}$ \\
 \hline \hline
$\ell$ & $\ell + f -$ & $2f + 2\ell - 3 e_{1} - $ & $3\ell - ( e_{1} +
 e_{2} +$ & $2f + 2(e_{1} + e_{9}) -$ \\ 
& $- 2e_{7} + e_{8} + e_{9}$ & $-
 e_{7} - e_{8} + 2e_{9}$ & $ + e_{3} + 2e_{7} +
 e_{8})$ & $ - (e_{2} + e_{3}) + e_{7}$ \\
\hline 
\end{tabular}
\end{center}
\caption{Action of $(-1)_{B}$, $\alpha_{B}$, $t_{\zeta}$ and $\tau_{B}$ on
$H^{\bullet}(B,{\mathbb Z})$}\label{table-auts}
\end{table}

\

\section{The cohomological Fourier-Mukai transform} \label{s-cohoFM} 
For the purposes of the spectral construction we will need also the
action of the relative Fourier-Mukai transform for $\beta : B \to
\cp{1}$ on the cohomology of $B$. By definition the Fourier-Mukai
transform is the exact functor on the bounded derived category $D^{b}(B)$ 
of $B$ given by the formula
\[
\xymatrix@R=4pt{
\FM_{B} : & D^{b}(B) \ar[r] & D^{b}(B) \\
& {\mathcal F} \ar@{|->}[r] &
R^{\bullet}p_{1*}(p_{2}^{*}{\mathcal
F}\stackrel{L}{\otimes} {\mathcal P}_{B}).
}
\]
Here $p_1,p_2$ are the projections of $B \times_{\cp{1}}B$ to its two
factors, and ${\mathcal P}_B$ is the Poincare sheaf:
\[
{\mathcal P}_B:= {\mathcal O}_B(\Delta-e \times_{\cp{1}}B -
B\times_{\cp{1}}e-q^*{\mathcal O}_{\cp{1}}(1)),
\]
with $q = \beta\circ p_{1} = \beta\circ p_{2}$. Using the zero section
$e : \cp{1} \to B$ we can identify $B$ with the
relative moduli space ${\mathcal M}(B/\cp{1})$ of semistable
(w.r.t. to a suitable polarization), rank one, degree zero  torsion
free sheaves along the fibers of $\beta : B \to \cp{1}$. Under this
identification, the sheaf ${\mathcal P}_{B} \to B\times_{\cp{1}} B =
B\times_{\cp{1}} {\mathcal M}(B/\cp{1})$ becomes the universal
sheaf. This puts us in the setting of
\cite[Theorem~1.2]{bridgeland-cy} and implies that $\FM_{B}$ is an
{\em autoequivalence} of $D^{b}(B)$. In particular we can view any vector
bundle $V \to B$ in two different ways - as $V$ and as the object
$\FM_{B}(V) \in D^{b}(B)$. 

\medskip

The cohomological Fourier-Mukai transform is defined as the unique
linear map 
\[
\fm_{B} : H^{\bullet}(B,{\mathbb Q}) \to
H^{\bullet}(B,{\mathbb Q})
\] 
satisfying:
\begin{equation} \label{eq-fm-vs-FM}
\fm_{B} \circ ch = ch \circ \FM_{B}.
\end{equation} 
Explicitly, 
\[
\fm_{B}(x) = \op{pr}_{2*}(\op{pr}_1^*(x)\cdot ch(j_{*}{\mathcal P})
\cdot td(B \times B)) \cdot 
td(B)^{-1},
\]
where $\op{pr}_i$ are the projections of $B \times B$ to its factors
and $j : B\times_{{\mathbb P}^{1}} B \hookrightarrow B\times B$ is the
natural inclusion.

We will need an explicit description of the cohomological spectral involution 
\[
\ct_{B} := \fm_{B}^{-1}\circ
\tau_{B}^{*} \circ \fm_{B}.
\] 
For this we proceed to calculate the action of $\fm_{B}$ and
$\fm_{B}^{-1}$ in the obvious basis in cohomology. 

Let $\op{pt} \in H^{4}(B,{\mathbb Z})$ denote the
class Poincare dual to the homology class of a point in $B$  and let
$1 \in H^{0}(B,{\mathbb Z})$ be the class which is Poincare dual to
the fundamental class of $B$. The classes $1$, $f$, $e_{1}$, \ldots,
$e_{9}$, $\op{pt}$ constitute a basis of
$H^{\bullet}(B,{\mathbb Q})$. 

To calculate $\fm_{B}$ we will use the identity \eqref{eq-fm-vs-FM}
together with a calculation of the action of $\FM_{B}$ on certain
basic sheaves, which is carried out in Lemma~\ref{lem-fm-sheaves} below.

The first observation is that there are two ways to lift a sheaf 
$G$ on ${\mathbb P}^{1}$ to a sheaf on $B$. First we may consider the
pullback $\beta^{*}(G)$. Second, for any section $\xi : {\mathbb
P}^{1} \to B$ of $\beta$ we may form the push-forward
$\xi_{*}G$. These two lifts behave quite differently. For example, if
$G$ is a line bundle, then $\beta^{*}G$ is a line bundle on $B$,
whereas $\xi_{*}G$ is a torsion sheaf on $B$ supported on $\xi$.
The action of $\FM_{B}$ interchanges these two types of sheaves (up to a
shift):

\begin{lem} \label{lem-fm-sheaves} 
For any sheaf $G$ on $\cp{1}$ and any section $\xi$ of $\beta$ we have:
\[
\begin{split}
\FM_{B}(\beta^*G) & =e_*(G\otimes{\mathcal O}_{\cp{1}}(-1))[-1] \\
\FM_{B}(\xi_*G) & =\beta^*G\otimes{\mathcal O}_B(\xi -e)
\otimes \beta^*{\mathcal O}_{\cp{1}}(-e\cdot\xi-1),
\end{split}
\]
where as usual for a complex $K^{\bullet} = (K^{i},d_{K}^{i})$ and an
integer $n \in {\mathbb Z}$ we put $K^{\bullet}[n]$ for the complex
having $(K[n])^{i} = K^{n+i}$ and $d_{K[n]} = (-1)^{n}d_{K}$. 
\end{lem}
{\bf Proof.} By definition we have $\FM_{B}(\beta^*G) =
Rp_{2*}(p_{1}^{*}\beta^{*}G\otimes {\mathbb P}_{B})$. But $\beta\circ
p_{1} = \beta\circ p_{2}$ and so by the projection formula we get 
$\FM_{B}(\beta^*G) = Rp_{2*}(p_{2}^{*}\beta^{*}G\otimes {\mathbb
P}_{B}) = \beta^{*}G\otimes Rp_{2*}{\mathcal P}_{B}$. In order to calculate 
$Rp_{2*}{\mathcal P}_{B}$, note first that $Rp_{2*}{\mathcal P}_{B}$ is
a complex concentrated in degrees zero and one since $p_{2}$ is a
morphism of relative dimension one. Next observe 
that $R^{0}p_{2*}{\mathcal P}_{B} =
0$.  Indeed, by definition ${\mathcal P}_{B}$ is a rank one torsion
free sheaf on $B\times_{\cp{1}} B$, and so $R^{0}p_{2*}{\mathcal
P}_{B}$ must be a torsion free sheaf on $B$. On the other hand, from
the definition of ${\mathcal P}_{B}$ we see that both $R^{0}p_{2*}{\mathcal
P}_{B}$ and $R^{1}p_{2*}{\mathcal P}_{B}$ are torsion sheaves on $B$
whose reduced support is precisely $e \subset B$. Therefore
$R^{0}p_{2*}{\mathcal P}_{B}$ is torsion and torsion free at the same
time and so $R^{0}p_{2*}{\mathcal P}_{B} = 0$ \label{torsion}. 
This implies that
$Rp_{2*}{\mathcal P}_{B} = R^{1}p_{2*}{\mathcal P}_{B}[-1]$. Now,
since $R^{2}p_{2*}{\mathcal P}_{B} = 0$ we can apply the cohomology
and base change theorem \cite[Theorem~12.11]{hartshorne} to conclude
that $R^{1}p_{2*}{\mathcal P}_{B}$ has the base change property for
arbitrary (i.e. not necessarily flat) morphisms. In particular
considering the base change diagram
\[
\xymatrix{
B = B\times_{\cp{1}} e \ar@{^{(}->}[r] \ar[d]_{\beta} & B\times_{\cp{1}} B
\ar[d]^{p_{2}} \\
\cp{1} \ar@{^{(}->}[r]_-{e} & B
}
\] 
we have that 
\[
\begin{split}
e^{*}R^{1}p_{2*}{\mathcal P}_{B} & =
R^{1}\beta_{*}({\mathcal P}_{B|B\times_{\cp{1}} e}) =
R^{1}\beta_{*}{\mathcal O}_{B} = (\beta_{*}\omega_{B/\cp{1}})^{\vee}
\\
&  = (\beta_{*}({\mathcal O}_{B}(-f)\otimes \beta^{*}{\mathcal
O}(2)))^{\vee} = {\mathcal O}_{\cp{1}}(-1).
\end{split}
\]
Since $e \subset B$ is the reduced support of $R^{1}p_{2*}{\mathcal
P}_{B}$ and $(R^{1}p_{2*}{\mathcal P}_{B})_{|e}$ is a line bundle, it
follows that $e\subset B$ is actually the scheme theoretic support of 
$R^{1}p_{2*}{\mathcal P}_{B}$ and so $R^{1}p_{2*}{\mathcal P}_{B} =
e_{*}{\mathcal O}_{\cp{1}}(-1)$, which finishes the proof of the first
part of the lemma.

Let now $\xi : \cp{1} \to B$ be a section of $\beta$. Then
$\FM_{B}(\xi_{*}G) = Rp_{2*}(p_{1}^{*}\xi_{*}G\otimes {\mathcal
P}_{B})$.  But $p_{1}^{*}\xi_{*}G$ is a sheaf on $B\times_{\cp{1}} B$ 
supported on
$\xi\times_{\cp{1}} B \subset B\times_{\cp{1}} B$ and is in fact the
extension by zero of the sheaf $\beta^{*}G$ on $B = \xi\times_{\cp{1}}
B$. Moreover by definition we have ${\mathbb P}_{B|\xi\times_{\cp{1}}
B} = {\mathcal O}_{B}(\xi - e - (e\cdot \xi + 1)f)$. Taking into
account that $p_{2} : \xi\times_{\cp{1}} B \to B$ is an isomorphism, we
get the second statement of the lemma. \hfill $\Box$

\bigskip
\bigskip

With all of this said we are now ready to derive the explicit formulas
for $\fm_{B}$. First, observe that $ch({\mathcal O}_{B}) = 1$ and so
by \eqref{eq-fm-vs-FM}  and Lemma~\ref{lem-fm-sheaves} we have
\[
\begin{split}
\fm_{B}(1) & = ch(\FM_{B}({\mathcal O}_{B})) =
ch(\FM_{B}(\beta^{*}{\mathcal O}_{\cp{1}})) \\
& = ch(e_{*}({\mathcal O}_{\cp{1}}(-1))[-1]) = - ch(e_{*}({\mathcal
O}_{\cp{1}}(-1)). 
\end{split}
\]
But from the short exact sequence of sheaves on $B$
\[
0 \to {\mathcal O}_{B}(-e - f) \to {\mathcal O}_{B}(-f) \to
e_{*}{\mathcal O}_{\cp{1}}(-1) \to 0
\] 
we calculate
\[
\begin{split}
ch(e_{*}({\mathcal O}_{\cp{1}}(-1)) & = ch({\mathcal O}_{B}(-f)) -
ch({\mathcal O}_{B}(-e - f)) \\
& = (1 - f + 0\cdot \op{pt})\cdot \left( 1 + (e - f) +
\frac{1}{2}\op{pt}\right) \\
& = e - \frac{1}{2}\op{pt}.
\end{split}
\]
In other words $\fm_{B}(1) = - e + (1/2)\op{pt} = -e_{9} + (1/2)\op{pt}$.

Next we calculate $\fm_{B}(\op{pt})$. Let $t \in \cp{1}$ be a fixed
point. Then $\op{pt} = ch({\mathcal O}_{e(t)}) = ch(e_{*}{\mathcal
O}_{t})$ and so
\[
\begin{split}
\fm_{B}(\op{pt}) & = ch(\FM_{B}(e_{*}{\mathcal O}_{t})) \\
& = ch({\mathcal O}_{f}) = ch({\mathcal O}_{B}) - ch({\mathcal
O}_{B}(-f)) \\
& = 1 - (1 - f + 0\cdot \op{pt}) = f.
\end{split}
\]
To calculate $\fm_{B}(f)$ note that $ch({\mathcal O}_{B}(f)) = 1 + f$
and so 
\[
\begin{split}
\fm_{B}(f) & = ch(\FM_{B}({\mathcal O}_{B}(f))) - \fm_{B}(1) \\
& = ch(\FM_{B}(\beta^{*}{\mathcal O}_{\cp{1}}(1))) - \fm_{B}(1) \\
& = ch(e_{*}{\mathcal O}_{\cp{1}}[-1]) - \left(-e +
\frac{1}{2}\op{pt}\right) \\
& = - [ch({\mathcal O}_{B}) - ch({\mathcal O}_{B}(-e))] + e -
\frac{1}{2}\op{pt} \\ 
& = - \left[ 1 - \left( 1 -  e - \frac{1}{2}\op{pt}\right)\right] + e -
\frac{1}{2}\op{pt} \\  
& = - \op{pt}.
\end{split}
\]
Finally we calculate $\fm_{B}(e_{i})$.  If $i = 1, \ldots, 7$, the
class $e_{i}$ is a class of a section $e_{i} : \cp{1} \to B$ of
$\beta$ and so we can apply Lemma~\ref{lem-fm-sheaves} to ${\mathcal
O}_{e_{i}}$. We have $ch({\mathcal O}_{e_{i}}) = e_{i} + (1/2)\op{pt}$
and hence
\[
\begin{split}
\fm_{B}(e_{i}) & = ch(\FM_{B}({\mathcal O}_{e_{i}})) -
\frac{1}{2}\fm_{B}(\op{pt}) \\
& = ch(\FM_{B}(e_{i*}{\mathcal O}_{\cp{1}})) -
\frac{1}{2}\fm_{B}(\op{pt}) \\
& = ch({\mathcal O}_{B}(e_{i} - e_{9} - f)) - \frac{1}{2}f \\
& = 1 + (e_{i} - e_{9} - f) - \op{pt} - \frac{1}{2}f \\
& = 1 + (e_{i} - e_{9} -\frac{3}{2}f) - \op{pt}.
\end{split}
\]
For $e_{9}$ we get in the same way
\[
\fm_{B}(e_{9}) = ch({\mathcal O}_{B}) - \frac{1}{2}f = 1 - \frac{1}{2}f,
\]
and so it only remains to calculate $\fm_{B}(e_{8})$. 

Unfortunately we can not use the same method for calculating
$\fm_{B}(e_{8})$ since $e_{8}$ is only a numerical section of $\beta$
and splits as a union of two irreducible curves $e_{8} = e_{9} +
n_{1}$. However, recall that the automorphism $\alpha_{B} : B \to B$
moves a section to a section. Consequently $\alpha_{B}(e_{7})$ will be
another section of $\beta$. Let $a : \cp{1} \to B$ denote the
map corresponding to $\alpha_{B}(e_{7})$. Then
\[
ch({\mathcal O}_{\alpha_{B}(e_{7})}) = ch({\mathcal O}_{B}) -
ch({\mathcal O}_{B}(-\alpha_{B}(e_{7})) = \alpha_{B}(e_{7}) +
\frac{1}{2}\op{pt}. 
\]
Thus 
\[
\fm_{B}(\alpha_{B}(e_{7})) = ch(\FM_{B}(a_{*}{\mathcal O}_{\cp{1}})) -
\frac{1}{2} f = ch({\mathcal O}_{B}(\alpha_{B}(e_{7}) - e_{9} -
(e_{9}\cdot \alpha_{B}(e_{7}) + 1)f) - \frac{1}{2} f.
\]
But according to Table~\ref{table-auts} we have 
$e_{9}\cdot \alpha_{B}(e_{7}) = e_{9}\cdot (2\ell -e_{1} - e_{2} -
e_{3} - e_{7} - e_{8}) =
0$ and so
\[
\fm_{B}(\alpha_{B}(e_{7})) = 1 + \alpha_{B}(e_{7}) - e_{9} -
\frac{3}{2}f - \op{pt}.
\]
In terms of $e_{8}$ this reads 
\[
\begin{split}
2\fm_{B}(\ell) - \fm_{B}(e_{8}) & = 1 + 2\ell -\sum_{i = 1}^{3}e_{i} 
- e_{7} - e_{8} - e_{9} -
\frac{3}{2}f - \op{pt} + \fm_{B}(\sum_{i = 1}^{3} e_{i} + e_{7}) \\
& = 1 + 2\ell  - \sum_{i=1}^{3}e_{i} - e_{7} - e_{8} - e_{9} -
\frac{3}{2}f - \op{pt} + \\
& \qquad \qquad\qquad + \left(4 + \sum_{i = 1}^{3} e_{i} + e_{7} 
- 4e_{9} - 6f - 4\op{pt}\right) \\
& = 5 + (2\ell - \frac{15}{2}f  - e_{8} -
5e_{9})  -5\op{pt}.
\end{split}
\]
Also from $\fm_{B}(f) = - \op{pt}$ we get 
\[
3\fm_{B}(\ell) - \fm_{B}(e_{8}) = 8 + (3\ell -12 f - e_{8} - 8e_{9}) -
8 \op{pt}. 
\]
Solving these two equations for $\fm_{B}(e_{8})$ results in 
\[
\fm_{B}(e_{8}) = 1 + (e_{8} - e_{9} - \frac{3}{2} f) - \op{pt},
\]
which completes the calculation of $\fm_{B}$. 

In summary, the action of $\ct$ and
the auxiliary actions of $\fm_{B}$ and $\fm_{B}^{-1}$ are recorded in tables 
\ref{table-t} and \ref{table-fm} respectively.

\begin{table}[!ht]
\begin{center}
\begin{tabular}{|l||c|c|} \hline
& $\fm_{B}$ & $\fm_{B}^{-1}$  \\ \hline\hline
$1$ & $- e_{9} + \frac{1}{2} \op{pt}$ & $e_{9} +  \frac{1}{2} \op{pt}$ \\
\hline
$\op{pt}$ & $f$ & $- f$ \\ \hline
$f$ & $-\op{pt}$ & $\op{pt}$ \\ \hline
$e_{i}$, & $1 + e_{i} - e_{9} - \frac{3}{2}f - \op{pt}$ & $-1 + e_{i}
- e_{9} - \frac{3}{2}f + \op{pt}$ \\ 
$i \neq 9$ & &  \\ \hline
$e_{9}$ & $1 - \frac{1}{2}f $ & $ -1 - \frac{1}{2}f$ \\ 
\hline
\end{tabular}
\end{center} 
\caption{Action of the cohomological Fourier-Mukai transform}
\label{table-fm}
\end{table}

\begin{table}[!ht]
\begin{center}
\begin{tabular}{|l||c|} \hline
 & $\ct_{B}$ \\ \hline \hline
$1$ & $1$ \\ \hline
$\op{pt}$ & $\op{pt}$ \\ \hline
$f$ & $f$ \\ \hline
$e_{j}$ & $2f + 2e_{9} - e_{j}  - 2\op{pt}$ \\
$j = 1,2,3$ & \\ \hline
$e_{i}$, & $2f - \ell + 2e_{9} + e_{7} + e_{i} - \op{pt}$ \\ 
$i = 4,5,6$ & \\ \hline
$e_{7}$ & $f + \ell - e_{1} - e_{2} - e_{3} + e_{9} - \op{pt}$ \\
\hline
$e_{8}$ & $2f - \ell + 2e_{9} + e_{7} + e_{8} - \op{pt}$ \\ \hline
$e_{9}$ & $e_{9}$ \\ \hline\hline
$\ell$ & $5f - e_{1} - e_{2} - e_{3} + e_{7} + 5e_{9} - 3 \op{pt}$ \\
\hline
\end{tabular} 
\end{center} 
\caption{Action of $\fm_{B}^{-1}\circ \tau_{B}^{*} \circ \fm_{B}$ on
cohomology}\label{table-t}
\end{table}

\section{Action on bundles} \label{ss-dp9-bundleaction}

In this section we show how the cohomological computations in the
previous section lift to actions of the Fourier-Mukai transform 
$\FM_{B}$ and the spectral
involution $\T_{B} := \FM_{B}^{-1}\circ
 \tau_{B}^{*} \circ \FM_{B}$ on (complexes of) sheaves on $B$.
Recall that the Chern character intertwines $\FM_{B}$ and $\fm_{B}$:
$\fm_{B} \circ ch = ch \circ \FM_{B}$. Similarly, it intertwines $\T_{B}$
and $\ct_{B}$:
$\ct_{B} \circ ch = ch \circ \T_{B}$.

Note that the Fourier-Mukai transform of a general sheaf ${\mathcal F}$
on $B$ is a complex of sheaves, not a single sheaf. Nevertheless,  all 
the sheaves we are interested in are taken by $\T_{B}$ again to sheaves.
To explain what is going on exactly we will need to introduce some
notation first. Put $c_{1} : D^{b}(B) \to \op{Pic}(B)$ for the first
Chern class map in Chow cohomology. In combination with $\T_{B}$, the
map $c_{1}$ induces a well defined map 
\begin{equation} \label{eq-T-on-Pic}
{\mathcal P}\!\op{ic}(B) \to \op{Coh}(B) \subset D^{b}(B)
\stackrel{\T_{B}}{\to}   D^{b}(B) \stackrel{c_{1}}{\to} \op{Pic}(B),
\end{equation}
where ${\mathcal P}\!\op{ic}(B)$ denotes the Picard category whose
objects are all line bundles on $B$ and whose morphisms are the
isomorphisms of line bundles. Since $\T_{B}$ is an autoequivalence,
the map \eqref{eq-T-on-Pic} descends to a well defined map of sets 
\[
\widetilde{\T}_{B} : \op{Pic}(B) = \pi_{0}({\mathcal P}\!\op{ic}(B)) \to
\op{Pic}(B). 
\]
If we identify $\op{Pic}(B)$ and $H^{2}(B,{\mathbb Z})$ via the first
Chern class map, we can describe $\widetilde{\T}_{B}$ alternatively as 
$\widetilde{\T}_{B}(-) = [\ct_{B}(\exp(c_{1}(-)))]_{2} \in
H^{2}(B,{\mathbb Z})$.

Denote by $\op{Pic}^{W}(B) \subset \op{Pic}(B)$ the subgroup generated
by $f$ and the classes of all sections of $\beta$ that meet the
neutral component of each fiber. A straightforward calculation shows that 
$\op{Pic}^{W}(B) = \op{Span}(f,e_{9},\{ f + e_{i} - e_{1} +
e_{9} \}_{i=2}^{6}, 2e_{7} - e_{8} + 2f)$ (note that $f + e_{i} - e_{1} +
e_{9}$ is the class of the section $[e_{i}] - [e_{1}]$ and $2e_{7} -
e_{8} + 2f$ is the class of the section $2[e_{7}]$) and that 
$\op{Span}(o_{1},o_{2})^{\perp} = \op{Span}(e_{9}, \{ e_{i} -
e_{1}\}_{i=2}^{6}, \ell - e_{7} - 2e_{1}, 2\ell - e_{8} - 4e_{1})$. In
particular $\op{Pic}^{W}(B)$ is a sublattice of index $3$ in
$\op{Span}(o_{1},o_{2})^{\perp}$. With this notation we have:

\begin{theo} \label{prop-T}
Let $L$ be a line bundle on $B$. Then 
\begin{itemize}
\item[{\em (i)}] The complex $\T_{B}(L) \in D^{[0,1]}(B)$ becomes a line 
bundle when restricted on the open set $B - (o_{1}\cup o_{2})$. More
precisely, the zeroth cohomology sheaf ${\mathcal H}^{0}(\T_{B}(L))$
is a line bundle on $B$ and the 
 first cohomology sheaf ${\mathcal H}^{1}(\T_{B}(L))$ is
supported on the divisor $o_{1} + o_{2}$.
\item[{\em (ii)}] The map $\widetilde{\T}_{B}$ satisfies
\[
 \widetilde{\T}_{B}(L)= 
\tau_B^*(L) \otimes{\mathcal O}_B((c_1(L) \cdot (e-\zeta))f + 
(c_1(L) \cdot f +1) (e-\zeta +f)).
\]
\item[{\em (iii)}] For every $L \in \op{Pic}^{W}(B)$ the image
$\T_{B}(L)$ is a line bundle on $B$ and so 
\[
\T_{B}(L) = \tau_B^*(L) \otimes{\mathcal O}_B((c_1(L) \cdot (e-\zeta))f + 
(c_1(L) \cdot f +1) (e-\zeta +f)).
\]
In particular $\T_{B} : \op{Pic}^{W}(B) \to (\op{Pic}^{W}(B) + (e-\zeta +f))
\subset \op{Pic}(B)$ is an affine isomorphism.
\end{itemize}
\end{theo}
{\bf Proof.} The proof of this proposition is rather technical and 
involves some elementary but long calculations in the derived category
$D^{b}(B)$. 

Since $\T_{B} = \FM_{B}^{-1}\circ \tau_{B}^{*} \circ \FM_{B}$
we need to understand $\FM_{B}^{-1}$. The following lemma is standard. 

\begin{lem} \label{lem-(FM-1)} The inverse $\FM_{B}^{-1}$ of the
Fourier-Mukai functor $\FM_{B}$ is isomorphic to the functor 
\[
\D_{B}\circ \FM_{B} \circ \D_{B} : D^{b}(B) \to D^{b}(B),
\]
where $\D_{B}$ is the (naive) Serre duality functor $\D_{B}(F) :=
\rhom(F,\omega_{B})$ with $\omega_{B}$ being the canonical line bundle on
$B$. 
\end{lem}
{\bf Proof.} It is well known (see e.g. \cite[Section~2]{orlov-k3})
that $\FM_{B}$  has left and right adjoint functors
$\FM_{B}^{*}$ and $\FM_{B}^{!}$ which are both isomorphic to 
$\FM_{B}^{-1}$. Furthermore, these adjoint functors can be defined by
explicit formulas, see \cite[Section~2]{orlov-k3}, e.g. the right
adjoint is given by:
\[
\FM_{B}^{!}(F) =
R\op{pr}_{1*}(\op{pr}_{2}^{*}F\stackrel{L}{\otimes} {\mathcal
P}^{\vee})){\otimes}\omega_{B}[2].
\] 
Here $\op{pr}_{i} : B\times B \to B$ are the projections onto the two
factors, ${\mathcal P} \to B\times B$ is the extension by zero of ${\mathcal
P}_{B}$  and  $K^{\vee} := \rhom(K,{\mathcal O}_{B\times B})$.
Using e.g. the formula for the right adjoint functor, the relative
duality formula \cite{hartshorne-rd} and the fact that $\omega_{B}$ is
a line bundle, one calculates 
\[
\begin{split}
\FM_{B}^{!}(F) & =
R\op{pr}_{1*}(\op{pr}_{2}^{*}F\stackrel{L}{\otimes} {\mathcal
P}^{\vee}){\otimes}\omega_{B}[2] \\
& = R\op{pr}_{1*}((\op{pr}_{2}^{*}F\stackrel{L}{\otimes} {\mathcal
P}^{\vee})\otimes \op{pr}_{2}^{*}\omega_{B}[2]\otimes
\op{pr}_{2}^{*}\omega_{B}^{-1} )\otimes \omega_{B} \\
&= R\op{pr}_{1*}(\op{pr}_{2}^{*}(F\otimes\omega_{B}^{-1})
\stackrel{L}{\otimes}{\mathcal P}^{\vee}\otimes
\op{pr}_{2}^{*}\omega_{B}[2]) \otimes \omega_{B} \\
& = (R\op{pr}_{1*}(\op{pr}_{2}^{*}(F^{\vee}\otimes\omega_{B})
\stackrel{L}{\otimes}{\mathcal P})^{\vee} \otimes \omega_{B} \\
& = (\FM_{B}(\D_{B}(F)))^{\vee}\otimes \omega_{B} \\
& = \D_{B}\circ\FM_{B}\circ \D_{B}(F).
\end{split}
\]
which proves the lemma.  \hfill $\Box$

\bigskip

Next observe that $\op{Pic}(B)$ is generated by 
all sections of $\beta$. Indeed $\op{Pic}(B)$ is generated by 
$\ell$ and $e_{1}, e_{2}, \ldots,
e_{9}$. The divisor classes  $e_{1}, \ldots, e_{7}$ and $e_{9}$ are already
sections of $\beta$. Also $\alpha_{B}(e_{1}) = \ell - e_{1} - e_{9}$
is a section and so $\ell$ is contained in the group generated by all
sections. Furthermore, $\alpha_{B}(e_{7}) = 2\ell - e_{1} - e_{2} -
e_{3} - e_{7} - e_{8}$ is a section and so $e_{8}$ is contained in
the group generated by all sections. 

In view of this it suffices to prove parts (i) and (ii) of the 
theorem for line bundles
of the form $L = {\mathcal O}_{B}(\sum a_{i}\xi_{i})$ where $a_{i} \in
{\mathbb Z}$ and $\xi_{i}$ are sections of $\beta$. 

Put $\cV_{0} := e_{*}{\mathcal O}_{\cp{1}}(-1)$. Consider the group
$\op{Ext}^{1}(\cV_{0},{\mathcal O}_{B})$ of extensions of 
$\cV_{0}$ by ${\mathcal O}_{B}$. 

Since $e^{2} = -1$ we have $\cV_{0} = e_{*}e^{*}{\mathcal O}_{B}(e)$
and so $\cV_{0}$ fits in a short exact sequence
\begin{equation} \label{eq-cV1}
0 \to {\mathcal O}_{B} \to {\mathcal O}_{B}(e) \to \cV_{0} \to 0.
\end{equation}
In particular we have a quasi-isomorphism $[{\mathcal O}_{B} \to
{\mathcal O}_{B}(e)] \widetilde{\to} \cV_{0}$ where in the complex 
\[
[{\mathcal O}_{B} \to  {\mathcal O}_{B}(e)],
\] 
the sheaf ${\mathcal O}_{B}$ is placed in degree $-1$ and ${\mathcal
O}_{B}(e)$ is placed in degree $0$. Thus we have
\[
\begin{split}
\op{Ext}^{1}(\cV_{0},{\mathcal O}_{B}) & =
\op{Hom}_{D^{b}(B)}(\cV_{0},{\mathcal O}_{B}[1]) = 
\op{Hom}_{D^{b}(B)}([{\mathcal O}_{B} \to {\mathcal O}_{B}(e)],
{\mathcal O}_{B}[1]) \\
& = {\mathbb H}^{0}(B, [{\mathcal O}_{B} \to
{\mathcal O}_{B}(e)]^{\vee}[1]) = {\mathbb H}^{0}(B, [{\mathcal
O}_{B}(-e)  \to {\mathcal O}_{B}]),
\end{split}
\]
where in the complex $[{\mathcal
O}_{B}(-e)  \to {\mathcal O}_{B}]$ the sheaf  ${\mathcal O}_{B}$ is
placed in degree zero. In particular we have a quasi-isomorphism 
$[{\mathcal O}_{B}(-e)  \to {\mathcal O}_{B}] \widetilde{\to}
e_{*}{\mathcal O}_{\cp{1}}$ and hence $\op{Ext}^{1}(\cV_{0},{\mathcal
O}_{B}) = H^{0}(B,e_{*}{\mathcal O}_{\cp{1}}) = {\mathbb C}$. This
shows that there is a unique (up to isomorphism) sheaf $\cV_{1}$ which
is a non-split extension of $\cV_{0}$ by ${\mathcal O}_{B}$.
But from \eqref{eq-cV1} we see that  the line
bundle ${\mathcal O}_{B}(e)$ is one such extension, i.e.
we must have $\cV_{1} \cong {\mathcal O}_{B}(e)$. 

Next consider the group of
extensions $\op{Ext}^{1}(\cV_{1},{\mathcal O}_{B}(f)) =
H^{1}(B,\cV_{1}^{\vee}\otimes {\mathcal O}(f))$. By the Leray spectral
sequence we have a short exact sequence
\[
0 \to H^{1}(\cp{1},(\beta_{*}\cV_{1}^{\vee})\otimes {\mathcal O}(1)) \to 
H^{1}(B,\cV_{1}^{\vee}\otimes {\mathcal O}(f)) \to H^{0}(\cp{1},
(R^{1}\beta_{*}\cV_{1}^{\vee})\otimes {\mathcal O}(1)) \to 0.
\]
But $\beta_{*}(\cV_{1}^{\vee}) =  \beta_{*}{\mathcal O}(-e) = 0$ and 
$R^{1}\beta_{*}(\cV_{1}^{\vee}) =  R^{1}\beta_{*}{\mathcal O}(-e) =
{\mathcal O}(-1)$. Thus $\op{Ext}^{1}(\cV_{1},{\mathcal O}_{B}(f)) =
H^{0}(\cp{1},{\mathcal O}) = {\mathbb C}$ and so there is a unique (up
to isomorphism) non-split extension 
\[
0 \to {\mathcal O}_{B}(f) \to \cV_{2} \to \cV_{1} \to 0.
\]
Arguing by induction we see that for every $a \geq 1$ there is a
unique up to isomorphism vector bundle $\cV_{a} \to B$ of rank $a$ on
$B$ satisfying $\beta_{*}(\cV_{a}^{\vee}) = 0$,
$R^{1}\beta_{*}(\cV_{a}^{\vee}) = {\mathcal O}(-a)$ and
$\op{Ext}^{1}(\cV_{a},{\mathcal O}_{B}(af)) = {\mathbb C}$ is generated
by the non-split short exact sequence
\[
0 \to {\mathcal O}_{B}(af) \to \cV_{a+1} \to \cV_{a} \to 0.
\]
Alternatively, for each positive integer we can consider the
vector bundle $\Psi_{a}$ of rank $a$ which is defined recursively as
follows:
\begin{itemize}
\item $\Psi_{1} := {\mathcal O}_{B}$, and
\item $\Psi_{a+1}$ is the unique non-split extension
\[
0 \to {\mathcal O}_{B}(af) \to \Psi_{a+1} \to \Psi_{a} \to 0.
\]
\end{itemize}

\bigskip
The fact that the $\Psi_{a}$'s are correctly defined can be checked
exactly as above. Moreover for each $a \geq 1$ 
$\cV_{a}$ can be identified with the unique
non-split extension 
\[
0 \to \Psi_{a} \to \cV_{a} \to e_{*}{\mathcal O}_{\cp{1}}(-1) \to 0.
\]
\

Let now $\xi : \cp{1} \to B$ be a section of $\beta$. The first step
in calculating $\T_{B}$ is given in the following lemma.

\begin{lem} \label{lem-FM-axi} For any integer $a$ we have
\[
\FM({\mathcal O}_{B}(a\xi)) = \begin{cases} \cV_{-a}\otimes {\mathcal
O}_{B}(\xi - e - (\xi\cdot e + 1)f)[-1], & \text{ for } a\leq 0 \\
 \cV_{a}^{\vee} \otimes {\mathcal
O}_{B}(-f)\otimes {\mathcal
O}_{B}(\xi - e - (\xi\cdot e + 1)f),  & \text{ for } a > 0
\end{cases}
\]
\end{lem}
{\bf Proof.} By Lemma~\ref{lem-fm-sheaves} we know that
$\FM_{B}({\mathcal O}_{B}) = e_{*}{\mathcal O}(-1)[-1]$ which gives
the statement of the lemma for $a= 0$. To prove the statement for $a =
-1 $ consider the short exact sequence 
\begin{equation} \label{eq-xi-ses}
0 \to {\mathcal O}_{B}(-\xi) \to {\mathcal O}_{B} \to \xi_{*}{\mathcal
O}_{\cp{1}} \to 0
\end{equation}
of sheaves on $B$. For an object $K \in D^{b}(B)$ let $\FM_{B}^{i}(K)$
denote the $i$-th cohomology sheaf of the complex $\FM_{B}(K)$. Since
$\FM_{B}$ is an exact functor on $D^{b}(B)$ it sends any short exact
sequence to a long exact sequence of cohomology sheaves. Applying 
$\FM_{B}$ to \eqref{eq-xi-ses}  and using
Lemma~\ref{lem-fm-sheaves} we get
\[
\xymatrix@R=6pt{
0 \ar[r] &  \FM_{B}^{0}({\mathcal O}_{B}(-\xi)) \ar[r] & 0 \ar[r] & 
{\mathcal O}_{B}(\xi - e - (1+\xi\cdot e)f)
\ar`r_l/5pt[lll] `^dr/5pt[lll][dll] & \\
& \FM_{B}^{1}({\mathcal O}_{B}(-\xi)) \ar[r] & e_{*}{\mathcal
O}(-1) \ar[r] & 0. &
}
\]
Thus $\FM_{B}^{0}({\mathcal O}_{B}(-\xi)) = 0$ and
$\FM_{B}^{1}({\mathcal O}_{B}(-\xi))$ fits in a short exact sequence 
\begin{equation} \label{eq-FM1(-xi)-ses}
0 \to {\mathcal O}_{B}(-e) \to \FM_{B}^{1}({\mathcal
O}_{B}(-\xi))\otimes {\mathcal O}(- \xi + (1 + \xi\cdot e)f) \to
e_{*}{\mathcal O}_{\cp{1}} \to 0. 
\end{equation}
Since \eqref{eq-xi-ses} is non-split and $\FM_{B}$ is an additive
functor, it follows that \eqref{eq-FM1(-xi)-ses} will not split. But
$\op{Ext}^{1}(e_{*}{\mathcal O}_{\cp{1}},{\mathcal O}_{B}(-e)) = 
\op{Ext}^{1}(e_{*}{\mathcal O}_{\cp{1}}(e),{\mathcal O}_{B}) =
{\mathbb C}$ as we saw above and therefore we must have 
\[
\FM_{B}({\mathcal O}_{B}(-\xi)) \cong {\mathcal O}_{B}(\xi - (1 +
\xi\cdot e)f)[-1] = \cV_{1}\otimes {\mathcal O}_{B}(\xi - e - (1 +
\xi\cdot e)f)[-1].
\]
Assume that the Lemma is proven for ${\mathcal O}_{B}(-a\xi)$ for
some positive $a$. Then we have a short exact sequence of sheaves on
$B$
\begin{equation} \label{eq-(-axi)-ses}
0 \to {\mathcal O}_{B}(-(a+1)\xi) \to {\mathcal O}_{B}(-a\xi) \to
\xi_{*}{\mathcal O}_{\cp{1}}(a) \to 0.
\end{equation}
Applying $\FM_{B}$ to \eqref{eq-(-axi)-ses} and using
Lemma~\ref{lem-fm-sheaves} we get
\[
\xymatrix@R=6pt{
0 \ar[r] &  \FM_{B}^{0}({\mathcal O}(-(a+1)\xi)) \ar[r] & 0 \ar[r] & 
{\mathcal O}(\xi - e  
+ (a -  1 -\xi\cdot e)f)
\ar`r_l/5pt[lll]`^dr/5pt[lll][dll] & \\
& \FM_{B}^{1}({\mathcal O}(-(a+1)\xi)) \ar[r] &
\FM_{B}^{1}({\mathcal O}(-a\xi)) \ar[r] & 0. &
}
\]
and so again $\FM_{B}^{0}({\mathcal O}_{B}(-(a+1)\xi)) =
0$. Furthermore, by the inductive hypothesis we have 
$\FM_{B}^{1}({\mathcal O}(-a\xi)) = \cV_{a}\otimes {\mathcal
O}_{B}(\xi - e - (1 + \xi\cdot e)f)$ and so by
the same reasoning as above the short exact sequence
\[
0 \to {\mathcal O}_{B}(af) \to \FM_{B}^{1}({\mathcal
O}_{B}(-(a+1)\xi))\otimes {\mathcal O}(e - \xi + (1 + \xi\cdot e)f) \to
\cV_{a} \to 0
\]
must be non-split. Since $\cV_{a+1}$ is the only such non-split
extension, we must have 
\[
\FM_{B}({\mathcal O}_{B}(-(a+1)\xi)) = \cV_{a+1}\otimes {\mathcal
O}_{B}(\xi - e - (1 + \xi\cdot e)f)[-1].
\]
This completes the proof of the lemma for all $a \leq 0$. 
The argument for $a > 0$ is exactly the same and is left as an
exercise. \hfill $\Box$

\bigskip

The next step is to calculate the action of $\T_{B}$ on line bundles
of the form ${\mathcal O}_{B}(a\xi)$. 

Due to Lemma~\ref{lem-(FM-1)} we have $\T_{B} = \D_{B}\circ \FM_{B}
\circ \D_{B} \circ \tau_{B}^{*} \circ \FM_{B}$.  Since $\tau_{B}$ is
an automorphism of $B$ 
we have $\D_{B}\circ \tau_{B}^{*} = \tau_{B}^{*}\circ \D_{B}$ and
so
\begin{equation} \label{eq-T}
\T_{B} = (\D_{B}\circ \FM_{B})\circ \tau_{B}^{*} \circ (\D_{B}\circ
\FM_{B}). 
\end{equation}
To calculate $\D_{B}(\FM_{B}({\mathcal O}_{B}(a\xi))$ we need to
distinguish two cases: $a = 0$ and $a \neq 0$. When $a= 0$, we have 
$\D_{B}((\FM_{B}({\mathcal O}_{B})) = \D_{B}(e_{*}{\mathcal
O}(-1)[-1])$. But as we saw above the short exact sequence
\eqref{eq-cV1} induces a quasi-isomorphism
\[
\left[ \begin{array}{c} {\mathcal O}_{B} \\ \downarrow \\ {\mathcal
O}_{B}(e) \end{array}\right] \!\!\!\! \begin{array}{c} 0 \\ \ \\ 1
\end{array}
\stackrel{\op{q.i.}}{\longrightarrow} e_{*}{\mathcal O}(-1)[-1].
\]
Applying duality one gets
\[
\D_{B}(e_{*}{\mathcal O}(-1)[-1])  = 
\left[ \begin{array}{c} {\mathcal O}_{B}(-e) \\ \downarrow \\ {\mathcal
O}_{B} \end{array}\right] \!\!\!\! \begin{array}{c} -1 \\ \ \\ 0
\end{array}\otimes {\mathcal O}_{B}(-f) 
 = \left[ \begin{array}{c} {\mathcal O}_{B}(-e - f) \\ \downarrow \\ {\mathcal
O}_{B}(-f) \end{array}\right] \!\!\!\! \begin{array}{c} -1 \\ \ \\ 0
\end{array} 
 = e_{*}{\mathcal O}_{\cp{1}}(-1).
\]
But for $a \neq 0$ the sheaves $\FM_{B}({\mathcal O}_{B}(a\xi))$ are
locally free and so we get 
\[
\D_{B}\circ \FM_{B}({\mathcal O}_{B}(a\xi)) = 
\begin{cases}
\cV_{-a}^{\vee}\otimes {\mathcal O}_{B}(e - \xi + (\xi\cdot e)f)[1],  &
\text{ for } a < 0 \\
\cV_{0},  & \text{ for } a = 0 \\
\cV_{a}\otimes {\mathcal O}_{B}(e - \xi + (1 + \xi\cdot e)f),  &
\text{ for } a > 0.
\end{cases}
\]
To apply $\tau_{B}^{*}$ next we need to calculate
$\tau_{B}^{*}\cV_{a}$. For this recall that $\cV_{a}$ is isomorphic to
the unique non-split extension 
\[ 
0 \to \Psi_{a} \to \cV_{a} \to e_{*}{\mathcal O}_{\cp{1}}(-1) \to 0.  
\]
Since $\tau_{B}(f) = f$ and $\Psi_{a}$ is built by successive
extensions of multiples of $f$, it follows that $\tau_{B}^{*}\Psi_{a}
\cong \Psi_{a}$ for every $a$. So $\cW_{a} := \tau_{B}^{*}\cV_{a}$ is
the unique non-split extension
\[
0 \to \Psi_{a} \to \cW_{a} \to \zeta_{*}{\mathcal O}_{\cp{1}}(-1) \to 0,
\]
where as before $\zeta = \tau_{B}^{*}(e)$. With this notation we have
\[
\tau_{B}^{*}\circ \D_{B}\circ \FM_{B}({\mathcal O}_{B}(a\xi)) = 
\begin{cases}
\cW_{-a}^{\vee}\otimes {\mathcal O}_{B}(\zeta - \tau_{B}^{*}(\xi) +
(\xi\cdot e)f)[1],   &
\text{ for } a < 0 \\
\zeta_{*}{\mathcal O}_{\cp{1}}(-1),  & \text{ for } a = 0 \\
\cW_{a}\otimes {\mathcal O}_{B}(\zeta - \tau_{B}^{*}(\xi) + 
(1 + \xi\cdot e)f),  &
\text{ for } a > 0.
\end{cases}
\]
Now to finish the calculation of $\T_{B}({\mathcal O}_{B}(a\xi))$ we
have to work out $\FM_{B}(\cW_{a}\otimes {\mathcal O}_{B}(\zeta -
\phi))$ and $\FM_{B}(\cW_{a}^{\vee} \otimes {\mathcal O}_{B}(\zeta -
\phi))$ for all $a > 0$ and all sections $\phi : \cp{1} \to B$ of
$\beta$. Again we proceed by induction in $a$. 

Let  $a = 1$. By definition $\cW_{1}$ is the unique non-split extension
\[
0 \to {\mathcal O}_{B} \to \cW_{1} \to \zeta_{*}{\mathcal O}(-1) \to 0,
\]
and hence $\cW_{1} = {\mathcal O}_{B}(\zeta)$ and $\cW_{1}^{\vee} 
= {\mathcal O}_{B}(-\zeta)$. In particular $\cW_{1}^{\vee}\otimes {\mathcal
O}_{B}(\zeta - \phi) = {\mathcal O}_{B}(-\phi)$. Consequently by
Lemma~\ref{lem-FM-axi} we get
\[
\FM_{B}(\cW_{1}^{\vee}\otimes {\mathcal O}_{B}(\zeta - \phi)) = 
\cV_{1}\otimes {\mathcal O}_{B}(\phi - e - (1 + \phi\cdot e)f)[-1] = 
{\mathcal O}_{B}(\phi - (1 + \phi\cdot e)f)[-1].
\]
Substituting $\phi = \tau_{B}^{*}(\xi)$ we get 
\[ 
\begin{split}
\FM\circ\tau_{B}^{*}\circ \D_{B}\circ \FM_{B} ({\mathcal O}_{B}(-\xi))
& = {\mathcal O}_{B}(\tau_{B}^{*}(\xi) - (1 + \tau_{B}^{*}(\xi)\cdot e
- \xi\cdot e)f) \\
& = {\mathcal O}_{B}(\tau_{B}^{*}(\xi) - (1 + \xi\cdot \zeta
- \xi\cdot e)f).
\end{split}
\]
\

Let now  $a = 2$. We have a short exact sequence
\[ 
0 \to {\mathcal O}_{B}(f)  \to \cW_{2} \to {\mathcal O}_{B}(\zeta) \to 0
\]
and so 
\begin{equation} \label{eq-cW2vee}
0 \to {\mathcal O}_{B}(- \phi)  \to \cW_{2}^{\vee}\otimes {\mathcal
O}_{B}(\zeta - \phi) \to {\mathcal O}_{B}(\zeta - \phi - f) \to 0.
\end{equation}
In particular we need to calculate $\FM_{B}({\mathcal
O}_{B}(\zeta - \phi))$. For this note that since ${\mathcal
O}_{B}(\zeta - \phi)$ is a line bundle which has degree zero on the
fibers of $\beta$, the sheaf $\FM_{B}^{0}({\mathcal
O}_{B}(\zeta - \phi))$ will have to be torsion free and torsion at
the same time and so $\FM_{B}^{0}({\mathcal
O}_{B}(\zeta - \phi)) = 0$ (see the argument on p.~\pageref{torsion}).
Consequently if we 
apply $\FM_{B}$ to the exact sequence
\[
0 \to {\mathcal O}_{B}(\zeta - \phi) \to {\mathcal O}_{B}(\zeta) \to
\phi_{*}{\mathcal O}_{\cp{1}}(\zeta\cdot \phi) \to 0,
\]
we will get a short exact sequence of sheaves
\[
0 \to {\mathcal O}_{B}(\zeta - 2e  - 2f) \to {\mathcal O}_{B}(\phi - e
- ( 1 + \phi\cdot e  - \phi\cdot \zeta)f) \to \FM_{B}^{1}({\mathcal
O}_{B}(\zeta - \phi)) \to 0. 
\]
In other words 
$\FM_{B}^{1}({\mathcal O}_{B}(\zeta - \phi))\otimes {\mathcal O}_{B}
( e - \phi + (1 + \phi\cdot e - \phi\cdot \zeta)f) = {\mathcal
O}_{D}$, where $D$ is an effective divisor in the linear system 
$|{\mathcal O}_{B}(\phi - \zeta + e + (1 - \phi\cdot e + \phi\cdot
\zeta)f)|$. 

To understand this linear system better consider the 
section $\mu : \cp{1} \to B$ for which $[\mu] = [\phi] - [\zeta]$ in
$\MW(B,e)$. Then as in section~\ref{ss-dp9-cohoaction} we can
write
\[
{\mathcal O}_{B}(\phi - \zeta) = {\mathcal O}_{B}(\mu - e + a f + b
n_{1} + c n_{2}). 
\]  
Taking into account that $\mu\cdot n_{i} = 1 - \phi\cdot n_{i}$ and
that $\mu^{2} = -1$ we can solve for $a$, $b$ and $c$ to get
\[
a = -1 +\phi\cdot e - \phi\cdot\zeta + \phi\cdot n_{1} + \phi\cdot
n_{2}, \qquad b = - \phi\cdot n_{1}, \qquad c = - \phi\cdot n_{2},
\]
which yields
\[
\begin{split}
{\mathcal O}_{B}(\phi - \zeta + e + (1 - \phi\cdot e + \phi\cdot
\zeta)f) & = {\mathcal O}_{B}(\mu + (\phi\cdot n_{1})o_{1} + (\phi\cdot
n_{2})o_{2}) \\
& = {\mathcal O}_{B}(\mu + (\mu\cdot o_{1})o_{1} + (\mu\cdot o_{2})o_{2}).
\end{split}
\]
Therefore, the numerical section $\mu + (\phi\cdot n_{1})o_{1} +
(\phi\cdot n_{2})o_{2}$ is the only effective divisor in the linear
system $|{\mathcal O}_{B}(\phi - \zeta + e + (1 - \phi\cdot e + \phi\cdot
\zeta)f)|$ and so $D = \mu + (\phi\cdot n_{1})o_{1} + (\phi\cdot
n_{2})o_{2}$ as divisors. Note that the fact that $\phi$ is a section
implies that $\phi\cdot n_{i}$ is either zero or one, and so $D$ is
always reduced.

This implies 
$\FM_{B}({\mathcal O}_{B}(\zeta - \phi)) = i_{D*}{\mathcal O}_{D}\otimes 
{\mathcal O}_{B}(\phi - e - (1 + \phi\cdot e - \phi\cdot
\zeta)f)[-1]$, where $i_{D} : D \hookrightarrow B$ is the natural inclusion. 
Next note that by definition of $\FM_{B}$ we have 
$\FM_{B}(K\otimes \beta^{*}M) =  \FM_{B}(K)\otimes \beta^{*}M$ for any
locally free sheaf $M \to \cp{1}$. Thus 
\[
\FM_{B}({\mathcal O}_{B}(\zeta - \phi - f)) = i_{D*}{\mathcal O}_{D}\otimes 
{\mathcal O}_{B}(\phi - e - (2 + \phi\cdot e - \phi\cdot
\zeta)f)[-1].
\]
We are now ready to apply $\FM_{B}$ to \eqref{eq-cW2vee}. The result
is 
\[
\xymatrix@R=6pt{
0 \ar[r] &  0 \ar[r] & {\mathcal S}^{0} \ar[r] & 0
\ar`r_l/5pt[lll]`^dr/5pt[lll][dll] & \\
& {\mathcal O}_{B}(\phi - (1 + \phi\cdot e)f)  \ar[r] &
 {\mathcal S}^{1}  \ar[r] & i_{D*}i_{D}^{*} 
{\mathcal O}_{B}(\phi - e - (2 + \phi\cdot e - \phi\cdot
\zeta)f) \ar[r] & 0,
}
\]
where ${\mathcal S}^{i} := \FM_{B}^{i}(\cW_{2}^{\vee}\otimes 
{\mathcal O}_{B}(\zeta - \phi))$. 

Writing  $\cL := {\mathcal O}_{B}(-e - (1 - \phi\cdot
\zeta)f)$ and ${\mathcal F} := {\mathcal S}^{1}\otimes {\mathcal
O}_{B}(-\phi + (1+\phi\cdot e)f)$, we find a non-split short exact sequence
\begin{equation} \label{eq-S1-ses}
0 \to {\mathcal O}_{B} \to {\mathcal F} \to i_{D*}i_{D}^{*}{\mathcal L} \to 0.
\end{equation}
Next we analyze the space of such extensions. We want to calculate
\[
\op{Ext}^{1}(i_{D*}i_{D}^{*}\cL,{\mathcal O}_{B}) =
Hom_{D^{b}(B)}(i_{D*}i_{D}^{*}\cL,{\mathcal O}_{B}[1]) = {\mathbb
H}^{0}(B, (i_{D*}i_{D}^{*}\cL)^{\vee}[1]).
\] 
As before, after tensoring the short exact sequence 
\[
0 \to {\mathcal O}_{B}(-D) \to {\mathcal O}_{B} \to i_{D*}{\mathcal
O}_{D} \to 0
\]
of the effective divisor $D$ by $\cL$ we get a quasi-isomorphism
\[
\left[ \begin{array}{c} \cL(-D) \\ \downarrow \\ 
\cL \end{array}\right] \!\!\!\! \begin{array}{c} -1 \\ \ \\ 0
\end{array}
\stackrel{\op{q.i.}}{\longrightarrow} i_{D*}i_{D}^{*}\cL,
\]
and so 
\[
(i_{D*}i_{D}^{*}\cL)^{\vee}[1] = \left[ \begin{array}{c} \cL^{\vee} 
\\ \downarrow \\ 
\cL^{\vee}(D) \end{array}\right] \!\!\!\! \begin{array}{c} -1 \\ \ \\ 0
\end{array} = i_{D*}i_{D}^{*}(\cL^{\vee}(D)).
\]
In particular $\op{Ext}^{1}(i_{D*}i_{D}^{*}\cL,{\mathcal O}_{B}) =
H^{0}(B, i_{D*}i_{D}^{*}(\cL^{\vee}(D))) =
H^{0}(D,i_{D}^{*}(\cL^{\vee}(D)))$. Since $D$ is a tree of smooth
rational curves, the dimension of the space of 
global sections of the line bundle $i_{D}^{*}(\cL^{\vee}(D))$ will
depend only on the degree of $\cL^{\vee}(D)$ on each component of
$D$. But $D = \mu + (\phi\cdot n_{1})o_{1} + (\phi\cdot n_{2})o_{2} = 
\mu + (\mu\cdot o_{1})o_{1} + (\mu\cdot o_{2})o_{2}$ and since $\mu$
is a section of $\beta$ we know that $\mu\cdot o_{i}$ is either $0$ or
$1$. We can distinguish three cases:
\begin{itemize}
\item[(a)] $\mu\cdot o_{1} = \mu\cdot o_{2} = 0$, i.e. $\mu \in
\op{Pic}^{W}(B)$ and $D = \mu$;
\item[(b)] $\mu$ intersects only one of the $o_{i}$'s, i.e. $D$ is the
union of $\mu$ and that $o_{i}$;
\item[(c)] $\mu\cdot o_{1} = \mu\cdot o_{2} = 1$ and so $D = \mu +
o_{1} + o_{2}$. 
\end{itemize}
Also since $D$ is linearly equivalent to $\phi - \zeta + e + ( 1 -
\phi\cdot e + \phi\cdot \zeta)f$ we find 
\[
\cL\cdot \mu = -1, \qquad \cL\cdot o_{1} = \cL\cdot o_{2} = 0.
\]
This gives the following answers for
$\op{Ext}^{1}(i_{D*}i_{D}^{*}\cL,{\mathcal O}_{B})$: 

\bigskip

\noindent
\underline{in case (a)}:  Since $D = \mu$ we have 
$(\cL^{\vee}(D))_{|D} =(\cL^{\vee}(\mu))_{|\mu}
= {\mathcal O}_{\mu}(1)\otimes {\mathcal O}_{\mu}(-1) = {\mathcal
O}_{\mu}$ and so $\op{Ext}^{1}(i_{D*}i_{D}^{*}\cL,{\mathcal O}_{B}) =
H^{0}(\mu,{\mathcal O}_{\mu}) = {\mathbb C}$.

\medskip

\noindent
\underline{in case (b)}: Say for concreteness $\mu\cdot o_{1} = 1$ and
$\mu\cdot o_{2} = 0$. Then $D = \mu + o_{1}$ is a normal crossing
divisor with a single singular point $\{x \} = \mu\cap o_{1}$. Then
$(\cL^{\vee}(D))_{|\mu} = {\mathcal O}_{\mu}(1)\otimes {\mathcal
O}_{\mu} = {\mathcal O}_{\mu}(1)$ and $(\cL^{\vee}(D))_{|o_{1}} = 
{\mathcal O}_{o_{1}}\otimes {\mathcal O}_{o_{1}}(-1) = {\mathcal
O}_{o_{1}}(-1)$. Hence $(\cL^{\vee}(D))_{|D}$ is the line bundle on
$D$ obtained by identifying the fiber $({\mathcal O}_{\mu}(1))_{x}$
with the fiber $({\mathcal O}_{o_{1}}(-1))_{x}$. Since
$H^{0}({\mathcal O}_{o_{1}}(-1)) = 0$ it follows that
$\op{Ext}^{1}(i_{D*}i_{D}^{*}\cL,{\mathcal O}_{B}) =
H^{0}(D,(\cL^{\vee}(D))_{|D})$ can be identified with the space of all
sections of ${\mathcal O}_{\mu}(1)$ that vanish at $x \in \mu$,
i.e. we again have $\op{Ext}^{1}(i_{D*}i_{D}^{*}\cL,{\mathcal O}_{B})
= {\mathbb C}$.

\medskip

\noindent
\underline{in case (c)}: 
The divisor $D = \mu + o_{1} + o_{2}$ is again a
normal crossings divisor but has now two singular points $x_{1}$ and
$x_{2}$, where $\{ x_{i} \} = \mu\cap o_{i}$ for $i = 1, 2$. In this
case we have $(\cL^{\vee}(D))_{|\mu} = {\mathcal O}_{\mu}(2)$ and 
$(\cL^{\vee}(D))_{|o_{i}} = {\mathcal O}_{o_{i}}(-1)$. Hence
$\op{Ext}^{1}(i_{D*}i_{D}^{*}\cL,{\mathcal O}_{B})$ gets identified
with the space of all sections in ${\mathcal O}_{\mu}(2)$ vanishing at
the points $x_{1}$ and $x_{2}$ and is therefore one dimensional. 
\bigskip

In other words we always  have a unique (up to isomorphism) choice for  the
sheaf ${\mathcal F}$. In fact, it is not hard to identify
the middle term of the non-split extension \eqref{eq-S1-ses}. Indeed, 
let $o := D - \mu$ be the union of the vertical components of $D$. We
have a short exact sequence:
\[
0 \to {\mathcal O}_{o}(-\mu) \to H^{0}(o,{\mathcal O}_{o}(\mu))\otimes
{\mathcal O}_{o} \to {\mathcal O}_{o}(\mu) \to 0.
\]
When we pull it back via
\[
{\mathcal O}_{B}(\mu) \to {\mathcal O}_{o}(\mu)
\]
we get a non-split sequence
\[
0 \to {\mathcal O}_{o}(-\mu) \to {\mathcal F}' \to {\mathcal
O}_{B}(\mu) \to 0. 
\]
Since we have already seen that such an extension is unique, we
conclude that ${\mathcal F}' = {\mathcal F}$.

We have shown that $\FM_{B}(\cW_{2}^{\vee}\otimes {\mathcal
O}_{B}(\zeta - \phi))[1]$ is 
a rank one sheaf on $B$ such that:
\begin{itemize}
\item The torsion in $\FM_{B}(\cW_{2}^{\vee}\otimes {\mathcal
O}_{B}(\zeta - \phi))[1]$ is ${\mathcal O}_{o}(-\mu)$.
\item $\FM_{B}(\cW_{2}^{\vee}\otimes {\mathcal O}_{B}(\zeta -
\phi))/(\op{torsion}) 
={\mathcal O}_{B}(2\phi - \zeta + e + (\phi\cdot \zeta - 2\phi\cdot
e)f - o)$.
\item The sheaf $\FM_{B}(\cW_{2}^{\vee}\otimes {\mathcal
O}_{B}(\zeta - \phi))[1]$ is the unique non-split extension of the
line bundle  ${\mathcal O}_{B}(2\phi - \zeta + e + (\phi\cdot \zeta -
2\phi\cdot  e)f - o)$  by the torsion sheaf ${\mathcal
O}_{o}(-\mu)$.
\end{itemize}

\bigskip

Let $a=3$. Then the short exact sequence 
\[
0 \to {\mathcal O}_{B}(2f) \to \cW_{3} \to \cW_{2} \to 0
\]
induces a short exact sequence
\[
0 \to \cW_{2}^{\vee}\otimes {\mathcal O}_{B}(\zeta - \phi) \to
\cW_{3}^{\vee}\otimes {\mathcal O}_{B}(\zeta - \phi) \to {\mathcal
O}_{B}(\zeta - \phi - 2f) \to 0.
\]
Applying $\FM_{B}$ one gets again that $\FM_{B}^{0}(\cW_{3}^{\vee}\otimes
{\mathcal O}_{B}(\zeta - \phi)) = 0$ and $\FM_{B}^{1}(\cW_{3}^{\vee}\otimes
{\mathcal O}_{B}(\zeta - \phi))$ fits in the non-split short exact sequence
\begin{equation} \label{eq-a=3}
0 \to \FM_{B}^{1}(\cW_{2}^{\vee}(\zeta - \phi)) \to
\FM_{B}^{1}(\cW_{3}^{\vee}(\zeta - \phi)) \to \FM_{B}^{1}({\mathcal
O}(\zeta 
- \phi))\otimes {\mathcal O}(-2f) \to 0. 
\end{equation}
\

\bigskip
Now recall that 
\[
\FM_{B}({\mathcal O}_{B}(\zeta - \phi)) = i_{D*}i_{D}^{*}{\mathcal
O}_{B}(\phi - e - (1 + \phi\cdot e - \phi\cdot\zeta)f),
\]
where $D = \mu + (\mu\cdot o_{1})o_{1} + (\mu\cdot o_{2})o_{2}$ is the
unique effective divisor in the linear system $|{\mathcal O}_{B}(\phi
- \zeta + e + (1 - \phi\cdot e + \phi\cdot \zeta)f)|$. 

In particular we have
\[
\begin{split}
\mu\cdot e & = \phi\cdot e - 1 + 1 - \phi\cdot e + \phi\cdot \zeta =
\phi\cdot \zeta \\
\mu\cdot \zeta & = \phi\cdot \zeta + 1 +  1 - \phi\cdot e + \phi\cdot
\zeta - \mu\cdot o_{1} - \mu\cdot o_{2} = 2 - \phi\cdot e + 2
\phi\cdot\zeta  - \mu\cdot o_{1} - \mu\cdot o_{2} \\
\mu\cdot \phi & = - 1 - \phi\cdot \zeta +  \phi\cdot e +  1 -
\phi\cdot e + \phi\cdot \zeta - (\mu\cdot o_{1})(\phi\cdot o_{1}) -
(\mu\cdot o_{2})(\phi\cdot o_{2}) = 0,
\end{split}
\]
and hence
\[
i_{\mu}^{*}{\mathcal
O}_{B}(\phi - e - (1 + \phi\cdot e - \phi\cdot\zeta)f) = {\mathcal
O}_{\mu}(-1 - \phi\cdot e).
\]
We are now ready to calculate $\FM_{B}({\mathcal O}_{B}(\zeta -
\phi))\otimes {\mathcal O}(-2f)$ for the three possible shapes of the
divisor $D$.

\bigskip

\noindent
\underline{Case (a)} \quad $D = \mu$ and so 
$\FM_{B}({\mathcal O}_{B}(\zeta - \phi))\otimes {\mathcal O}(-2f) =
{\mathcal O}_{\mu}(-3 - \phi\cdot e)$. Furthermore we showed that in
this case we have $\FM_{B}^{1}(\cW_{2}^{\vee}(\zeta - \phi)) = {\mathcal
O}_{B}(2\phi - \zeta + e + (\phi\cdot\zeta - 2 \phi\cdot e)f)$ and so
after twisting \eqref{eq-a=3} by ${\mathcal
O}_{B}(2\phi - \zeta + e + (\phi\cdot\zeta - 2 \phi\cdot e)f)^{-1}$ we
get a non-split short exact sequence
\[
0 \to {\mathcal O}_{B} \to ? \to {\mathcal O}_{\mu}(a) \to 0,
\]
where 
\[
? = \FM_{B}^{1}(\cW_{3}^{\vee}(\zeta - \phi))\otimes {\mathcal
O}_{B}(2\phi - \zeta + e + (\phi\cdot\zeta - 2 \phi\cdot e)f)^{-1},
\]
and
\[
a = - 3 - \phi\cdot e + \mu\cdot (-2\phi + \zeta - e - (\phi\cdot\zeta
- 2 \phi\cdot e)f) = - 1.
\]
Therefore we must have $? = {\mathcal O}_{B}(\mu) = {\mathcal
O}_{B}(\phi - \zeta + e + (1 - \phi\cdot e + \phi\cdot \zeta)f)$ and
so
\[
\FM_{B}^{1}(\cW_{3}^{\vee}(\zeta - \phi)) = {\mathcal O}_{B}(3\phi -
2\zeta + 2e + (1 - 3\phi\cdot e + 2\phi\cdot \zeta)f).
\]
\

\bigskip

\noindent
\underline{Case (b)} \quad In this case $\mu$ intersects exactly one
of the $o_{i}$, say $o_{1}$. Then $D = \mu + o_{1}$ and so
$\FM_{B}^{1}({\mathcal O}_{B}(\zeta - \phi))\otimes {\mathcal
O}_{B}(-2f) = {\mathcal O}_{\mu}(- 3 - \phi\cdot e)\cup_{x}
{\mathcal O}_{o_{1}}$ Moreover the torsion in
$\FM_{B}^{1}(\cW_{2}^{\vee}(\zeta - \phi))$ is ${\mathcal
O}_{o_{1}}(-1)$ and $\FM_{B}^{1}(\cW_{2}^{\vee}(\zeta -
\phi))/(\op{torsion}) = {\mathcal
O}_{B}(2\phi - \zeta + e + (\phi\cdot\zeta - 2 \phi\cdot e)f -
o_{1})$. Tensoring \eqref{eq-a=3} with ${\mathcal O}_{o_{1}}$ and
taking into account the fact that $\FM_{B}^{1}(\cW_{2}^{\vee}(\zeta -
\phi))_{|o_{1}} = {\mathbb C}^{2}\otimes {\mathcal O}_{o_{1}}$ we get
a long exact sequence of ${\mathcal T}\! or$ sheaves
\[
\xymatrix@R=7pt{
&   & \ldots \ar[r] 
& {\mathcal T}\!or_{1}^{{\mathcal O}_{B}}({\mathcal O}_{\mu}(-3 -
\phi\cdot e)\cup_{x} {\mathcal O}_{o_{1}},{\mathcal O}_{o_{1}}) 
\ar`r_l/5pt[lll]`^dr/5pt[lll][dll] & \\
&  {\mathbb C}^{2}\otimes {\mathcal O}_{o_{1}} \ar[r] &
\FM_{B}^{1}(\cW_{3}^{\vee}(\zeta -
\phi))_{|o_{1}} \ar[r] & {\mathcal O}_{o_{1}} \ar[r] & 0.
}
\]
Next we calculate 
 ${\mathcal T}\!or_{1}^{{\mathcal O}_{B}}({\mathcal O}_{\mu}(-3 -
\phi\cdot e)\cup_{x} {\mathcal O}_{o_{1}},{\mathcal
 O}_{o_{1}})$. 

\begin{lem} \label{lem-Tor}
${\mathcal T}\!or_{1}^{{\mathcal O}_{B}}({\mathcal O}_{\mu}(-3 -
\phi\cdot e)\cup_{x} {\mathcal O}_{o_{1}},{\mathcal
 O}_{o_{1}}) = 0$
\end{lem}
{\bf Proof.}  Recall that for any integer $a$ we have the following 
short exact sequence of sheaves on $B$:
\[
0 \to {\mathcal O}_{o_{1}}(-1) \to {\mathcal O}_{\mu}(a)\cup_{x}
{\mathcal O}_{o_{1}} \to {\mathcal O}_{\mu}(a) \to 0.
\]
Tensoring this sequence with ${\mathcal O}_{o_{1}}$ we obtain a long
exact sequence of ${\mathcal T}\! or$ sheaves:
\[
\xymatrix@R=7pt@C=8pt{
&  {\mathcal T}\!or_{1}^{{\mathcal O}_{B}}({\mathcal
O}_{o_{1}}(-1),{\mathcal O}_{o_{1}})  \ar[r]&  {\mathcal
T}\!or_{1}^{{\mathcal O}_{B}}({\mathcal O}_{\mu}(a)\cup_{x} {\mathcal
O}_{o_{1}},{\mathcal O}_{o_{1}}) \ar[r]  
&  {\mathcal
T}\!or_{1}^{{\mathcal O}_{B}}({\mathcal O}_{\mu}(a), 
{\mathcal O}_{o_{1}})
\ar`r_l/5pt[lll]`^dr/5pt[lll][dll] & \\
&  {\mathcal O}_{o_{1}}(-1) \ar[r] &
{\mathcal O}_{o_{1}} \ar[r] & {\mathcal O}_{x} \ar[r] & 0.
}
\]
In order to calculate the sheaves ${\mathcal T}\!or_{1}^{{\mathcal
O}_{B}}({\mathcal O}_{o_{1}}(-1),{\mathcal O}_{o_{1}})$ and 
${\mathcal  T}\!or_{1}^{{\mathcal O}_{B}}({\mathcal O}_{\mu}(a), 
{\mathcal O}_{o_{1}})$ recall that we have 
${\mathcal  T}\!or_{i}^{{\mathcal O}_{B}}(K,M) = {\mathcal
H}^{-i}(K\stackrel{L}{\otimes}_{{\mathcal O}_{B}} M)$ for any two
objects $K, M \in 
D^{b}(B)$. Now note that ${\mathcal O}_{o_{1}}(-1) = {\mathcal
O}_{o_{1}}\otimes {\mathcal O}_{B}(-\mu)$ and that 
\[
{\mathcal O}_{o_{1}} \stackrel{\text{q.i.}}{=} 
\left[ \begin{array}{c} {\mathcal O}_{B}(-o_{1}) \\ \downarrow \\ 
{\mathcal O}_{B} \end{array}\right] \!\!\!\! \begin{array}{c} -1 \\ \ \\ 0
\end{array}, \qquad 
{\mathcal O}_{\mu}(a) \stackrel{\text{q.i.}}{=}
\left[ \begin{array}{c} {\mathcal O}_{B}(af - \mu) \\ \downarrow \\ 
{\mathcal O}_{B}(af) \end{array}\right] \!\!\!\! \begin{array}{c} -1 \\ \ \\ 0
\end{array}, 
\]
and so 
\[
{\mathcal O}_{\mu}(a)\stackrel{L}{\otimes} {\mathcal O}_{o_{1}}
\stackrel{\text{q.i.}}{=} 
\left[ \begin{array}{c} {\mathcal O}_{B}(- \mu -o_{1}) \\ \downarrow \\ 
{\mathcal O}_{B}(-\mu)\oplus 
{\mathcal O}_{B}(-o_{1})  \\ \downarrow \\ 
{\mathcal O}_{B} \end{array}\right] \!\!\!\! \begin{array}{c}
-2 \\ \ \\ -1 \\ \ \\ 0 
\end{array} \otimes {\mathcal O}_{B}(af).
\]
Similarly 
\[
{\mathcal O}_{o_{1}}(-1)\stackrel{L}{\otimes} {\mathcal O}_{o_{1}}
\stackrel{\text{q.i.}}{=}
\left[ \begin{array}{c} {\mathcal O}_{B}(- 2o_{1}) \\ \downarrow \\ 
{\mathcal O}_{B}(-o_{1})\oplus 
{\mathcal O}_{B}(-o_{1})  \\ \downarrow \\ 
{\mathcal O}_{B} \end{array}\right] \!\!\!\! \begin{array}{c}
-2 \\ \ \\ -1 \\ \ \\ 0 
\end{array} \otimes {\mathcal O}_{B}(-\mu).
\]
Consequently ${\mathcal T}\!or_{i}^{{\mathcal
O}_{B}}({\mathcal O}_{o_{1}}(-1),{\mathcal O}_{o_{1}}) = 
{\mathcal  T}\!or_{i}^{{\mathcal O}_{B}}({\mathcal O}_{\mu}(a), 
{\mathcal O}_{o_{1}}) = 0$ for all $i \neq 0$. This proves the
lemma. \hfill $\Box$ 

\bigskip

\noindent
The previous lemma 
implies that $\FM_{B}^{1}(\cW_{3}^{\vee}(\zeta - \phi))_{|o_{1}} =
{\mathbb C}^{3}\otimes {\mathcal O}_{o_{1}}$ and that
$\FM_{B}^{1}(\cW_{3}^{\vee}(\zeta - \phi))$ fits in the commutative
diagram 
\[ 
\xymatrix@R=7pt@C=7pt{ & 0 \ar[d] & 0 \ar[d] & 0 \ar[d] & \\
0 \ar[r] & {\mathcal O}_{B}(2\phi - \zeta + e + (\phi\zeta - 2\phi e)
f - 2o_{1}) \ar[r] \ar[d] & ? \ar[r] \ar[d] & {\mathcal O}_{\mu}(-4 -
\phi e) \ar[r] \ar[d] & 0 \\
0 \ar[r] & \FM_{B}^{1}(\cW_{2}^{\vee}(\zeta - \phi)) \ar[r] \ar[d] & 
\FM_{B}^{1}(\cW_{3}^{\vee}(\zeta - \phi)) \ar[r] \ar[d] & 
{\mathcal O}_{\mu}(-3 -\phi e)\cup_{x} {\mathcal O}_{o_{1}} \ar[r]
\ar[d] & 0 \\
0 \ar[r] & {\mathbb C}^{2}\otimes {\mathcal O}_{o_{1}} \ar[r] \ar[d] & 
{\mathbb C}^{3}\otimes {\mathcal O}_{o_{1}} \ar[r] \ar[d] & 
{\mathcal O}_{o_{1}}\ar[r] \ar[d] & 0 \\
& 0 & 0 & 0 & 
}
\]
where ? is a non-split extension of ${\mathcal O}_{\mu}(-4 -
\phi e)$ by ${\mathcal O}_{B}(2\phi - \zeta + e + (\phi\zeta - 2\phi e)
f - 2o_{1})$. This implies that $? = {\mathcal O}_{B}(2\phi - \zeta +
e + (\phi\zeta - 2\phi e)f - 3o_{1})$ and that
$\FM_{B}^{1}(\cW_{3}^{\vee}(\zeta - \phi))$ fits in a short exact
sequence
\[
0 \to {\mathcal O}_{B}(3\phi - 2\zeta +
2e + (1 + 2\phi\zeta - 3\phi e)f - 3o_{1}) \to
\FM_{B}^{1}(\cW_{3}^{\vee}(\zeta - \phi)) \to {\mathbb C}^{3}\otimes
{\mathcal O}_{o_{1}} \to 0.
\]
In particular we see that the torsion in
$\FM_{B}^{1}(\cW_{3}^{\vee}(\zeta - \phi))$  is supported on $o_{1}$.

The same reasoning applied to the restriction of \eqref{eq-a=3} to
$\mu$ instead of  $o_{1}$ implies that $\FM_{B}^{1}(\cW_{3}^{\vee}(\zeta
- \phi))/(\op{torsion})$ is isomorphic to
 the line bundle ${\mathcal O}_{B}(3\phi - 2\zeta +
2e + (1 + 2\phi\zeta - 3\phi e)f - 2o_{1})$. Since ${\mathcal
O}_{B}(3\phi - 2\zeta +  2e + (1 + 2\phi\zeta - 3\phi e)f -
2o_{1})_{|o_{1}} = {\mathcal O}_{o_{1}}(2)$ we conclude that the
torsion in $\FM_{B}^{1}(\cW_{3}^{\vee}(\zeta - \phi))$ is isomorphic to
the kernel of the natural map ${\mathbb C}^{3}\otimes {\mathcal
O}_{o_{1}} \cong H^{0}(o_{1},{\mathcal O}_{o_{1}}(2x))\otimes {\mathcal
O}_{o_{1}} \to {\mathcal O}_{o_{1}}(2x) \cong {\mathcal
O}_{o_{1}}(2)$. In particular we see that the torsion in
$\FM_{B}^{1}(\cW_{3}^{\vee}(\zeta - \phi))$ is a rank two vector bundle
on $o_{1}$, which has no sections and is of degree $-2$, i.e. is
isomorphic to ${\mathcal O}_{o_{1}}(-1)\oplus {\mathcal
O}_{o_{1}}(-1)$.

\bigskip

\noindent
\underline{Case (c)} \quad In this case $D = \mu + o_{1} + o_{2}$. 
An analysis, analogous to the one used in case (b), now shows that 
the torsion in 
$\FM_{B}^{1}(\cW_{3}^{\vee}(\zeta - \phi))$ is isomorphic to ${\mathcal
O}_{o_{1}}(-1)^{\oplus 2}\oplus {\mathcal O}_{o_{2}}(-1)^{\oplus 2}$
and that $\FM_{B}^{1}(\cW_{3}^{\vee}(\zeta - \phi))/(\op{torsion})$ is
isomorphic to the line bundle ${\mathcal O}_{B}(3\phi - 2\zeta +
2e + (1 + 2\phi\zeta - 3\phi e)f - 2o_{1} - 2o_{2})$.

Continuing inductively we get that for every $a \geq 1$ the object 
$\FM_{B}(\cW_{a}^{\vee}\otimes {\mathcal O}_{B}(\zeta - \phi))[1]$ is
a rank one sheaf on $B$ such that
\begin{itemize}
\item The torsion in $\FM_{B}(\cW_{a}^{\vee}\otimes
{\mathcal O}_{B}(\zeta - \phi))[1]$ is isomorphic to 
\[
{\mathcal
O}_{o_{1}}^{\oplus (a-1)(\phi\cdot n_{1})}(-1)\oplus
{\mathcal O}_{o_{2}}^{\oplus (a-1)(\phi\cdot n_{2})}(-1).
\] 
(In this
formula it is tacitly understood that the direct sum of
zero copies of a sheaf is the zero sheaf.)
\item The sheaf $\FM_{B}(\cW_{a}^{\vee}\otimes
{\mathcal O}_{B}(\zeta - \phi))[1]/(\op{torsion})$ is isomorphic to
the line bundle 
\[
\begin{split}
{\mathcal O}_{B}
( a\phi & - (a-1)\zeta + (a-1)e + ((a-2) + (a-1)\phi\cdot \zeta - a\phi\cdot
e)f \\
& - (a-1)(\phi\cdot n_{1})o_{1} - (a-1)(\phi\cdot n_{2})o_{2}).
\end{split}
\]
\end{itemize}
\

Now by substituting $\phi = \tau_{B}^{*}(\xi)$ in the above formula and by
noticing that $D({\mathcal O}_{o_{i}}(-1)) = {\mathcal O}_{o_{i}}(-1)[-1]$
we obtain
\[
\begin{split}
{\mathcal H}^{0}\T_{B}(& {\mathcal O}_{B}(-a\xi)) = \\
& = {\mathcal O}_{B}(\tau_{B}^{*}(-a\xi)
+ ((-a\xi)\cdot(e - \zeta))f  + (1-a)(e - \zeta + f) \\
& \qquad\qquad\qquad\qquad   + (a-1)(\xi\cdot
o_{1})o_{1} + (a-1)(\xi\cdot o_{2})o_{2}), \\
{\mathcal H}^{1}\T_{B}(& {\mathcal O}_{B}(-a\xi)) =
{\mathcal O}_{o_{1}}^{\oplus (a-1)(\xi\cdot o_{1})}(-1)
\oplus
{\mathcal O}_{o_{2}}^{\oplus (a-1)(\xi\cdot o_{2})}(-1),
\end{split}
\]
for all $a \geq 1$. We have already analyzed the case $a=0$ above and
so this proves the theorem for $L = {\mathcal
O}_{B}(-a\xi)$ and $a \geq 0$. The cases $L = {\mathcal
O}_{B}(a\xi)$ with $a >0$ or $L = {\mathcal
O}_{B}(\sum a_{i}\xi_{i})$ with different $\xi_{i}$'s are 
analyzed in exactly the same way. \hfill $\Box$

\bigskip
\bigskip

\begin{rem} \label{rem-T} \ {\bf (i)} The calculation of $\T_{B}(L)$ in the
proof of Theorem~\ref{prop-T} works equally well on 
a rational elliptic surface in the five dimensional
family from Corollary~\ref{cor-5parameters} (with the choice of
$\zeta$ as in Remark~\ref{rem-find-zeta}). Since in this case 
$\op{Pic}^{W}(B) = \op{Pic}(B)$, we see that for a general $B$ in the
five dimensional family we have $\T_{B|\op{Pic}(B)} =
\widetilde{T}_{B}$. In particular $\T_{B}$ sends all line bundles to
line bundles and induces an affine automorphism on $\op{Pic}(B)$.

\medskip

\noindent
{\bf (ii)} In the proof of Theorem~\ref{prop-T} we also showed
that the statement of Theorem~\ref{prop-T}(iii) admits a partial
inverse. Namely, we showed that if $L$ is a multiple of a section,
then $\T_{B}(L)$ is a line bundle if {\em and only if} $L \in
\op{Pic}^{W}(B)$. 
\end{rem}
\

\bigskip

The previous theorem shows that the  $\T_{B}$ action on $\op{Pic}(B)$ is 
somewhat complicated. If we work modulo the exceptional curves $o_1, o_2$, the
formulas simplify considerably. (Working modulo $o_1, o_2$ amounts to
contracting these two curves.)

\begin{cor} \label{cor-TmodO}
The action of $\widetilde{\T}_{B}$ induces an affine automorphism of
$\op{Pic}(B) /  
({\mathbb Z}o_1 \oplus {\mathbb Z}o_2)$, namely:
\[
\widetilde{\T}_{B}(L) = \alpha_B^*(L) \otimes{\mathcal O}_B(e-\zeta +f) \quad
\mod(o_1,o_2).
\]
\end{cor}
{\bf Proof.} Apply Theorem~\ref{prop-T} together with
\eqref{eq-i2-fibers}. \ \hfill $\Box$

\bigskip

Using these two results we can now describe  the action of $\T_{B}$ on
sheaves supported on curves in $B$. Let $C \subset B$ be a curve
which is finite over $\cp{1}$. Denote by
$i_{C} : C \hookrightarrow B$ the inclusion map. For the purposes of
the spectral construction we will need to calculate the action of the
spectral involution $\T_{B}$
on sheaves of the form $i_{C*}i_{C}^{*}L$ for some $L \in
\op{Pic}(B)$:

\begin{prop} \label{prop-Tonspectral} Let $C \subset B$ be a curve
which is finite over $\cp{1}$ and such that ${\mathcal O}_{B}(C) \in
\op{Pic}^{W}(B)$ (for example we may take $C$ in the linear system $|re +
kf|$ for some integers $r$, $k$). Let $L \in \op{Pic}(B)$. Put $D :=
\alpha_{B}(C)$. Then
\begin{itemize}
\item[{\em (a)}] $\T_{B}(i_{C*}i_{C}^{*}L) = i_{D*}i_{D}^{*}(\T_{B}(L))$.
\item[{\em (b)}] $\T_{B}(i_{C*}i_{C}^{*}L) =
i_{D*}i_{D}^{*}(\alpha_B^*(L) \otimes{\mathcal O}_B(e-\zeta +f))$.
\end{itemize}
\end{prop} 
{\bf Proof.}  Since $C$ is assumed to be finite over $\cp{1}$ it
follows that $i_{C}^{*}L$ will be flat over $\cp{1}$ and so $V =
\FM_{B}(L)$ will be a vector bundle on $B$ of rank $r = C\cdot f$, 
which is semistable and of
degree zero on every fiber of $\beta$. But then $\tau_{B}^{*}V$ will be
again a vector bundle of this type. Moreover if $f_{t}$ is a general fiber
of $\beta$ then we can write $V_{|f_{t}} \cong a_{1}\oplus \ldots
\oplus a_{r}$, where $a_{i}$ are line bundles of degree zero on
$f_{t}$. In fact if we put $\{p_{1}, \ldots, p_{r} \} = C\cap f_{t}$
for the intersection points of $C$ and $f_{t}$ we have $a_{i} 
= {\mathcal O}_{f_{t}}(p_{i} - e(t))$. Now $\tau_{B}$ induces an
isomorphism $\tau_{B} : f_{\tau_{\cp{1}}(t)} \to f_{t}$ and 
\[
(\tau_{B}^{*}V)_{|f_{\tau_{\cp{1}}(t)}} = \tau_{B}^{*}a_{1}\oplus
\ldots \oplus \tau_{B}^{*}a_{r}.
\]
By definition $\tau_{B} = t_{\zeta}\circ \alpha_{B}$. Since every
translation on an elliptic curve induces the identity on
$\op{Pic}^{0}$ it follows that $\tau_{B}^{*}a_{i} = \alpha_{B}^{*}a_{i} =
{\mathcal O}_{f_{\tau_{\cp{1}}(t)}}(\alpha_{B}(p_{i}) -
e(f_{\tau_{\cp{1}}(t)}))$. This shows that $\FM_{B}^{-1}(\tau_{B}^{*}V)$
will be a line bundle supported on $D = \alpha_{B}(C)$ and so to prove
(a) we only need to identify this line bundle explicitly. 

Consider the short exact sequence 
\[
0 \to L(-C) \to L \to i_{C*}i_{C}^{*}L \to 0.
\]
Applying the exact functor $\T_{B}$ we get a long exact sequence of
sheaves
\[
\xymatrix@R=6pt{
0 \ar[r] & {\mathcal H}^{0}\T_{B}(L(-C)) \ar[r] & {\mathcal
H}^{0}\T_{B}(L) \ar[r] & 
\T_{B}(i_{C*}i_{C}^{*}L)
\ar`r_l/5pt[lll]`^dr/5pt[lll][dll] & \\
&  {\mathcal H}^{1 }\T_{B}(L(-C))\ar[r] & {\mathcal H}^{1
}\T_{B}(L)\ar[r] 
& 0. &
}
\]
However, by parts (i) and (iii) of 
Theorem~\ref{prop-T}  we have 
\[
{\mathcal H}^{1 }\T_{B}(L(-C)) =  {\mathcal H}^{1}\T_{B}(L)
\]
and so $\T_{B}(i_{C*}i_{C}^{*}L)$ fits in a short exact sequence 
\[
0 \to {\mathcal H}^{0}\T_{B}(L(-C)) \to {\mathcal H}^{0}\T_{B}(L) \to
\T_{B}(i_{C*}i_{C}^{*}L) \to 0.
\]
But in the proof of Theorem~\ref{prop-T} we showed that
for any line bundle $K \in \op{Pic}(B)$ one has 
\[
{\mathcal H}^{0}\T_{B}(K) = \widetilde{T}_{B}(K)\otimes {\mathcal
O}_{B}((c_{1}(K)\cdot o_{1})o_{1} + (c_{1}(K)\cdot o_{2})o_{2}).
\]
Taking into account that ${\mathcal O}(o_{i})_{|C} = {\mathcal O}_{C}$
we can twist the above exact sequence by 
\[
{\mathcal
O}_{B}(-(c_{1}(K)\cdot o_{1})o_{1} - (c_{1}(K)\cdot o_{2})o_{2})
\] 
to
obtain
\[
0 \to \widetilde{\T}_{B}(L(-C)) \to \widetilde{\T}_{B}(L)\to
\T_{B}(i_{C*}i_{C}^{*}L) \to 0. 
\]
To calculate $\widetilde{\T}_{B}(L(-C))$ 
let $\Omega : \op{Pic}(B) \to \op{Pic}(B)$ denote the
linear part of the affine map $\widetilde{\T}_{B}$. In other words
$\Omega(L) = \tau_{B}^{*}(L) + (c_{1}(L)\cdot(e - \zeta))f +
(c_{1}(L)\cdot f)(e - \zeta + f)$ and $\widetilde{\T}_{B}(L) =
\omega(L) + (e - \zeta + f)$.  Then  $\widetilde{\T}_{B}(L(-C)) =
\widetilde{\T}_{B}(L)\otimes {\mathcal O}_{B}(-\Omega(C))$. 

Using the formula describing
$\Omega$ one checks immediately that $\Omega$ is a linear involution of 
$\op{Pic}(B)$ which preserves the intersection pairing. Also we have 
$\Omega(o_{1}) = - o_{2}$ and $\Omega(o_{2}) = - o_{2}$ and so
$\Omega$ preserves $\op{Span}(o_{1},o_{2})^{\perp}$. But according to
Corollary~\ref{cor-TmodO} the restriction of $\Omega$ to
$\op{Span}(o_{1},o_{2})^{\perp} \supset \op{Pic}^{W}(B)$ coincides
with the restriction of $\alpha_{B}^{*}$, which yields 
\[
\widetilde{\T}_{B}(L(-C)) =
\widetilde{\T}_{B}(L)\otimes {\mathcal O}_{B}(-\Omega(C)) =
\widetilde{\T}_{B}(L)\otimes {\mathcal O}_{B} (-\alpha_{B}^{*}(C)) =
\widetilde{\T}_{B}(L)\otimes {\mathcal O}_{B} (-D).
\]
Consequently $\T_{B}(i_{C*}i_{C}^{*}L)$ fits in the exact sequence
\[
0 \to \widetilde{\T}_{B}(L)\otimes {\mathcal O}_{B} (-D) \to
\widetilde{\T}_{B}(L) \to \T_{B}(i_{C*}i_{C}^{*}L) \to 0.
\]
But as we saw above $\T_{B}(i_{C*}i_{C}^{*}L)$ is the extension by
zero of some line bundle on $D$ and so we must have
$\T_{B}(i_{C*}i_{C}^{*}L) =
i_{D*}i_{D}^{*}\widetilde{\T}_{B}(L)$. Finally note that
$\alpha_{B}^{*}$ preserves $\op{Pic}^{W}(B)$ since $\alpha_{B}^{*}(o_{1})
= o_{2}$. Therefore $D$ is disjoint from $o_{1}$ and $o_{2}$ and so 
the restriction of $\widetilde{\T}_{B}(L)$ to $D$ will be the same as
the restriction of the projection of $\widetilde{\T}_{B}(L)$ onto 
$\op{Span}(o_{1},o_{2})^{\perp}$. Applying again
Corollary~\ref{cor-TmodO} we get that
$i_{D*}i_{D}^{*}\widetilde{\T}_{B}(L) = i_{D*}i_{D}^{*}{\mathcal
O}_{B}(\alpha_{B}^{*}L + (e - \zeta + f))$. The Proposition is
proven. \hfill $\Box$


\end{document}